\documentclass[12pt]{article}

\usepackage{amsmath}
\usepackage{xcolor}
\usepackage{algorithm}
\usepackage{algpseudocode}
\usepackage{overpic}
\usepackage{multirow}

\usepackage{empheq,mathtools}

\DeclareMathOperator*{\argmin}{argmin}

\DeclareMathOperator*{\diam}{diam}
\DeclareMathOperator*{\supp}{supp}
\DeclareMathOperator*{\card}{card}
\DeclareMathOperator*{\Tau}{\mathcal{T}}

\newcommand{\TauFHrom}{\mathcal{T}_{\textrm{f}}}      
\newcommand{\TauAHrom}{\mathcal{T}_{\textrm{a}}}

\newcommand{\spanOp}[1]{\big\{#1\big\}}

\newcommand{\norm}[2]{\big\lVert#1\big\rVert_{#2}}

\newcommand{\myDer}{\partial}
\newcommand{\der}[2]{\frac{\myDer#1}{\myDer#2}}
\newcommand{\derSecond}[3]{\frac{\myDer^#3#1}{\myDer#2^#3}}

\newcommand{\qedwhite}{\hfill \ensuremath{\Box}}

\usepackage{bbm}
\usepackage{varioref}
\usepackage{hyperref}
\usepackage{soul}

\usepackage[hyperref=true,%
    backend=bibtex,%
    bibstyle=abbrv,%
    maxbibnames=99,%
    style=numeric]{biblatex}
\usepackage{parskip}
\usepackage{amssymb}

\usepackage{bm}
\usepackage{subcaption}
\usepackage{adjustbox}

\usepackage{tikz}

 \usepackage{geometry}
\definecolor{myGreen}{HTML}{2E8B57}

\usepackage[amsmath]{ntheorem}
\usepackage[nameinlink]{cleveref}
\usepackage{xfrac}

\usepackage{siunitx}
\newtheorem{lemma}{Lemma}[section]
\newtheorem{definition}[lemma]{Definition}
\newtheorem{theorem}[lemma]{Theorem}

\Crefname{definition}{Definition}{definitions}
\Crefname{theorem}{Theorem}{theorems}
\Crefname{assumption}{Assumption}{assumptions}
\Crefname{remark}{Remark}{remarks}  

\usepackage[normalem]{ulem}
\newcommand{\stkout}[1]{\ifmmode\text{\sout{\ensuremath{#1}}}\else\sout{#1}\fi}

\newbool{commenti}
\newbool{slowCompilation}

\booltrue{commenti}

\booltrue{slowCompilation}

\usepackage{soul}
\newcommand{\IM}[1]{\ifbool{commenti}{\textcolor{red}{#1}}{}}
\definecolor{myPurple}{HTML}{800080}
\newcommand{\Ar}[1]{\ifbool{commenti}{\textcolor{blue}{#1}}{}}

\usepackage{MnSymbol,wasysym}
\usepackage{booktabs}
\usepackage{authblk}

\title{A hybrid finite volume - spectral element method  for aeroacoustic problems}

\author[1,*]{Alberto Artoni}
\author[1]{Paola F. Antonietti}
\author[1]{Ilario Mazzieri}
\author[1]{Nicola Parolini}
\author[2]{Daniele Rocchi}

\affil[1]{MOX-Laboratory for Modeling and Scientific Computing, Department of Mathematics, Politecnico
di Milano, 20133 Milan, Italy}
\affil[2]{Department of Mechanical Engineering, Politecnico di Milano, 20156 Milan, Italy}

\affil[*]{\small Corresponding author: \texttt{alberto.artoni@polimi.it}}

\begin{document}
	
	\maketitle
 
	\section*{Abstract} 
	We propose a hybrid Finite Volume (FV) - Spectral Element Method (SEM) for modelling aeroacoustic phenomena based on the Lighthill's acoustic analogy. First the fluid solution is computed employing a FV method. Then, the sound source term is projected onto the acoustic grid and the inhomogeneous Lighthill's wave equation is solved employing the SEM. The novel projection method computes offline the intersections between the acoustic and the fluid grids in order to preserve the accuracy. The proposed intersection algorithm is shown to be robust, scalable and able to efficiently compute the geometric intersection of arbitrary polyhedral elements. 
    We then analyse the properties of the projection error, showing that if the fluid grid is fine enough we are able to exploit the accuracy of the acoustic solver and we numerically assess the obtained theoretical estimates. 
    Finally, we address two relevant aeroacoustic benchmarks, namely the corotating vortex pair and the noise induced by a laminar flow around a squared cylinder, to demonstrate in practice the effectiveness of the projection method when dealing with high order solvers. The flow computations are performed with OpenFOAM \cite{OpenFOAM}, an open-source finite volume library, while the inhomogeneous Lighthill's wave equation is solved with SPEED \cite{SPEED2013}, an open-source spectral element library.

   \section{Introduction}
	
\begin{figure}[h]
		\begin{center} 
			\includegraphics[width=1\textwidth]
    {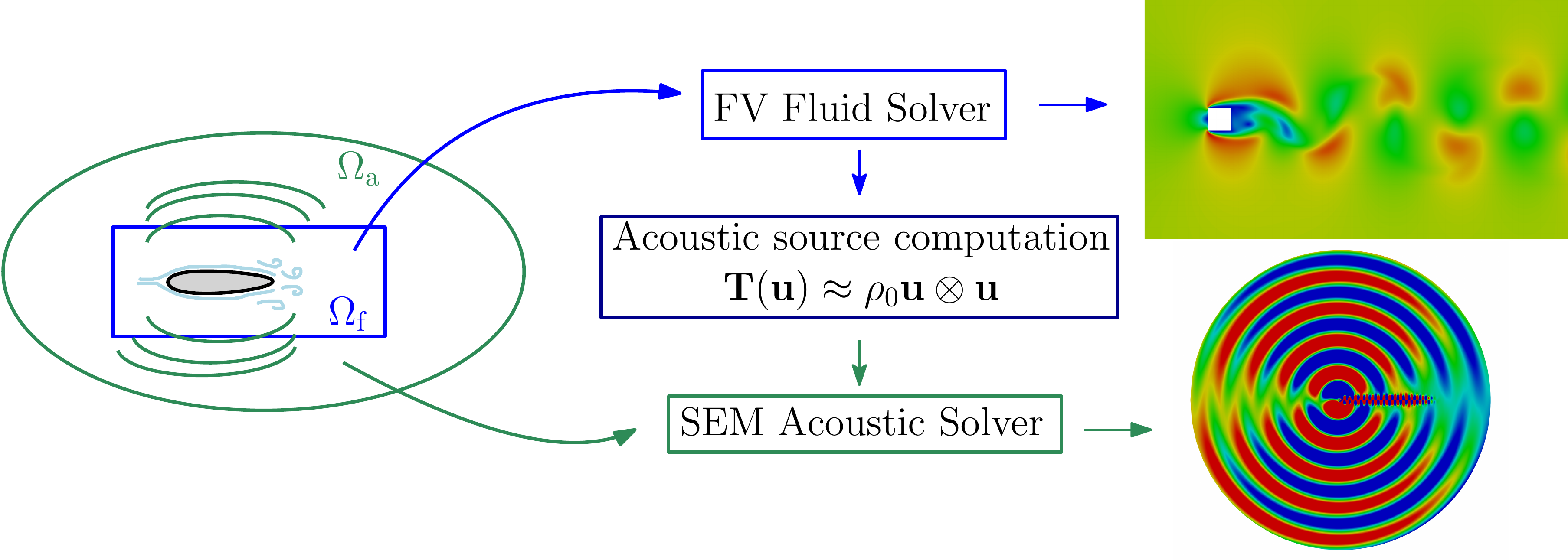}
			\caption{Computational domain for the aeroacoustic problem. First, the fluid problem is solved on $\Omega_{\textrm{f}}$. Then the acoustic source term is computed from the flow velocity. Finally, an inhomogeneous acoustic wave equation is solved on $\Omega_{\textrm{a}}$.} \label{fig:domainAeroProblem}
		\end{center}
	\end{figure}

	Aeroacoustics studies the propagation of noise generated by fluid flows. A typical problem of interest can be the noise induced by a car side view mirror. For a car moving at $\SI{40}{\meter\per\second}$, the corresponding Reynolds number is of the order of $10^6$. The typical mesh size required to capture the fluid length scales (even when working with turbulence models) are of the order of $10^{-3}\SI{}{\meter}$ or $10^{-4}\SI{}{\meter}$ \cite{Wang2017}, far from the involved acoustic scales that range from $\SI{0.05}{\meter}$ up to $\SI{5}{\meter}$, see for instance \cite{Munz2016}.  Due to the multiscale nature of the involved length scales, a widely employed class of Computational Aeroacoustics (CAA) methodologies separate the flow field from the acoustic computations, in a hybrid approach, \Ar{see for instance the recent reviews in \cite{Kaltenbacher2020} or \cite{Lele2004}}.
	Those methods are based on aeroacoustic analogies, namely rearrangement of the mass and momentum conservation laws of the flow, and are well suited for hybrid computations.
	The main idea is to feed in a one-way coupling the sound noise source induced by the flow field to an acoustic transport problem, see for instance Figure~\ref{fig:domainAeroProblem}. 
	Since the first development of aeroacoustics, hybrid methods have been established as a practical method for fast and accurate predictions for certain flow problems. 
 %

    In this work, we consider problems with high speed flow velocity, low Mach number and in which we assume that there is no interaction between the fluid pressure and the acoustic pressure, namely, where the one-way coupling hypothesis holds. Examples of physical relevance where such assumptions are verified and where recently hybrid methods have been successfully applied are \cite{Escobar2008, Kaltenbacher2013} for subsonic flows, \cite{Alipour2011, Huppe2014, Falk2021} for human phonation and \cite{Ewert2003} for airframe noise at low Mach.
	The greatest advantage of hybrid computational strategies for aeroacoustic problems is the possibility of choosing the optimal computational grids and optimal numerical methods for both the acoustic problem and the flow problem. 
    In the flow problem, sufficiently fine computational grid to resolve the turbulence or wall scales must be employed in order to properly describe the underlying physics of the flow. For the acoustic problem, the domain is usually much bigger and the involved acoustic pressure length scales are larger. 
    Moreover, the sound generation mechanism can be often confined only to the fluid region, allowing the acoustic problem to be modelled as a pure wave propagation problem. With this setting, optimal computational methods can be chosen for each problem independently. 
    In this work, we propose to employ Finite Volumes (FV) schemes for the fluid problem, since it is largely employed in the industrial framework. The corresponding FV solution is then employed to compute the source term for the acoustic wave propagation problem, which is then solved employing Spectral Element Methods (SEM).
    High order approximations have already been employed in computational aeroacoustic, see for instance \cite{HuppePHD} and \cite{Schlottke2016}.
	A critical step in hybrid methods is the coupling between the fluid and the acoustic domain, and how the noise source field is interpolated between the computational grids. 
    We remark that the fluid and acoustic grids might have very different granularities in order to capture the underlying (different) physics.
    While simple nearest neighbour interpolation fails to compute the acoustic sources accurately \cite{Caro2009}, conservative interpolation schemes lack of important analytical framework, although have been successfully employed, see for instance \cite{Escobar2010}. 
	The $L^2$-projection method has been employed in \cite{Schlottke2016} limited to nested Cartesian.
	In this work, we generalize the method to arbitrary polyhedral grids, proposing an efficient algorithm to compute the intersections between the fluid and acoustic elements. We then compute the projection of the sound source term computed as a post-process of the flow solution onto the acoustic grid by employing a quadrature free method on polyhedral elements. The proposed coupling strategy is flexible and acts as a black box, requiring only the sound source term at the cell centre of the fluid cell and hence it is well suited to be plugged onto any finite volume solver. Furthermore, it is naturally fitted for high order approximations since the employed quadrature formula integrates exactly arbitrary polynomials. We provide a rigorous theoretical analysis quantifying the effect of the projection error. This allows us to state that if the fluid grid is fine enough, we can exploit the accuracy of the spectral solver. The theoretical results are then validated by means of numerical experiments. The flow computations are performed with OpenFOAM \cite{OpenFOAM}, an open-source finite volume library, while the inhomogeneous acoustic wave equation is solved with SPEED \cite{SPEED2013}, an open-source spectral element library. \\
	The paper is structured as follows. 
    In Section~\ref{sec:AeroAcousticModel} we introduce the aeroacoustic hybrid problem, and we propose our strategy to solve the inhomogeneous Lighthill's wave equation. We focus our attention on the coupling between the fluid and acoustic problem. 
    We develop in Section~\ref{sec:Analysis} the theoretical analysis on the projection method and we then discuss in Section~\ref{sec:ImplementationAspects} the challenging implementation aspects.
    In Section~\ref{sec:IntersectionResults} we test the proposed intersection algorithm and we verify in Section~\ref{sec:ConvergenceResults} the theoretical estimates for the projection error. In Section~\ref{sec:AeroResults} we apply the proposed computational strategy on benchmark aeroacoustic problems.

	\section{The aeroacoustic model problem} \label{sec:AeroAcousticModel}
	It is possible to find in literature a wide variety of aeroacoustic models: from semi-analytical strategies based on employing suitable Green functions that led to the popular Curle \cite{Curle1955} and Ffowcs Williams Hawkings analogies \cite{FWH1969}, to more recent models that aim to solve the acoustic perturbed equations (APE), see for instance \cite{Ewert2003}.
	Most of the approximation methods proposed for these models rely on a hybrid strategy: first, they compute the fluid flow solution and then, they solve the acoustic problem using the latter to compute the sound source. This is the principle upon which the Lighthill's wave equation is based on.
	
	\subsection{Lighthill's wave equation}\label{sec:lighthill_wave}
	Let $\Omega \subset \mathbb{R}^3$, be a connected open bounded domain with sufficiently smooth boundary $\partial\Omega$. We denote by $\mathbf{x} \in \Omega$ the vector of spatial coordinates, and by $t \in (0,T]$ the time coordinate, being $T>0$ a final observation time.
	We consider in $\Omega \times (0,T]$ the compressible unsteady Navier-Stokes equations:
	\begin{align}
		\der{\rho}{t} + \nabla\cdot(\rho \mathbf{u}) &= 0 \label{eq:massBalance},\\
		\der{\rho\mathbf{u}}{t} + \nabla\cdot (\rho\mathbf{u}\otimes\mathbf{u}) &= -\nabla p + \nabla\cdot \bm{\sigma}  \label{eq:momentumBalance}, \\
		\frac{\partial \rho E}{\partial t} + \nabla \cdot [(\rho E + p)\mathbf{u}] &= \nabla \cdot (\bm{\sigma}\mathbf{u} - \mathbf{q}), \label{eq:energyBalance}
	\end{align}
  supplemented with suitable boundary, initial conditions and a state equation that will be detailed later on. Equations \eqref{eq:massBalance}, \eqref{eq:momentumBalance} and  \eqref{eq:energyBalance} are the mass, momentum and energy balance equations, respectively, where $\bm{\sigma}$ denotes the viscous stress tensor, $\rho$ is the fluid density, $p$ is the pressure, $\mathbf{u}$ is the fluid velocity, $\rho E$ is the total energy and $\mathbf{q}$ is the heat flux.
  \Ar{The Navier-Stokes equations \eqref{eq:massBalance}-\eqref{eq:energyBalance} could, in principle, be adopted to describe the full model for aeroacoustic problems. 
  However, due to the different scales in the acoustic and flow problems, directly solving the system of equations \eqref{eq:massBalance}-\eqref{eq:energyBalance} is still computationally unaffordable for far field noise resolution, see for instance \cite{Sagaut2007}. Hence, we need to propose a different strategy.
  }
We define the \textcolor{blue}{adimensional} Mach ($\textrm{Ma}$) and Reynolds ($\textrm{Re}$) numbers as 
	\begin{equation*}
		\textrm{Ma} = \displaystyle \frac{U}{c_0}, \qquad \textrm{Re} = \frac{UL}{\nu},  
	\end{equation*}
 where $U$ is the characteristic speed of the flow, $c_0$ is the speed of sound, $L$ is the characteristic length of the flow problem, and $\nu$ is the kinematic viscosity. 
 We derive the Lighthill's wave equation, see \cite{Lighthill1952}, by taking the time derivative of \eqref{eq:massBalance} and subtracting the divergence of the momentum equation \eqref{eq:momentumBalance}. Then, we have 
	\begin{equation*}
		\der{^2\rho}{t^2} = \nabla\cdot\nabla\cdot(\rho \mathbf{u}\otimes\mathbf{u} + p\mathbf{I} - \bm{\sigma}).
	\end{equation*}
	Summing and subtracting in the above equation the term $c_0^2 \Delta \rho$, we obtain the following wave equation:
	\begin{equation} \label{eq:FullLighthill}
		\der{^2\rho}{t^2} -c_0^2 \Delta \rho = \nabla\cdot \nabla\cdot \mathbf{T},
	\end{equation}
	where the right-hand side has been reformulated by introducing the so-called Lighthill's tensor 
 \begin{equation}\label{eq:Lighthill_full}
 \mathbf{T} = \rho \mathbf{u}\otimes\mathbf{u} + (p-c_0^2\rho)\mathbf{I} - \bm{\sigma},      
 \end{equation} 
 being $\mathbf{I}$  the identity tensor.
 The model \eqref{eq:FullLighthill}-\eqref{eq:Lighthill_full} can be further simplified depending on the problem of interest.  
 For a sufficiently high Reynolds number,
 it is possible to neglect the viscous source term in the Lighthill's tensor.
 Assuming a low Mach number and no combustion effects, the fluid can be considered isentropic, leading to $p = c_0^2 \rho$. Under these assumptions, the Lighthill's tensor in \eqref{eq:Lighthill_full} reduces to $\mathbf{T} = \rho_0 \mathbf{u}\otimes\mathbf{u}$, where $\rho_0$ is a reference density for the fluid.  This leads to the following wave equation:
	\begin{equation} \label{eq:FullLighthill2}
		\der{^2\rho}{t^2} -c_0^2 \Delta \rho = \nabla\cdot\nabla\cdot(\rho_0 \mathbf{u} \otimes \mathbf{u}),
	\end{equation}
which describes the evolution of a density wave in a quiescent material where the speed of propagation is given by $c_0$, the fluid \Ar{reference} density is given by $\rho_0$ and the sound source is given by the approximation of the Lighthill's stress tensor.
Equation~\ref{eq:FullLighthill2} is then supplemented with suitable initial and boundary conditions, as detailed in the following.
 
\subsection{The hybrid coupled model problem}

\begin{figure}[h]
		\begin{center} 
			\includegraphics[width=0.5\textwidth]
    {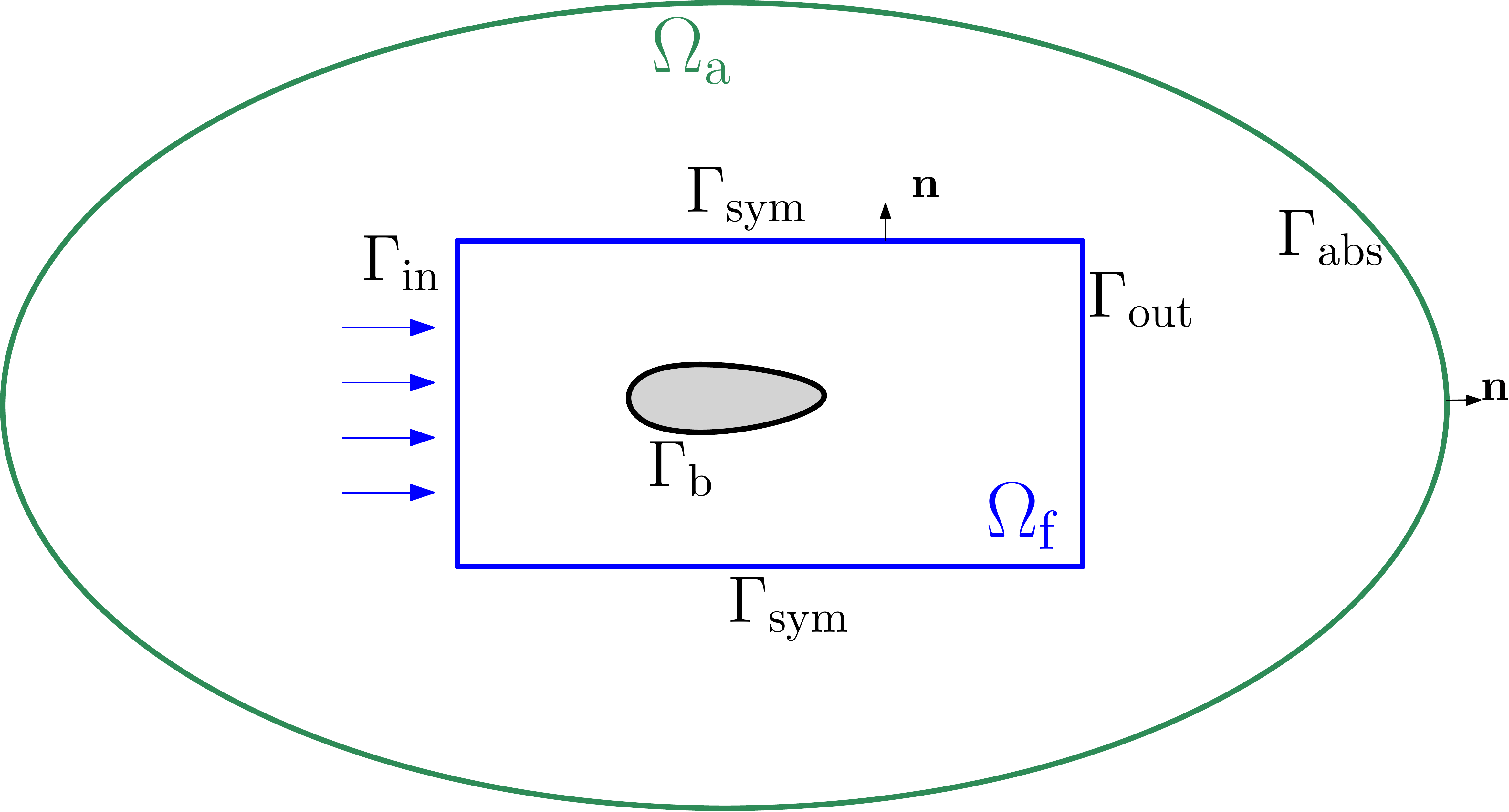}
			\caption{Computational domain for the aeroacoustic problem and sketch of the hybrid numerical strategy. First, the fluid problem is solved on $\Omega_{\textrm{f}}$. Then the acoustic source term is computed from the flow velocity. Finally, an inhomogeneous acoustic wave equation is solved on $\Omega_{\textrm{a}}$.} \label{fig:domainAeroProblem2}
		\end{center}
	\end{figure}

 With the aim of studying aeroacoustic problems related to the noise generated by external flows around bodies, we consider the following setup. We assume acoustic compactness, which means that the size of the flow source structures that generate the acoustic field are small compared to the acoustic generated wavelength. This hypothesis is inherently fulfilled for low Mach number applications.
Next, we consider a connected domain $\Omega_{\textrm{f}}$, having sufficiently regular boundary $\partial \Omega_{\textrm{f}}$, embedded in a connected domain $\Omega_{\textrm{a}}$, with sufficiently regular boundary $\partial \Omega_{\textrm{a}}$, see Figure~\ref{fig:domainAeroProblem2}.

The hybrid algorithm requires to solve the following sequence of problems:\\
\textit{a. Flow Problem}. For the flow problem we consider the incompressible Navier-Stokes equations, that read as
for $t \in (0,T]$, find $\mathbf{u}(\mathbf{x},t):\Omega_{\textrm{f}} \times (0,T]\rightarrow \mathbb{R}^3$ and $p(\mathbf{x},t):\Omega_{\textrm{f}} \times (0,T]\rightarrow \mathbb{R}$ such that
\begin{equation}\label{eq:NavierStokes}
		\begin{aligned} 
			\der{\mathbf{u}}{t} + \nabla\cdot(\mathbf{u}\otimes\mathbf{u}) - \nabla\cdot(\nu \nabla \mathbf{u}) + \frac{1}{\rho_0}\nabla p &= 0, \quad \text{in } \Omega_{\textrm{f}}\times(0,T],\\
			\nabla\cdot\mathbf{u} &= 0, \quad \text{in } \Omega_{\textrm{f}}\times(0,T],\\
			\mathbf{u}(\mathbf{x},0) &= \textbf{0}, \quad \text{in } \Omega_{\textrm{f}}, \\
			\mathbf{u} &= \textbf{0}, \quad \text{on } \Gamma_{\textrm{b}}, \\
			\mathbf{u} &= \mathbf{g}, \quad \text{on } \Gamma_{\textrm{in}}, \\
			\nu \nabla \mathbf{u}\cdot\mathbf{n} - p\mathbf{n} &= \textbf{0}, \quad \text{on } \Gamma_{\textrm{out}}, \\
 \mathbf{u}\cdot\mathbf{n} &= \mathbf{0},  \quad \text{on }  \Gamma_{\textrm{sym}}, \\ 
  \nabla(\mathbf{u} - (\mathbf{u}\cdot\mathbf{n})\mathbf{n})\cdot\mathbf{n} &= \mathbf{0},  \quad \text{on }  \Gamma_{\textrm{sym}},
		\end{aligned}
\end{equation}
	
	where $\mathbf{n}$ is the outward unit normal vector to $\partial\Omega_{\textrm{f}}$, $\nu$ is the kinematic viscosity, $\rho_0$ is the fluid density and $\mathbf{g}$ is the inlet Dirichlet datum. Here, we suppose the fluid boundary can be decomposed \Ar{ in the pairwise disjoint portions $\Gamma_{\textrm{in}},\Gamma_{\textrm{out}},\Gamma_{\textrm{b}},\Gamma_{\textrm{sym}}$, such that  $\partial\Omega_{\textrm{f}} = \Gamma_{\textrm{in}}\cup\Gamma_{\textrm{out}}\cup\Gamma_{\textrm{b}}\cup\Gamma_{\textrm{sym}}$}. \st{$\Gamma_{\textrm{in}}\cap\Gamma_{\textrm{out}}\cap\Gamma_{\textrm{b}}\cap\Gamma_{\textrm{sym}} = \emptyset$,}.

\textit{b. Acoustic Source.} From the fluid velocity $\mathbf{u}$ we define the Lighthill's tensor as 
\begin{equation}
\mathbf{T} = \begin{cases}\rho_0 \mathbf{u}\otimes\mathbf{u} \quad &\text{if } \mathbf{x} \in \Omega_{\textrm{f}}, \\
	\mathbf{0} \quad &\text{if } \mathbf{x} \in \Omega_{\textrm{a}}\backslash\Omega_{\textrm{f}}.
	\end{cases}
\end{equation}
The Lighthill's tensor has support only on the fluid domain $\Omega_{\textrm{f}} \subseteq \Omega_{\textrm{a}}$, and it depends on the solution $\mathbf{u}$ of problem \eqref{eq:NavierStokes}, being the coupling term between the flow problem \eqref{eq:NavierStokes} and the acoustic problem \eqref{eq:LighthillWaveContCont}.
 
\textit{c. Acoustic Problem.} We consider in $\Omega_{\textrm{a}}$  the following non-homogeneous acoustic problem based on the Lighthill's wave equation, cf. Section~\ref{sec:lighthill_wave}:
	for $t \in (0,T]$, find  $\rho(\mathbf{x},t):\Omega_{\textrm{a}} \times (0,T]\rightarrow \mathbb{R}$ such that
	\begin{equation}\label{eq:LighthillWaveContCont}
		\begin{aligned} 
			\der{^2\rho}{t^2} -c_0^2 \Delta \rho = \nabla\cdot \nabla\cdot \mathbf{T}, &\qquad \text{in } \Omega_{\textrm{a}} \times (0,T),   \\
			c_0^2 \der{\rho}{\mathbf{n}} = 0, & \qquad \text{on } \Gamma_{\textrm{b}}\times (0,T),\\
			\frac{1}{\rho_0}\der{\rho}{\mathbf{n}} = -\frac{1}{\rho_0 c_0}\der{\rho}{t}(\mathbf{x},t), & \qquad \text{on } \Gamma_{\textrm{abs}}\times (0,T), \\
			\rho(\mathbf{x},0) = 0, & \qquad \mathbf{x} \in \Omega_{\textrm{a}}, \\
			\der{\rho}{t}(\mathbf{x},0) = 0, & \qquad \mathbf{x} \in \Omega_{\textrm{a}},
		\end{aligned}
	\end{equation}	
where $c_0$ is the speed of propagation of the wave and $\rho_0$ is the fluid density.	The boundary $\partial\Omega_{\textrm{a}}$ has been split as  $\partial\Omega_{\textrm{a}} = \Gamma_{\textrm{abs}} \cup \Gamma_{\textrm{b}}$. On the external boundary $\Gamma_{\textrm{abs}}$, cf. Figure~\ref{fig:domainAeroProblem}, we apply non-reflective boundary conditions, see \cite{Engquist1977}, while on $\Gamma_{\textrm{b}}$ we set a sound hard boundary condition, modelling a rigid wall. Initial conditions are set to zero. 
We are aware that the validity of this hybrid strategy and the underlying one-way coupling assumption is strongly problem-specific, depending on the geometry of the problem and the flow features. However, this approach is widely used in the context of aeroacoustics simulations, see for instance \cite{Escobar2008, Kaltenbacher2013, Huppe2014, Ewert2003}.


\subsection{Discretization of the incompressible Navier-Stokes equations}\label{sec:disc_NS}
	
The fluid flow problem is solved by employing the library OpenFOAM \cite{OpenFOAM}, an open-source library based on the cell centered finite volume method \cite{Peric1999}.
We consider a polyhedral tessellation $\TauFHrom$ of the domain $\Omega_\textrm{f}$ and we indicate with $\mathbf{x}_0$ the barycentre of the convex polyhedral cell $ K_\textrm{f}\in\TauFHrom$. Then, we introduce the space of piecewise constant functions $V_\textrm{f} = \left\{ v_\textrm{f} \in L^2(\Omega_\textrm{f}):v_\textrm{f}|_{ K_\textrm{f}} \in \mathbb{P}^0( K_\textrm{f}), \forall  K_\textrm{f} \in \TauFHrom \right\}$, \Ar{where $\mathbb{P}^0( K_\textrm{f})$ is the space of the constant functions on the element $K_\textrm{f}$}, and with $N_\textrm{f} = \dim{V_\textrm{f}}$ and we denote with $\bm{V}_\textrm{f} = [V_\textrm{f}]^3$ the vector valued discrete space.
In order to obtain a finite volume discretization of problem \eqref{eq:NavierStokes}, we integrate the momentum equation over the polyhedron $ K_\textrm{f} \in \TauFHrom$, getting
\begin{align}
	\int_{ K_\textrm{f}} \der{\mathbf{u}}{t} d\mathbf{x} + \int_{ K_\textrm{f}} \nabla\cdot(\mathbf{u}\otimes\mathbf{u}) d\mathbf{x} - \int_{ K_\textrm{f}} \nabla\cdot(\nu \nabla \mathbf{u})  d\mathbf{x}+ \int_{ K_\textrm{f}} \nabla \left(\frac{p}{\rho_0}\right) d\mathbf{x} &= 0,\label{eq:finite_vol_1}\\
	\int_{ K_\textrm{f}} \nabla\cdot\mathbf{u} d\mathbf{x} &= 0,
\end{align}
and then proceed by discussing the discretization of each term, introducing $\mathbf{u}_h \in \bm{V}_\textrm{f}$ and $\mathbf{p}_h \in V_\textrm{f}$.
The spatial approximation of the first integral in \eqref{eq:finite_vol_1} is straightforward, namely, 
\begin{equation}\label{eq:FV_time_derivative}
	\int_{ K_\textrm{f}}\frac{\partial \mathbf{u}}{\partial t}  d\mathbf{x}\approx \int_{ K_\textrm{f}}\frac{\partial \mathbf{u}_h}{\partial t}  d\mathbf{x} = | K_\textrm{f}| \frac{\partial \mathbf{u}_h}{\partial t}(\mathbf{x}_0,
\end{equation}
where $| K_\textrm{f}|$ is the volume of the element $ K_\textrm{f}$ and where a mid-point quadrature rule is employed.
Next, being $\nu$ constant, we approximate  the third term of \eqref{eq:finite_vol_1} as follows  
\begin{equation}  \label{eq:DiffusionFV}
	\int_{ K_\textrm{f}} \nabla\cdot(\nu \nabla \mathbf{u})d\mathbf{x} = \int_{\partial K_\textrm{f}} (\nu \nabla \mathbf{u}) \mathbf{n}d s \approx \sum_{ F \in \partial K_\textrm{f}} \nu \nabla \mathbf{u}_{ F} \mathbf{n} | F|,
\end{equation}
where $\nabla\mathbf{u}_{ F} = \nabla\mathbf{u}(\mathbf{x}_{ F})$, being $\mathbf{x}_{ F}$ the face cell barycenter. Note that in the last step, we use a mid-point quadrature rule on the face $ F$. 
Now, if the face cell $ F$ is shared by two elements $ K_\textrm{f}^+$ and $ K_\textrm{f}^-$, we reconstruct linearly $\nabla\mathbf{u}_{ F} \mathbf{n}$, \Ar{(see appendix \ref{appendix:discFV}, eq. \eqref{eq:DiffusionGrad})}. 
Concerning the convective term in  \eqref{eq:finite_vol_1}, integrating by parts, we get: 
\begin{equation} \label{eq:convective}
	\int_{ K_\textrm{f}} \nabla\cdot(\mathbf{u}\otimes\mathbf{u}) d\mathbf{x} = \int_{\partial K_\textrm{f}} \mathbf{u} (\mathbf{u}\cdot\mathbf{n}) d s\approx \sum_{ F \in \partial K_\textrm{f}} \mathbf{u}_ F (\mathbf{u}_ F\cdot\mathbf{n} )| F|,
\end{equation}
where we applied a mid-point quadrature rule on the face $ F$.
\Ar{The flux term is discretized with a linear upwind scheme (see appendix \ref{appendix:discFV}, eq. \eqref{eq:linearUpwind})}.
Finally, the pressure gradient term is discretized similarly, by observing that $\nabla p$ = $\nabla\cdot(p\mathbf{I})$ and by applying the Gauss theorem. For the time discretization, we first divide the temporal interval $(0, T]$ into $N$ subintervals, such that $T = N\Delta t$, setting $t^n = n \Delta t $, with $n = 0,\dots,N$. We consider a backward differentiation formula of second order (BDF2) discretization scheme for \eqref{eq:FV_time_derivative},  namely $\displaystyle \frac{\partial \mathbf{u}_h}{\partial t} \approx \frac{3\mathbf{u}^{n+1}_h-4\mathbf{u}_h^{n} + \mathbf{u}^{n-1}_h}{2\Delta t}$.
Finally, we remark that we compute at any time $t^n$ the aeroacoustic sound source term as a post-process of the fluid solution $\mathbf{u}_h^n$, i.e., $\nabla\cdot \mathbf{T} \approx \nabla\cdot (\rho_0 \mathbf{u}_h^n\otimes\mathbf{u}_h^n)$, see \eqref{eq:BilinearContinuosProblem} and also \eqref{eq:convective}.

\subsection{Discretization of the Lighthill's wave equation}\label{sec:sem-wave}

 We start by considering the variational formulation of the acoustic problem \eqref{eq:LighthillWaveContCont}:
\textrm{for $t \in (0;T]$, find $\rho(\mathbf{x},t)\in H^1(\Omega_{\textrm{a}})$ such that $\forall w \in H^1(\Omega_{\textrm{a}})$:} 
	\begin{align}
		\left(\der{^2\rho}{t^2}, w\right)_{\Omega_{\textrm{a}}} + c_0^2(\nabla \rho ,\nabla w)_{\Omega_{\textrm{a}}}  + c_0 \int_{\Gamma_{\textrm{abs}}}\frac{\partial \rho}{\partial t} w\ ds= - (\nabla\cdot\mathbf{T}, \nabla w)_{\Omega_{\textrm{a}}}, 
  \label{eq:BilinearContinuosProblem}  
\end{align}
 with initial conditions $\rho = \displaystyle\frac{\partial \rho}{\partial t} = 0$ in $\Omega_{\textrm{a}} \times \{0\}$,
 being $(\cdot,\cdot)_{\Omega_{\textrm{a}}}$ the $L^2$ product over the domain $\Omega_{\textrm{a}}$.
 %
 \Ar{We remark that we integrated by parts the term $ (\nabla\cdot\nabla\cdot\mathbf{T}, w)_{\Omega_{\textrm{a}}}$ and the resulting boundary terms are null both on $\Gamma_{\textrm{abs}}$ and on $\Gamma_{\textrm{b}}$, as it discussed in \cite{Escobar2010}.	
} 
 %
%
Next, we discretize problem \eqref{eq:BilinearContinuosProblem} by means of the SEM as follows. We introduce a conforming decomposition $\TauAHrom$ of the domain $\Omega_{\textrm{a}}$ made by hexahedral elements $ K_{\textrm{a}}$.
We denote by $\widehat{ K}$ the reference element $[-1,1]^3$, and we suppose that for any mesh element $ K_{\textrm{a}} \in \TauAHrom$ there exists a suitable trilinear invertible map $\bm{\theta}_{ K_{\textrm{a}}}:\widehat{ K}\rightarrow  K_{\textrm{a}}$ with positive Jacobian $\mathbf{J}_{ K_{\textrm{a}}}$. We define the characteristic mesh dimension as $\displaystyle h_{\textrm{a}} = \max_{ K_{\textrm{a}} \in \TauAHrom} h_{ K_{\textrm{a}}}$, being $h_{ K_{\textrm{a}}}$ the diameter of the element $ K_{\textrm{a}}$. Next, we introduce  
the finite-dimensional space: $V_{\textrm{a}} = \left\{ v \in C^0(\overline{\Omega}_\textrm{a})\cap H^1(\Omega_{\textrm{a}}): v|_{ K_{\textrm{a}}}\circ \bm{\theta}_{ K_{\textrm{a}}}^{-1}\in \mathbb{Q}_r(\widehat{ K}), \forall  K_{\textrm{a}} \in \TauAHrom\right\}$, where  $\mathbb{Q}_r(\widehat{ K})$ is the space of  polynomials of degree less than or equal to $r\geq 1$ in each coordinate direction, and we denote by $N_\textrm{a}$ the dimension of $V_{\textrm{a}}$.
Next, for any $u,w \in V_{\textrm{a}}$, we introduce  
the following bilinear form by means of the Gauss-Legendre-Lobatto (GLL) quadrature rule:
\begin{equation}
	(u,w)_{ K_{\textrm{a}}}^{\textrm{NI}}  = \sum_{i,j,k=0}^{r}  u(\bm{\theta}_{ K_{\textrm{a}}}(\bm{\xi}_{i,j,k}^{\textrm{GLL}}))w(\bm{\theta}_{ K_{\textrm{a}}}(\bm{\xi}_{i,j,k}^{\textrm{GLL}})) \omega_{i,j,k}^{\textrm{GLL}} |\det(\mathbf{J})| \approx (u,w)_{ K_{\textrm{a}}}
\end{equation}
where $\bm{\xi}^{\textrm{GLL}}$ are the GLL quadrature nodes, and $\omega^{\textrm{GLL}}$ their corresponding weights, defined in $[-1,1]^3$ (cf. \cite{Quarteroni2009}) and $\textrm{NI}$ stands for numerical integration. Moreover, we define
\begin{equation*}
	(u,w)_{\TauAHrom}^{\textrm{NI}} = \sum_{ K_{\textrm{a}} \in \TauAHrom} (u,w)_{ K_{\textrm{a}}}^{\textrm{NI}} \quad \forall \, u,w \in V_{\textrm{a}}.
\end{equation*}
The semi-discrete spectral element formulation of problem \eqref{eq:BilinearContinuosProblem} with numerical integration (SEM-NI) reads:
\textrm{for any time $t\in(0;T]$ find $\rho_h \in V_{\textrm{a}}$ such that:}
	\begin{equation}
		(\derSecond{\rho_h}{t}{2},w_h)_{\TauAHrom}^{\textrm{NI}}+  c_0^2(\nabla \rho_h ,\nabla w_h)_{\TauAHrom}^{\textrm{NI}}  + c_0  \int_{\Gamma_{\textrm{abs}}}\frac{\partial \rho_h}{\partial t} w_h \ ds = -(\nabla\cdot\mathbf{T},\nabla w_h)_{\TauAHrom}^{\textrm{NI}} \quad \forall w_h \in V_{\textrm{a}},
		\label{eq:SemiDiscreteSEM}
	\end{equation} 
with $\rho_h = \displaystyle\frac{\partial \rho_h}{\partial t} = 0$ in $\Omega_{\textrm{a}} \times \{ 0 \}$.	We recall that the term $\nabla \cdot \mathbf{T}$ is an external source that in our case is obtained from a numerical solution of problem \eqref{eq:NavierStokes} as described in Section~\ref{sec:disc_NS}.
In the next section, we detail how to compute effectively the right-hand side of eq.\eqref{eq:SemiDiscreteSEM}, i.e., how to approximate a field defined on the fluid mesh $\TauFHrom$ with a field defined on the acoustic grid $\TauAHrom$.
	
\subsection{$L^2-$projection of the acoustic source}\label{sec:L2Projection}

Let $q_{\textrm{f}} \in V_{\textrm{f}}$ be a function defined on the fluid grid $\TauFHrom$ such that $q_{\textrm{f}} = \sum_{i=1}^{N_\textrm{f}} \widehat{q}_{\textrm{f},i} \phi_{\textrm{f},i}$, where  $\spanOp{\phi_{\textrm{f},i}}_i^{N_\textrm{f}}$ is the set of $N_\textrm{f}$ basis functions associated to $V_{\textrm{f}}$, and $\widehat{q}_{\textrm{f},i}$ are the corresponding expansion coefficients. 
We define the $L^2$-projection of the field $q_{\textrm{f}}$ into $V_{\textrm{a}}$ as 
\begin{equation} \label{eq:ProjMin}
	 q_{\textrm{a}} =\argmin_{q \in V_{\textrm{a}}}\norm{q_{\textrm{f}}-q}{{L^2(\TauAHrom)}},
\end{equation}
\Ar{where $q_{\textrm{f}}$ has been extended by zero also on $\Omega_{\textrm{a}}$}.
Problem (\ref{eq:ProjMin}) is equivalent to the following: \textrm{find $q_{\textrm{a}} \in V_{\textrm{a}}$ s.t. } 
\begin{align}
	( q_{\textrm{a}}, \phi_{\textrm{a},i})_{\TauAHrom} &= ( q_{\textrm{f}},  \phi_{\textrm{a},i})_{\TauAHrom}  \quad \forall \phi_{\textrm{a},i} \in V_{\textrm{a}}, \label{def:ProjectionA}
\end{align}
where $q_{\textrm{a}} \in V_{\textrm{a}}$ is a function defined on the acoustic grid $\TauAHrom$ such that $q_{\textrm{a}} = \sum_{i=1}^{N_\textrm{a}} \widehat{q}_{\textrm{a},i} \phi_{\textrm{a},i}$, where  $\spanOp{\phi_{\textrm{a},i}}_i^{N_\textrm{a}}$ is the set of $N_\textrm{a}$ basis functions, and $\widehat{q}_{\textrm{a},i}$ are the corresponding expansion coefficients.
Motivated by the solution method used in Section~\ref{sec:disc_NS} we address the case where $q_{\textrm{f}}$ is a piecewise constant \Ar{function} over $\TauFHrom$ \st{, namely $q_{\textrm{f}} \in V_{\textrm{f}}$}. 
Then, problem (\ref{def:ProjectionA}) can be recast as follows:
\begin{align}  \label{eq:ProjectionProblemDiscrete}
\sum_{ K_{\textrm{a}} \in \TauAHrom} ( q_{\textrm{a}}, \phi_{\textrm{a},i})_{ K_{\textrm{a}}} = \sum_{ K_{\textrm{a}} \in \TauAHrom} ( \sum_{\ell=1}^{N_\textrm{f}} \widehat{q}_{\textrm{f},\ell} \phi_{\textrm{f},\ell} , \phi_{\textrm{a},i})_{ K_{\textrm{a}}} = \sum_{ K_{\textrm{a}} \in \TauAHrom} \sum_{\ell=1}^{N_\textrm{f}}  \widehat{q}_{\textrm{f},\ell} ( 1, \phi_{\textrm{a},i})_{ K_{\textrm{a}}\cap  K_{\textrm{f},\ell} }, 
\end{align}
where we have used that $ K_{\textrm{f},\ell} =  \supp(\phi_{\textrm{f},\ell})$.
The discrete algebraic counterpart of  \eqref{eq:ProjectionProblemDiscrete} becomes
\begin{equation} \label{eq:AlgebraicProblemProjection}
	\textrm{M}^{\textrm{aa}} \widehat{\mathbf{q}}_\textrm{a} = \textrm{M}^{\textrm{af}} \widehat{\mathbf{q}}_{\textrm{f}},
\end{equation}
where $\textrm{M}^{\textrm{aa}} \in \mathbb{R}^{N_\textrm{a}\times N_\textrm{a}}$ is the acoustic mass matrix, i.e., 
\begin{equation} \label{eq:massProjectionProblem}
	\textrm{M}^{\textrm{aa}}_{i,j} = \sum_{ K_{\textrm{a}} \in \TauAHrom} ( \phi_{\textrm{a},j},\phi_{\textrm{a},i})_{ K_{\textrm{a}}}, \qquad i,j=1,\dots,N_\textrm{a}, 
\end{equation}
while $\textrm{M}^{\textrm{af}} \in \mathbb{R}^{N_\textrm{a}\times N_\textrm{f}}$  is the coupling mass defined as
\begin{equation}
\textrm{M}^{\textrm{af}}_{i,\ell} =  \sum_{ K_{\textrm{a}} \in \TauAHrom} \sum_{\ell=1}^{N_\textrm{f}}  ( 1, \phi_{\textrm{a},i})_{ K_{\textrm{a}}\cap  K_{\textrm{f},\ell} }, \qquad \ell=1,\dots,N_\textrm{f}. \label{eq:rhsProj}
\end{equation}
The coupling mass has been computed with a suitable quadrature formula that will be described in Section~\ref{sec:quadrature_free}.
The vector $\widehat{\mathbf{q}}_\textrm{a}$ in \eqref{eq:AlgebraicProblemProjection} collects all the expansion coefficients of the acoustic field $q_{\textrm{a}}$, while $\widehat{\mathbf{q}}_\textrm{f}$ collects all the expansion coefficients of the fluid field $q_{\textrm{f}}$.

\section{Error analysis for the acoustic source} \label{sec:Analysis}

It is evident that the accuracy of the numerical solution $\rho_h$  in \eqref{eq:SemiDiscreteSEM} strongly depends on the approximation of the acoustic source, namely, $\nabla \cdot \mathbf{T}$. In our case the latter is obtained as a post-process of the numerical solution $\mathbf{u}_h$ of the flow problem. Quantifying the projection error between the acoustic and fluid grids is therefore of paramount importance. 

However, before presenting the main result of the section we need to introduce some preliminary results. 

\begin{lemma}(Interpolation error on \textnormal{GLL} nodes). \label{lemma:InterpolantGLL}
	Given $f \in H^s(\Omega)$ for some $s\ge 1$, consider the Lagrangian interpolant $I_{\textnormal{a}}^{\textnormal{GLL}} f$ at the Gauss Legendre Lobatto nodes, where $r$ denotes the polynomial degree of the interpolant function and $h$ is the mesh size of $\Tau_h$ tessellation of $\Omega$. Assuming $h$ to be quasi uniform, we have that:
	\begin{equation}
		\norm{f - I_{\textnormal{a}}^{\textnormal{GLL}}f}{L^2(\Omega)} \lesssim h^{\min(r+1,s)} \left(\frac{1}{r}\right)^s \norm{f}{H^s(\Omega)}.
	\end{equation} 
\end{lemma}
\st{For more details} \Ar{The proof of Lemma \ref{lemma:InterpolantGLL} can be found in \cite[Section (5.4.3)]{Canuto2007}}.

\begin{lemma}(hp-inverse inequality)
\label{remark:PoincareInverse} 
	Assume now that $ K = \bm{\theta}_{ K}(\widehat{ K})$ is a hexahedral element s.t. $ K \subset \mathbb{R}^3$, where $\bm{\theta}_{ K}$ is a suitable trilinear map. Then we have that:
	\begin{equation}
		\norm{ \nabla v}{L^2( K)} \lesssim  \frac{r^{2}}{h_ K} \norm{v}{L^2( K)}, \quad \forall v \circ \bm{\theta}_{ K}^{-1} \in \mathbb{Q}_r(\widehat{ K}),
	\end{equation}
	where $h_ K$ is $\diam{( K)}$.
\end{lemma}
\st{For more details} \Ar{For the proof}, see \cite[Theorem 4.76]{Schwab1998}.
Finally, we recall this Poincaré-Friedrich like inequality:
\begin{lemma} \label{remark:PoincareFriedrich}
	Given $\displaystyle u \in H^1( K)$, where $ K$ is an open bounded convex domain in $\mathbb{R}^d$ and $\displaystyle u \in L_0^2( K) = \left\{ v \in L^2( K) : \int_{ K}v = 0 \ \text{d}\mathbf{x} \right\}$ then we have that:
	\begin{equation}
		\norm{u}{L^2( K)} \lesssim \frac{\diam( K)^{1+d/2}}{| K|^{1/2}}\norm{\nabla u}{L^2( K)}.
	\end{equation}
\end{lemma}
We refer the reader to \cite[Corollary 3.4]{Zheng2005} and to  \cite[Remark 5.8]{Antonietti2019} for further details \Ar{on the proof}.

\begin{figure}[h]
	\centering
	\includegraphics[width=0.8\textwidth]{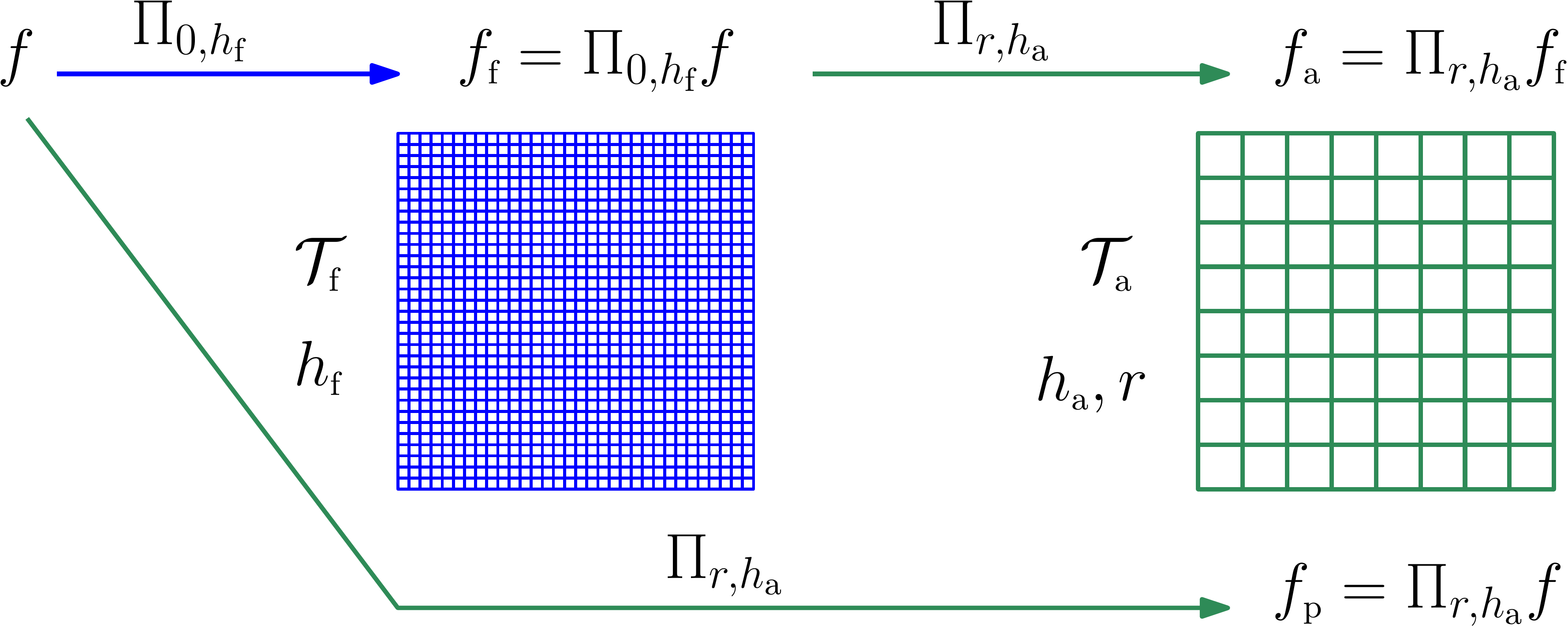}
	\caption{Schematic representation of the computational grids and the corresponding projection operators. For a given function $f \in L^2(\Omega)$, $f_\textrm{f}$ is the projection $\Pi_{0,h_{\textrm{f}}}f$, $f_\textrm{a}$ is the projection $\Pi_{r,h_{\textrm{a}}}f_\textrm{f}$ and $f_\textrm{p}$ is the projection $\Pi_{r,h_{\textrm{a}}}f$. The computational grids $\TauAHrom$ and $\TauFHrom$ are assumed to be nested. \label{fig:notationApprox}}
\end{figure}

Next, for the sake of the presentation, we consider the following setup: let $\Omega = \Omega_{\textrm{f}}=\Omega_{\textrm{a}}$ be a polygonal domain and let $\TauFHrom$ and $\TauAHrom$ be two nested grids of $\Omega$ as shown in Figure~\ref{fig:notationApprox}, namely
for all elements $ K_{\textrm{a}} \in \TauAHrom$ we assume that there exists a set of index $\mathcal{L}_{ K_{\textrm{a}}}$ such that $  \displaystyle K_{\textrm{a}} = \bigcup_{l\in\mathcal{L}_{ K_{\textrm{a}}}}   K_{\textrm{f},l}$.

We now introduce the following projection operators: $\Pi_{0,{h_{\textrm{f}}}} : L^2(\Omega) \rightarrow V_{\textrm{f}}$ and $\Pi_{r,{h_{\textrm{a}}}} : L^2(\Omega) \rightarrow V_{\textrm{a}}$, and we define the following functions:
\begin{equation}
    f_\textrm{f} = \Pi_{0,{h_{\textrm{f}}}}f,
\end{equation} 
that is the $L^2$ projection of $f \in L^2(\Omega)$ onto the space $V_{\textrm{f}}$,
\begin{equation}f_\textrm{p} = \Pi_{r,{h_{\textrm{a}}}} f,
\end{equation}
that is the $L^2$ projection of $f \in L^2(\Omega)$ onto the space $V_{\textrm{a}}$,
\begin{equation}f_\textrm{a} = \Pi_{r,{h_{\textrm{a}}}} f_\textrm{f},
\end{equation}
that is the $L^2$ projection of $f_\textrm{f} \in V_{\textrm{f}}$ onto the space $V_{\textrm{a}}$,
see Figure~\ref{fig:notationApprox}. \\
Now, we can state the following result.

\begin{theorem} (Approximation Theorem). \label{theorem:approximationTheorem}
	Let $\TauAHrom$ and $\TauFHrom$ be two grids of the same computational domain $\Omega_{\textnormal{f}} = \Omega_{\textnormal{a}} = \Omega$ made by hexahedral elements, such that $\TauFHrom$ is nested to $\TauAHrom$, namely, for every element $ K_{\textnormal{f}}$ there exists $ K_{\textnormal{a}}$ such that $ K_{\textnormal{f}} \subset  K_{\textnormal{a}}$. Given $f \in H^s(\Omega)$ with $s \ge 1$, let $f_\textnormal{f} = \Pi_{0,h_{\textnormal{f}}} f$ be the projection of $f$ onto the space $V_{\textnormal{f}}$ and let $f_\textnormal{a} = \Pi_{r,h_{\textnormal{a}}} f_\textnormal{f}$, namely the projection of $f_\textnormal{f}$ onto the space $V_{\textnormal{a}}$. Then, it holds:
	\begin{equation}
		\norm{f-f_\textnormal{a}}{L^2(\Omega)} \lesssim  h_{\textnormal{a}}^{\min(r+1,s)} \left(\frac{1}{r}\right)^s \norm{f}{H^s(\Omega)} + \frac{h_{\textnormal{f}}^2}{h_{\textnormal{a}}}r^2\norm{f}{H^s(\Omega)}. \label{eq:estimateApprox}
	\end{equation}
\end{theorem}

\textit{Proof}. Let $f_\textrm{p} = \Pi_{r,h_{\textrm{a}}} f$, see for instance Figure~\ref{fig:notationApprox}. By triangular inequality we have that
\begin{equation}
	\norm{f-f_\textrm{a}}{L^2(\Omega)} \le 
	\norm{f-f_\textrm{p}}{L^2(\Omega)} + \norm{f_\textrm{p}-f_\textrm{a}}{L^2(\Omega)}.
\end{equation}
The first term on the right hand side can be estimated by employing Lemma~\ref{lemma:InterpolantGLL}, i.e.:
\begin{equation} \label{eq:ApplicationInterpolantGLL}
	\begin{aligned}
		\norm{f-f_\textrm{p}}{L^2(\Omega)} &= \min_{\varphi \in V_{\textrm{a}}}  \norm{f-\varphi}{L^2(\Omega)} \lesssim \norm{f-I_{\textrm{a}}^{\textrm{GLL}}f}{L^2(\Omega)} \lesssim h_{\textrm{a}}^{\min(r+1,s)} \left(\frac{1}{r}\right)^s \norm{f}{H^s(\Omega)}. 
	\end{aligned}
\end{equation}
Next, we observe that by definition of the $L^2$-projection we get
\begin{equation}
	(f_\textrm{a},\phi)_{L^2(\Omega)} = (f_\textrm{f},\phi)_{L^2(\Omega)} \quad \forall \phi \in V_{\textrm{a}} \label{eq:ProiezioneFA},
\end{equation}
\begin{equation}
	(f_\textrm{p},\phi)_{L^2(\Omega)} = (f,\phi)_{L^2(\Omega)} \quad \forall \phi \in V_{\textrm{a}} \label{eq:ProiezioneFP}.
\end{equation}
Then, by subtracting \eqref{eq:ProiezioneFA} to \eqref{eq:ProiezioneFP}, we obtain
\begin{equation*}
	(f_\textrm{p} - f_\textrm{a},\phi)_{L^2(\Omega)} = (f - f_\textrm{f},\phi)_{L^2(\Omega)}, \quad \forall \phi \in V_{\textrm{a}}.
\end{equation*}
Furthermore, since $f_\textrm{p} - f_\textrm{a} \in V_{\textrm{a}}$, we can write
\begin{equation*}
	\begin{aligned}
		\norm{f_\textrm{p}-f_\textrm{a}}{L^2(\Omega)}^2 &= \int_{\Omega} (f - f_\textrm{f})(f_\textrm{p}-f_\textrm{a}) \ \text{d}\mathbf{x} = \sum_{ K_{\textrm{f}}} \int_{ K_{\textrm{f}}}  (f - f_\textrm{f})(f_\textrm{p}-f_\textrm{a}) \ \text{d}\mathbf{x}.
	\end{aligned}
\end{equation*}
and notice that 
\begin{equation*}
	(f - f_\textrm{f}, \varphi)_{L^2( K_{\textrm{f}})} = 0 \quad \forall \varphi \in \mathbb{P}^0( K_{\textrm{f}}),
\end{equation*}
where $\mathbb{P}^0( K_{\textrm{f}})$ is the space of the constant functions over $ K_{\textrm{f}}$.
By taking $\varphi = \Pi_{0,h_{\textrm{f}}} (f_\textrm{p} - f_\textrm{a})$ in the above equation yields to 
\begin{equation}
	\begin{aligned}
		\sum_{ K_{\textrm{f}}} \int_{ K_{\textrm{f}}}  (f - f_\textrm{f})(f_\textrm{p}-f_\textrm{a}) \ \text{d}\mathbf{x} &= 
		\sum_{ K_{\textrm{f}}} \int_{ K_{\textrm{f}}}  (f - f_\textrm{f})(f_\textrm{p}-f_\textrm{a} - \Pi_{0,h_{\textrm{f}}} (f_\textrm{p} - f_\textrm{a})) \ \text{d}\mathbf{x}. \\
		& \lesssim \sum_{ K_{\textrm{f}}} \norm{f - f_\textrm{f}}{L^2( K_{\textrm{f}})}\norm{f_\textrm{p} - f_\textrm{a}  - \Pi_{0,h_{\textrm{f}}} (f_\textrm{p} - f_\textrm{a})}{L^2( K_{\textrm{f}})} \\
		& \lesssim h_{\textrm{f}}^2 \sum_{ K_{\textrm{f}}} \norm{\nabla f}{L^2( K_{\textrm{f}})} \norm{\nabla (f_\textrm{p}-f_\textrm{a}) }{L^2( K_{\textrm{f}})} ,
	\end{aligned}
\end{equation}
where in the last inequality we employ Lemma~\ref{remark:PoincareFriedrich}. By linearity of the integral, noticing that by hypothesis $ K_{\textrm{f}} \subset  K_{\textrm{a}}$ and using  that $ \norm{\nabla f}{L^2( K_{\textrm{f}})} \le  \norm{\nabla f}{L^2( K_{\textrm{a}})}$ we obtain
\begin{equation*}
	\begin{aligned}
		\sum_{ K_{\textrm{f}}} h_{\textrm{f}}^2 \norm{\nabla f}{L^2( K_{\textrm{f}})} \norm{\nabla (f_\textrm{p}-f_\textrm{a}) }{L^2( K_{\textrm{f}})}  
		& \lesssim  \sum_{ K_{\textrm{a}}} h_{\textrm{f}}^2 \norm{\nabla f}{L^2( K_{\textrm{a}})}\norm{\nabla(f_\textrm{p}-f_\textrm{a})}{L^2( K_{\textrm{a}})} \\
		& \lesssim  \sum_{ K_{\textrm{a}}} \frac{h_{\textrm{f}}^2 r^2}{h_{\textrm{a}}} \norm{\nabla f}{L^2( K_{\textrm{a}})} \norm{f_\textrm{p}-f_\textrm{a} }{L^2( K_{\textrm{a}})} \\
		& \lesssim  \frac{h_{\textrm{f}}^2 r^2}{h_{\textrm{a}}} \norm{\nabla f}{L^2(\Omega)} \norm{f_\textrm{p}-f_\textrm{a} }{L^2(\Omega)},
	\end{aligned}
\end{equation*}
where in the last step we used the inverse inequality of Lemma~\ref{remark:PoincareInverse}.
Finally, we get
\begin{equation*}
	\norm{f_\textrm{p}-f_\textrm{a}}{L^2(\Omega)}^2 \lesssim  \frac{h_{\textrm{f}}^2 r^2}{h_{\textrm{a}}} \norm{\nabla f}{L^2(\Omega)} \norm{f_\textrm{p}-f_\textrm{a} }{L^2(\Omega)},
\end{equation*}
or equivalently,
\begin{equation}\label{eq:fp-fa}
	\norm{f_\textrm{p}-f_\textrm{a}}{L^2(\Omega)} \lesssim \frac{h_{\textrm{f}}^2 r^2}{h_{\textrm{a}}}  \norm{\nabla f}{L^2(\Omega)}
\end{equation}
and, since $\norm{\nabla f}{L^2(\Omega)} \lesssim \norm{f}{H^s(\Omega)}$, that concludes the proof.
\qedwhite

\section{Implementation aspects} \label{sec:ImplementationAspects}
An accurate solution of the projection problem \eqref{eq:ProjectionProblemDiscrete} requires computing the intersection between the elements $ K_{\textrm{a}} \in \TauAHrom$ and $ K_{\textrm{f}} \in \TauFHrom$. This operation is in general very expensive, but in many applications, it is crucial to compute it accurately in order to have reliable solutions. 
\textcolor{blue}{In this section, we consider $\TauAHrom$ and $\TauFHrom$ (not necessary nested),} and we present our strategy to compute the intersection between two elements $ K_{\textrm{a}}$ and $ K_{\textrm{f}}$ and we show that it is robust and scalable. We recall that the intersection $ K =  K_{\textrm{a}} \cap  K_{\textrm{f}}$ is in general a polyhedron in the three-dimensional space. 
Moreover, we describe the employed quadrature-free algorithm to compute the integral of polynomials functions over $ K$, cf. Equation~\eqref{eq:rhsProj}.

\subsection{Intersection algorithm} \label{sec:IntersectionAlgorithm}

\begin{figure}

\begin{center}
	\begin{subfigure}{\textwidth}
	\begin{center}
	\begin{overpic}[width=0.6\textwidth]{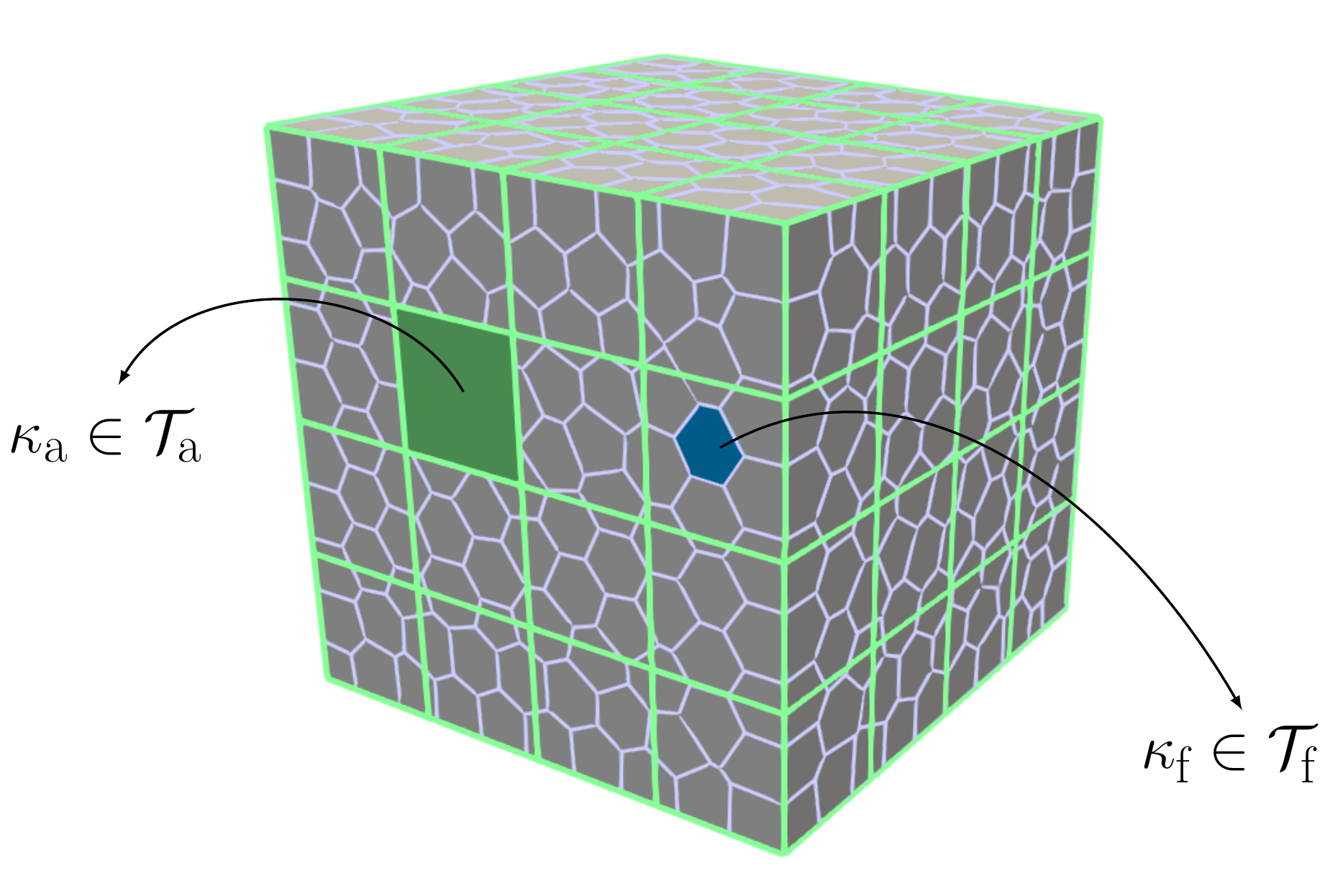}
\end{overpic}
	\end{center}
	\caption{Example of acoustic $\TauAHrom$ and fluid $\TauFHrom$ tessellations made of hexahedral and polyhedral elements, respectively, for the domain $\Omega = \Omega_{\textrm{a}} =\Omega_{\textrm{f}}$.\label{fig:domainExample}}
	\end{subfigure}
\end{center}
\begin{subfigure}{\textwidth}
	\adjustbox{valign=m}{
	\includegraphics[width=0.4\textwidth]{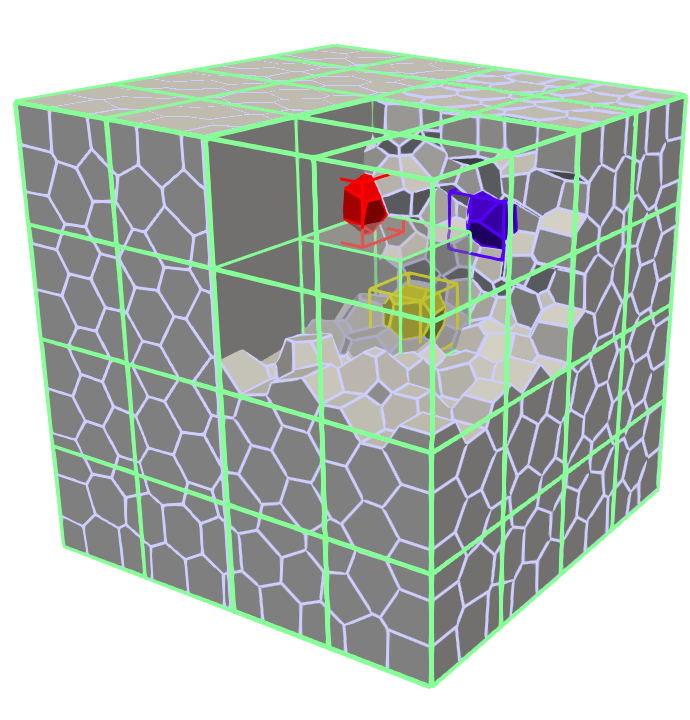}} \hfill
	\adjustbox{valign=m}{
 \includegraphics[width=0.4\textwidth]{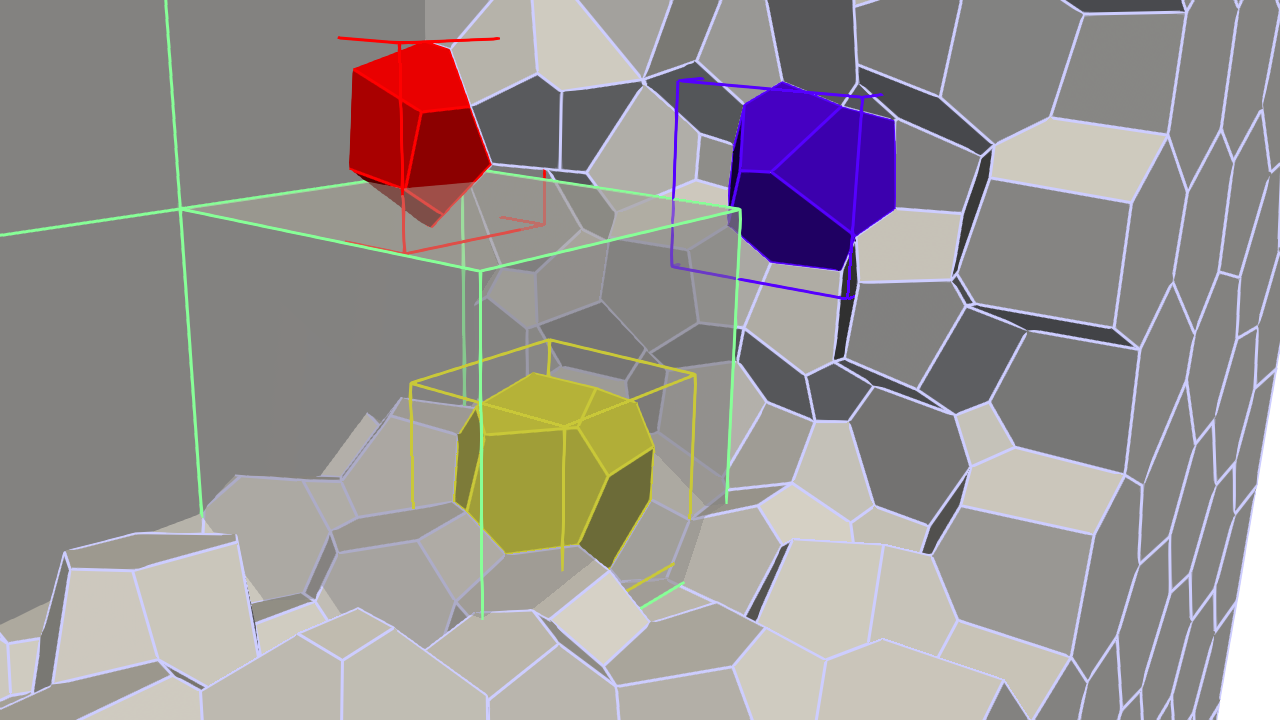}
	}
\caption{\textit{Bounding Box search}. Left, overview on the whole domain. Right, zoom on a selected acoustic element $ K_{\textrm{a}}$. At this stage, all the bounding box of the fluid elements intersecting $\mathcal{B}( K_{\textrm{a}})$ are selected and collected in the set $\mathcal{K}_\textrm{a}$. We show some of the selected fluid elements (yellow \textcolor{yellow}{$\bullet$}, red \textcolor{red}{$\bullet$} and blue \textcolor{blue}{$\bullet$}) and their respective bounding box. \label{fig:BBOXsearch}}
\end{subfigure}
\caption{Schematic representation of Algorithm~\ref{alg:cutCell}. (\textcolor{myGreen}{$\bullet$}) Acoustic element $ K_{\textrm{a}}$. (\textcolor{red}{$\bullet$}) Example of a fluid element $ K_{\textrm{f}} \in \mathcal{I}_\textrm{a}$. (\textcolor{blue}{$\bullet$}) Example of a fluid element $ K_{\textrm{f}} \in \mathcal{K}_\textrm{a}$ but not intersecting with the selected element $ K_{\textrm{a}}$. (\textcolor{yellow}{$\bullet$}) Example of a fluid element $ K_{\textrm{f}} \in \mathcal{C}_\textrm{a}$.
}
\end{figure}
\begin{figure}
\begin{subfigure}{\textwidth}
	\begin{minipage}{0.25\textwidth}
	\includegraphics[width=\textwidth]{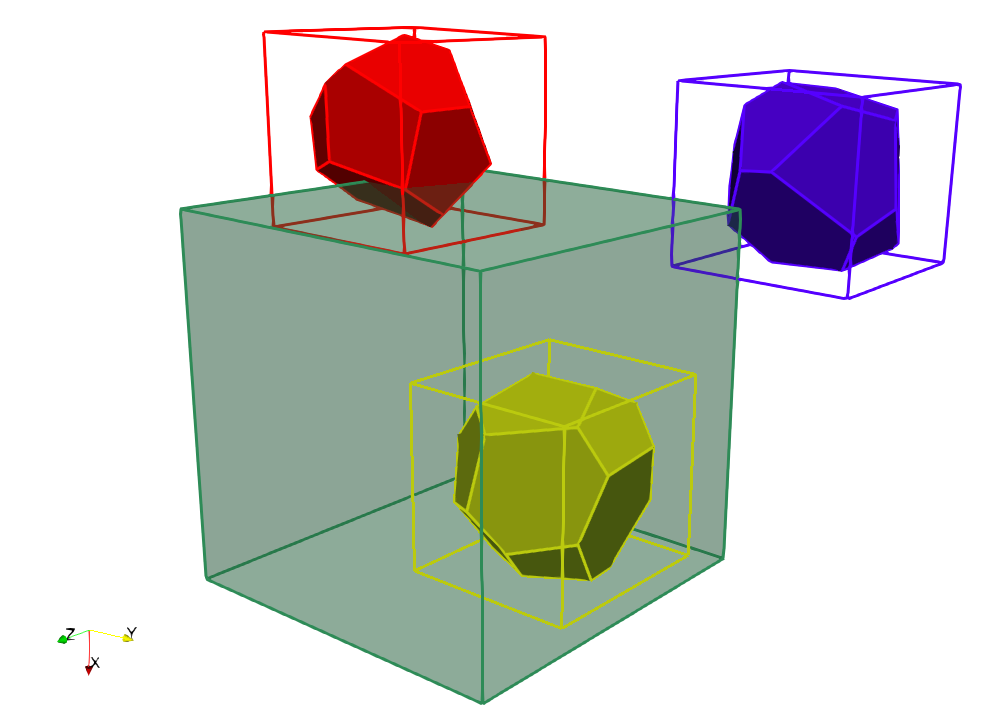}
\end{minipage}\hfill
\begin{minipage}{0.7\textwidth}
	\caption{\textit{Bounding Box selection}. At this stage, we check if the bounding box of a fluid element is contained inside $ K_{\textrm{a}}$. For example, the yellow (\textcolor{yellow}{$\bullet$}) fluid element is contained in $ K_{\textrm{a}}$, and hence it is added to $\mathcal{C}$. \label{fig:Richard}}
\end{minipage}
\end{subfigure}
\begin{subfigure}{\textwidth}
	\begin{minipage}{0.25\textwidth}
	\includegraphics[width=\textwidth]{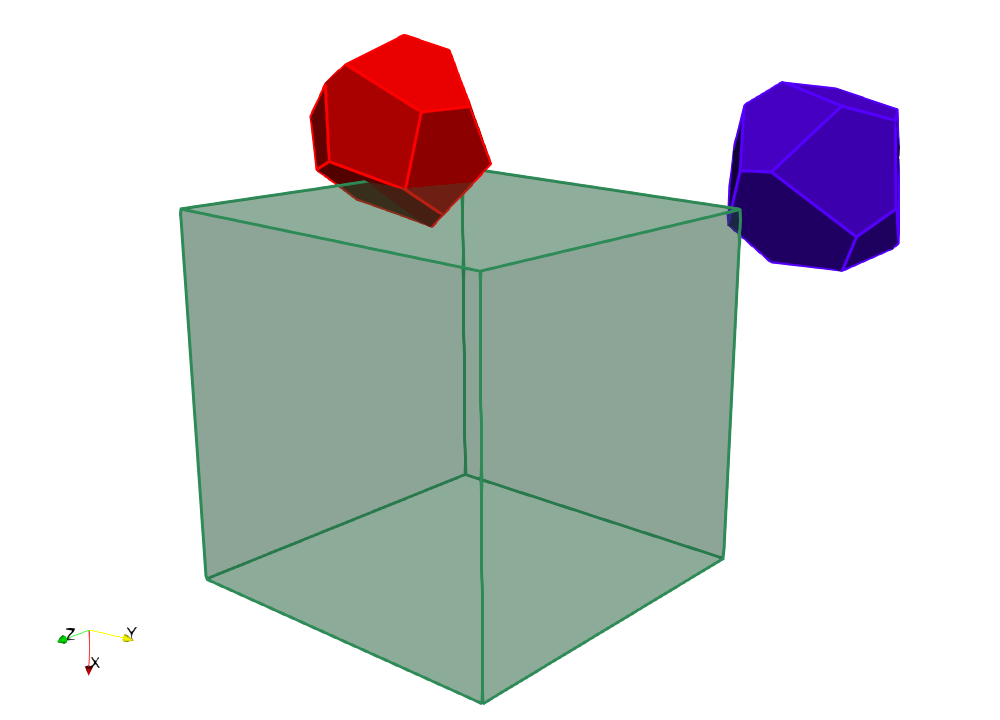}
	\end{minipage}\hfill
\begin{minipage}{0.7\textwidth}
	\caption{\textit{Separating Axis Theorem selection}. At this stage, we check if the fluid elements in $\mathcal{K}_\textrm{a} \setminus \mathcal{C}_\textrm{a}$ are intersecting with $ K_{\textrm{a}}$. We see that the red element (\textcolor{red}{$\bullet$})  is intersecting, and hence it is pushed in $\mathcal{I}$. The blue element (\textcolor{blue}{$\bullet$}) is not intersecting, so it is discarded. \label{fig:SAT}}
\end{minipage}
\end{subfigure}	
\begin{subfigure}{\textwidth}
	\begin{minipage}{0.25\textwidth}
	\includegraphics[width=\textwidth]{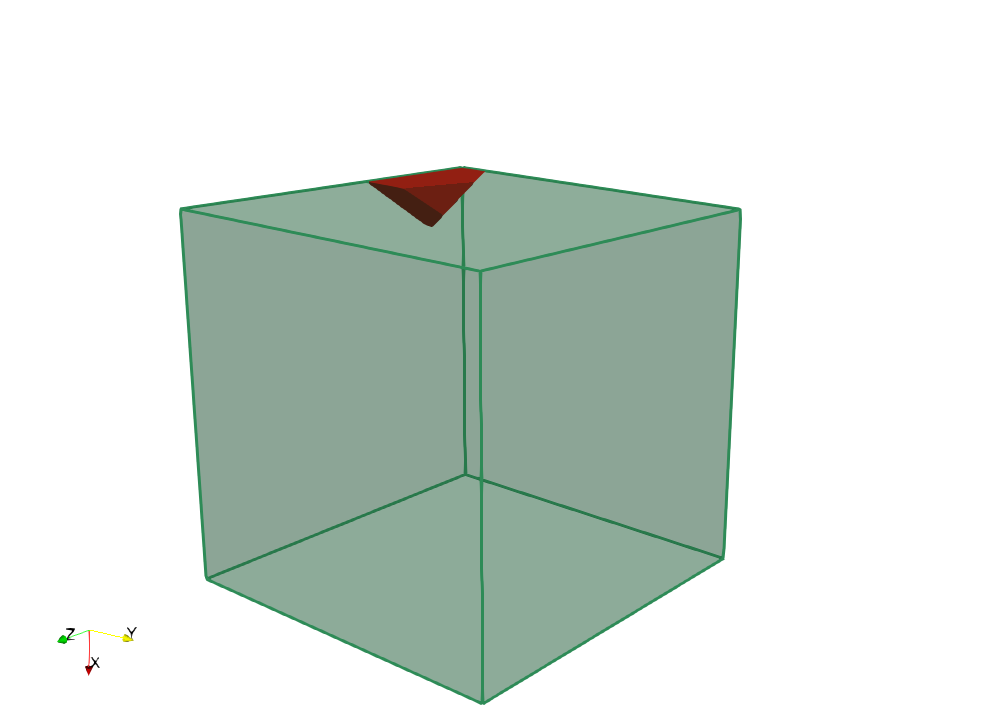}
	\end{minipage} \hfill
	\begin{minipage}{0.7\textwidth}
\caption{\textit{Intersection computation}. The intersection is explicitly computed only for the red element (\textcolor{red}{$\bullet$}). A new object is obtained, which is then stored and used for the projection computation. \label{fig:CGALinter}}
	\end{minipage}
\end{subfigure}	
\caption{Schematic representation of Algorithm~\ref{alg:cutCell}. 
(\textcolor{myGreen}{$\bullet$}) Acoustic element $ K_{\textrm{a}}$. (\textcolor{red}{$\bullet$}) Example of a fluid element $ K_{\textrm{f}} \in \mathcal{I}_\textrm{a}$. (\textcolor{blue}{$\bullet$}) Example of a fluid element $ K_{\textrm{f}} \in \mathcal{K}_\textrm{a}$ but not intersecting with the selected element $ K_{\textrm{a}}$. (\textcolor{yellow}{$\bullet$}) Example of a fluid element $ K_{\textrm{f}} \in \mathcal{C}_\textrm{a}$.
}
\end{figure}
The benefits of computing explicitly the intersection elements when projecting have been already explored in the context of low-order finite elements for aeroacoustics on tetrahedral meshes, see for instance \cite{Schoder2021}. The new grid obtained after the intersection is nested both with respect to the fluid grid and with respect to the acoustic grid, hence we can apply the analysis of Section~\ref{sec:Analysis}. Here, we propose a geometrical algorithm for computing the intersections between generic polyhedral grids. 
In particular, we consider a polyhedral tessellation $\TauFHrom$ for the fluid domain $\Omega_{\textrm{f}}$, while a hexahedral tessellation $\TauAHrom$ for the acoustic domain $\Omega_{\textrm{a}}$, see for instance Figure~\ref{fig:domainExample}. This choice is inherited from the numerical scheme that we apply to the aeroacoustic problem (see Section~\ref{sec:disc_NS}-~\ref{sec:sem-wave}), even if  the proposed algorithm is valid for generic polyhedral grids. 
Depending on the characteristic wave-length of the problem and on the numerical schemes considered, we assume that the number of fluid elements $ K_{\textrm{f}}$ is greater than the number of acoustic ones $ K_{\textrm{a}}$ and that the elements are all convex polyhedra. 
Considering polyhedra elements allows us to use the Separating Axis Theorem (SAT) for detecting if two elements have non-empty intersection. The main idea of the SAT is that, if two elements have empty intersection, then there exists a plane that separates them. Only a few directions depending on the normals to the faces of the elements and the edge elements have to be checked. A detailed description of the SAT algorithm can be found in  \cite[Chapter 8]{Eberly2020}. Before presenting the algorithm for computing the intersections between $ K_{\textrm{a}}$ and $ K_{\textrm{f}}$ we introduce some definitions and notations. 

\begin{definition}(Cartesian Bounding Box) Given a polyhedral element $ K \subset \mathbb{R}^3$, we denote with $\mathbf{v}_i = (x_i,y_i,z_i)$ with $i=1,\dots,n_v$ the $n_v$ vertices of $ K$.
\st{We indicate with $\mathcal{B}( K)$ his Cartesian bounding box:} 	
	\begin{equation}
	\stkout{ 	\mathcal{B}( K) = \Pi_{\alpha \in \{x,y,z\}}[\alpha_{min},\alpha_{max}],  }
	\end{equation} 
\st{where $\displaystyle \alpha_{min} = \min_{i=1,\dots,n_v} \alpha_i$, $\displaystyle \alpha_{max} = \max_{i=1,\dots,n_v} \alpha_i$ and where $\alpha_i \in \{x_i,y_i,z_i\}$.}
\Ar{We indicate with $\mathcal{B}( K)$ his Cartesian bounding box: 	
	\begin{equation}
		\mathcal{B}( K) = [x_{min},x_{max}] \times [y_{min},y_{max}] \times [z_{min},z_{max}],
	\end{equation} 
	where $\displaystyle x_{min} = \min_{i=1,\dots,n_v} x_i$, $\displaystyle x_{max} = \max_{i=1,\dots,n_v} x_i$, $\displaystyle y_{min} = \min_{i=1,\dots,n_v} y_i$, $\displaystyle y_{max} = \max_{i=1,\dots,n_v} y_i$, $\displaystyle z_{min} = \min_{i=1,\dots,n_v} z_i$, $\displaystyle z_{max} = \max_{i=1,\dots,n_v} z_i$. }
\end{definition}

For any element $ K_{\textrm{a}} \in \TauAHrom$, we define:
\begin{itemize}
	\item the set $\mathcal{K}_\textrm{a}$ collecting all the fluid elements $ K_{\textrm{f}}$ whose bounding box $\mathcal{B}( K_{\textrm{f}})$ intersects the bounding box $\mathcal{B}( K_{\textrm{a}})$, i.e., $\mathcal{B}( K_{\textrm{f}})\cap \mathcal{B}( K_{\textrm{a}}) \neq \emptyset$;
	\item the set $\mathcal{C}_\textrm{a}$ collecting all the fluid elements $ K_{\textrm{f}}$ whose bounding box $\mathcal{B}( K_{\textrm{f}})$ is strictly contained inside $ K_{\textrm{a}}$, i.e. $\mathcal{B}( K_{\textrm{f}})\subset  K_{\textrm{a}}$;
	\item \st{the set $\mathcal{I}_\textrm{a}$ collecting all the fluid elements $ K_{\textrm{f}}$ that have to be explicitly intersected with $ K_{\textrm{a}}$}.
	\item \Ar{the set $\mathcal{I}_\textrm{a}$ collects all the remaining fluid elements intersecting with $ K_{\textrm{a}}$. The fluid elements $ K_{\textrm{f}}$ do intersect with $ K_{\textrm{a}}$, but their bounding box is not fully contained inside $ K_{\textrm{a}}$ and hence the intersection has to be computed explicitly}.
\end{itemize}
We remark that the cardinality of $\mathcal{K}_\textrm{a}$ is strictly greater than the cardinality of $\mathcal{C}_\textrm{a}\cup\mathcal{I}_\textrm{a}$.
\Ar{Furthermore, note that not all the elements intersecting with $ K_\textrm{a}$ are in $\mathcal{I}_\textrm{a}$, since part of them is contained in $\mathcal{C}_\textrm{a}$.}
 Algorithm \ref{alg:cutCell} computes the intersections between $\TauAHrom$ and $\TauFHrom$ proceeding as follows: for any element $ K_{\textrm{a}} \in \TauAHrom$,
\begin{itemize}
	\item[1.] \textit{Bounding Box search}: search over the intersecting bounding boxes of the fluid elements $\mathcal{B}( K_{\textrm{f}})$ and the bounding box $\mathcal{B}( K_{\textrm{a}})$ of the acoustic element $ K_{\textrm{a}}$. If the intersection is not empty, the element $ K_{\textrm{f}}$ is added to the set $\mathcal{K}_\textrm{a}$, see Figure~\ref{fig:BBOXsearch}.
	\item[2.] \textit{Bounding Box selection}: 
	map the vertices of $\mathcal{B}( K_{\textrm{f}})$ via a Newton-Raphson algorithm by employing the trilinear map $\bm{\theta}_{ K_{\textrm{a}}}$. If all the vertices are  inside the reference element $\widehat{ K}_A$, then $ K_{\textrm{f}}$ is added to $\mathcal{C}_\textrm{a}$, see Figure~\ref{fig:Richard}.
	\item[3.] \textit{Separating Axis Theorem selection}: apply the SAT collision detection algorithm in order to understand if the intersections have to be computed. In fact, there might be fluid elements in $\mathcal{K}_\textrm{a}$ that are not effectively intersecting $ K_{\textrm{a}}$, see for instance Figure~\ref{fig:SAT}. The intersecting elements are added to $\mathcal{I}_\textrm{a}$.
	\item[4.] \textit{Intersection computation}:  compute explicitly the intersection between the fluid elements in $\mathcal{I}_\textrm{a}$ and the acoustic element $ K_{\textrm{a}}$, cf. Figure~\ref{fig:CGALinter}.
\end{itemize}

We remark that Step 2. of the algorithm is justified by the assumption that fluid elements are much smaller then the acoustic ones. The final intersection (Step 4.) is computed by employing the Computational Geometry Algorithms Library (CGAL) \cite{CGAL}, and it is based on the Nef implementation \cite{CGAL-Nef2022} that allows performing Boolean operations between solids.

\begin{algorithm}
	\begin{algorithmic}[1]
    \Statex{ $[\TauAHrom \cap \TauFHrom]$ =  \texttt{compute intersection}$(\TauAHrom, \TauFHrom)$} 
		\For{ $ K_{\textrm{a}} \in \TauAHrom$}
		\State Compute $\mathcal{B}( K_{\textrm{a}})$.
		\For{ $ K_{\textrm{f}} \in \TauFHrom$}
		\State Compute $\mathcal{B}( K_{\textrm{f}})$.
		\If{ $\mathcal{B}( K_{\textrm{a}}) \cap \mathcal{B}( K_{\textrm{f}})$}
			\State Add $ K_{\textrm{f}}$ in $\mathcal{K}_\textrm{a}$.
		\If{	$\mathcal{B}( K_{\textrm{f}})  \subset  K_{\textrm{a}} $ }
		\State $ K_{\textrm{f}} \cap  K_{\textrm{a}} =  K_{\textrm{f}}$ and $ K_{\textrm{f}} \in \mathcal{C}_\textrm{a}$.
		\State \Ar{Add $ K_{\textrm{f}}$ to  $\TauAHrom \cap \TauFHrom$.}
		\State \Ar{\textbf{break}}
		\EndIf
		\If{$ K_{\textrm{a}} \cap  K_{\textrm{f}}$}
		\State Add $ K_{\textrm{f}}$ in $\mathcal{I}_\textrm{a}$.
		\Else
		\State Elements are not intersecting.
		\EndIf
		\EndIf 
		\EndFor
		
		\For{ ( $ K_{\textrm{f}} \in \mathcal{I}_\textrm{a}$ )}
			\State Compute intersection $( K_{\textrm{a}} \cap 	 K_{\textrm{f}})$ with CGAL.
			\State \Ar{Add $ K_{\textrm{f}} \cap  K_{\textrm{a}}$ to  $\TauAHrom \cap \TauFHrom$.}
		\EndFor
		\EndFor
	\end{algorithmic}
	\caption{Algorithm to compute the intersection $\TauAHrom \cap \TauFHrom$ between the grids $\TauAHrom$ and $\TauFHrom$. \label{alg:cutCell}}
\end{algorithm}

\subsection{A quadrature-free method for integral evaluation}\label{sec:quadrature_free}

In this section, we explain how to compute numerically the integrals defined on the right-hand side of \eqref{eq:ProjectionProblemDiscrete}. In the aeroacoustic solver this technique is used for
computing the right-hand side of  
\eqref{eq:SemiDiscreteSEM}. 
We remark that if in $V_{\textrm{a}}$ we consider only linear polynomials in each space direction, i.e., $r=1$,  and if the maps $\bm{\theta}_{ K_{\textrm{a}}}$ are linear for any $ K_{\textrm{a}} \in \TauAHrom$, then it is convenient to use a mid-point quadrature method. In this case,  \eqref{eq:AlgebraicProblemProjection} becomes
\begin{equation}\label{eq:projKaltenbacher}
\sum_{ K_{\textrm{a}} \in \TauAHrom} ( q_{\textrm{a}}, \phi_{\textrm{a},i})_{ K_{\textrm{a}}}  = \sum_{ K_{\textrm{a}} \in \TauAHrom} \sum_{\ell=1}^{N_\textrm{f}}  \widehat{q}_{\textrm{f},\ell} ( 1, \phi_{\textrm{a},i})_{ K_{\textrm{a}}\cap  K_{\textrm{f},\ell} } 
 =  \sum_{ K_{\textrm{a}} \in   \TauAHrom} \sum_{\ell=1}^{N_\textrm{f}} \widehat{q}_{\textrm{f},\ell} \phi_{\textrm{a},i}(\mathbf{x}_b) \lvert K_{\textrm{a}} \cap  K_{\textrm{f},\ell} \rvert,
\end{equation}
where $\mathbf{x}_b$ is the barycentre of the intersection element $ K_{\textrm{a}} \cap  K_{\textrm{f},\ell}$, and $\lvert K_{\textrm{a}} \cap  K_{\textrm{f},\ell}\rvert$ is the volume of the intersection. The cut-volume cell-based interpolation that was proposed in \cite{Schoder2021} can be interpreted exactly as this mid point quadrature projection method. In fact in the latter work, the projection is evaluated by computing the intersections between a tetrahedral acoustic grid and a tetrahedral fluid grid and then using a mid-point quadrature rule on the intersected elements. However, when considering higher-order polynomials, i.e., $r>1$ in $V_{\textrm{a}}$,  or generic trilinear maps, leads to inexact quadrature integration that deteriorates the quality of the projection, as it will be numerically assessed in Sec.~\ref{sec:ConvergenceResults}. For this reason, we look for a quadrature formula that is able to integrate high-order polynomials on generic polyhedral elements (intersection of fluid and acoustic elements). When integrating polynomials over a polyhedral domain, one of the most popular choices is to sub-tessellate the polyhedral domain and then apply therein a standard quadrature formula over the tetrahedral mesh. This is  in general computationally expensive. For that reason we employed a Laserre-like integration \cite{Lasserre2015}, that has already been successfully applied in the context of discontinuous Galerkin methods, see for instance \cite{Antonietti2018}. The employed quadrature formula is able to integrate exactly homogeneous functions over general polyhedra $ K$.
We report here for completeness the main feature of the quadrature method, and refer to \cite{Antonietti2018} for further details. 
Let the polyhedron $ K \subset \mathbb{R}^3$  be a closed polytope, whose boundary $\partial K$ is defined by $m$ faces $ F_i \in \mathbb{R}^2$, with $i=1,\dots,m$. To each face $ F_i$ we associate a normal vector $\mathbf{n}_i$. Also, each face $ F_i$ lies on a hyperplane $\mathcal{H}_i$, and hence to each face $ F_i$ we associate a scalar $b_i$ such that $\forall \mathbf{x} \in \mathcal{H}_i$ we have that $\mathbf{n}_i\cdot \mathbf{x} = b_i$. 
Moreover, we split the polyhedron boundary as the union of $m$ faces, i.e., $\displaystyle \partial  K = \bigcup_i^m  F_i$, and the boundary of each face $ F_i$ as the union of $m_i$ edges, i.e.,  $\displaystyle \partial  F_i = \bigcup_{j}^{m_i}  F_{ij}$. Finally, the $m_{ij}$ vertices of each edge $ F_{ij}$ are denoted by  $\displaystyle \partial  F_{ij} = \bigcup_{k}^{m_{ij}}  F_{ijk}$.
Let $g$ to be homogeneous of degree $q>0$, namely, 
\begin{equation}\label{def:homo-function}
	q g(\mathbf{x}) = \nabla g(\mathbf{x}) \cdot \mathbf{x} \quad \forall \mathbf{x} \in \mathcal{ K},
\end{equation} 
and recall the generalized Stokes' theorem, see  \cite{Taylor1996}:
\begin{equation}
	\int_{ K} (\nabla\cdot \mathbf{V}(\mathbf{x}))g(\mathbf{x}) \text{d}\mathbf{x} + \int_{ K} \nabla g(\mathbf{x})\cdot \mathbf{V}(\mathbf{x}) \text{d}\mathbf{x}= \int_{\partial  K} \mathbf{V}(\mathbf{x})\cdot\mathbf{n}g(\mathbf{x}) d\sigma,
\end{equation}
where $\mathbf{V}: K \rightarrow \mathbb{R}^3$ is a generic vector field. By selecting $\mathbf{V}(\mathbf{x}) = \mathbf{x}$, and by applying \eqref{def:homo-function} we have 
\begin{equation} \label{eq:StokesUno}
	\int_{ K} g(\mathbf{x}) \text{d}\mathbf{x} = \frac{1}{3+q}\int_{\partial K} \mathbf{x}\cdot\mathbf{n}g(\mathbf{x}) d\sigma = \frac{1}{3+q}\sum_{i=1}^{m} b_i \int_{ F_i} g(\mathbf{x}) d\sigma.
\end{equation}
Next, by applying recursively integration by parts on \eqref{eq:StokesUno}, we obtain the following quadrature formula for computing the integral of a homogeneous function over a polyhedron $ K$:
\begin{equation} \label{eq:FormulaQuadratura}
	\begin{aligned}
		\int_{ K} g(\mathbf{x}) \text{d}\mathbf{x} = &\frac{1}{q+3} \sum_{i=1}^{m} \frac{b_i}{2+q} \left( \sum_{j=1}^{m_i} d_{ij} \int_{ F_{ij}} g(\mathbf{x}) d\nu + \int_{ F_i} \mathbf{x}_{0,i}\cdot\nabla g(\mathbf{x}) d\sigma \right), \\
		& \int_{ F_{ij}} g(\mathbf{x})d\nu =  \frac{1}{1+q} \left(\sum_{k=1}^{m_{ij}} d_{ijk} \int_{ F_{ijk}} g(\mathbf{x}) d\xi + \int_{ F_{ij}} \mathbf{x}_{0,ij}\cdot\nabla g(\mathbf{x}) d\nu\right),
	\end{aligned}
\end{equation}
where $d_{ij}$ is the Euclidean distance between the arbitrary point $\mathbf{x}_{0,i}$ and the edge $ F_{ij}$  and $d_{ijk}$ is the Euclidean distance between the arbitrary point $\mathbf{x}_{0,ij} \in  F_{ij}$ and the vertex $ F_{ijk}$.
We now apply the quadrature free rule described by 
Equation \eqref{eq:FormulaQuadratura} to \eqref{eq:ProjectionProblemDiscrete}. \st{since the integrated function is a polynomial, namely, it is a homogeneous function of degree $r$  }. \Ar{Each polynomial can be seen as sum of monomials, that are homogeneous function in the sense of eq.~\eqref{def:homo-function}. }
\st{Moreover, since we are employing spectral element methods, we usually integrate over a family of monomials.} To speed up the whole algorithm, the integrated monomials over $ K$ are stored and reused upon need. For further details on the implementation, we refer to Algorithm 2 in \cite{Antonietti2018}.

\section{Computational aspect of the  intersection algorithm}\label{sec:IntersectionResults}

In this section, we investigate some computational aspects of the algorithm presented 
 in Section~\ref{sec:IntersectionAlgorithm}. First, we verify the intersection algorithm in terms of accuracy and scalability.
Then, we use Algorithm~\ref{alg:cutCell} together with the quadrature-free method in Section~\ref{sec:quadrature_free} to compute integrals of polynomials over the domain $\Omega$. 

\begin{figure}
\begin{subfigure}{0.3\textwidth}
    \begin{overpic}[width=\textwidth]{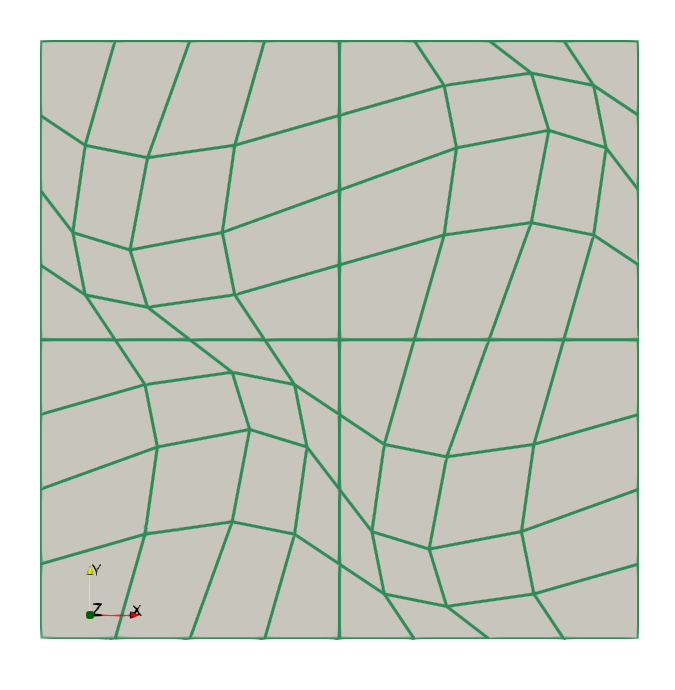}  
    \put(80,97){$\Large \TauAHrom^1$}
    \end{overpic}
\end{subfigure}
\begin{subfigure}{0.3\textwidth}
    \begin{overpic}[width=\textwidth]{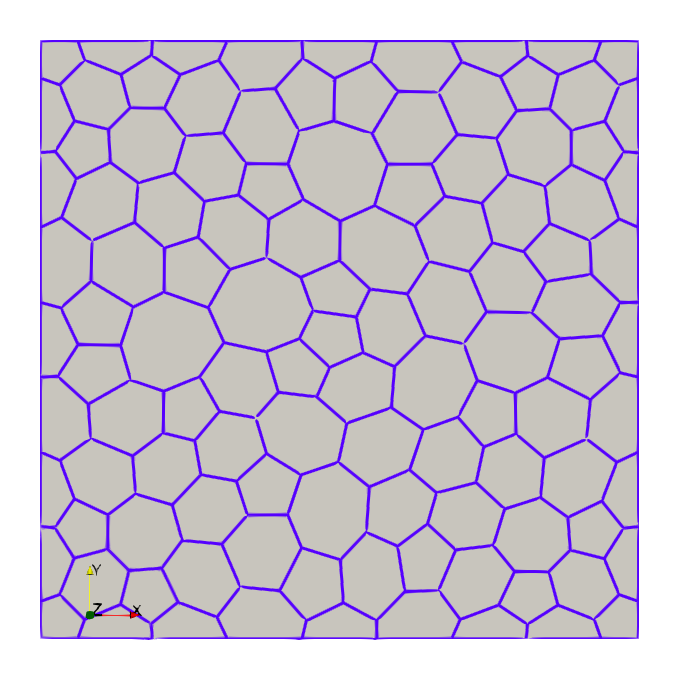}
        \put(80,97){$\Large \TauFHrom^1$}
    \end{overpic}
\end{subfigure} \hfill
\begin{subfigure}{0.35\textwidth} 
    \begin{overpic}[width=\textwidth]{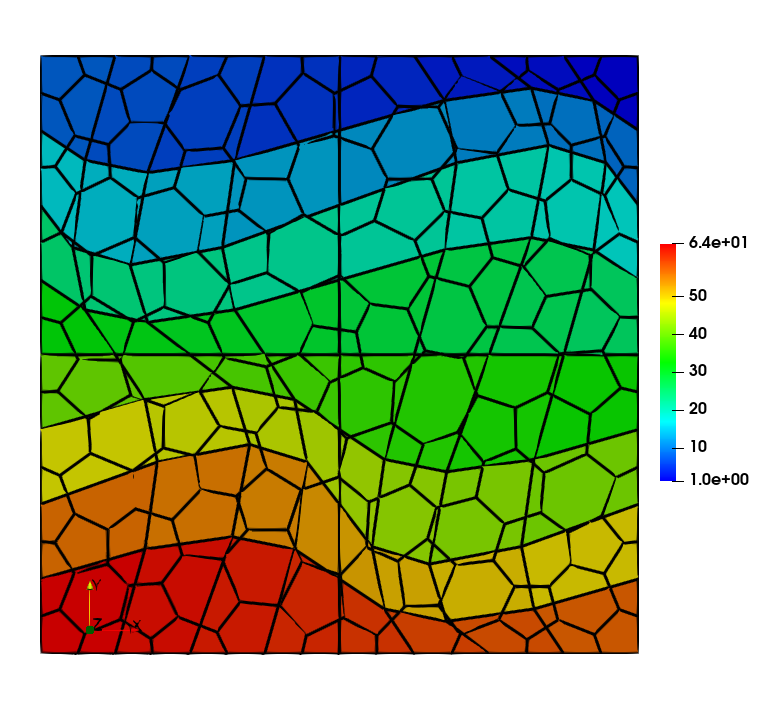}
    \put(85,15){$\Large \TauAHrom^1 \cap \TauFHrom^1$}
    \end{overpic}
    \end{subfigure}
    \caption{Two-dimensional view of Intersection between $\TauAHrom^1$ (left) and $\TauFHrom^1$ (center). The considered computational grids are first generated in two dimensions and then extruded, with only one element in the vertical direction. All the intersections between elements in $\TauFHrom^1$ and a single element in $\TauAHrom^1$ have the same color.}
    \label{fig:Test2D}
\end{figure}

\begin{figure}
\begin{subfigure}{0.3\textwidth}
    \begin{overpic}[width=\textwidth]{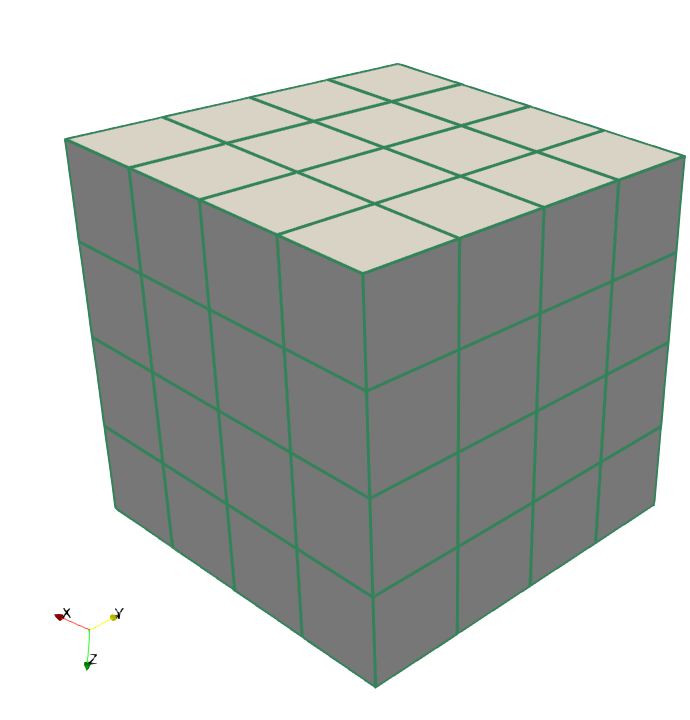}    
     \put(80,15){$\Large \TauAHrom^2$}
    \end{overpic}
\end{subfigure}
\begin{subfigure}{0.3\textwidth}
    \begin{overpic}[width=\textwidth]{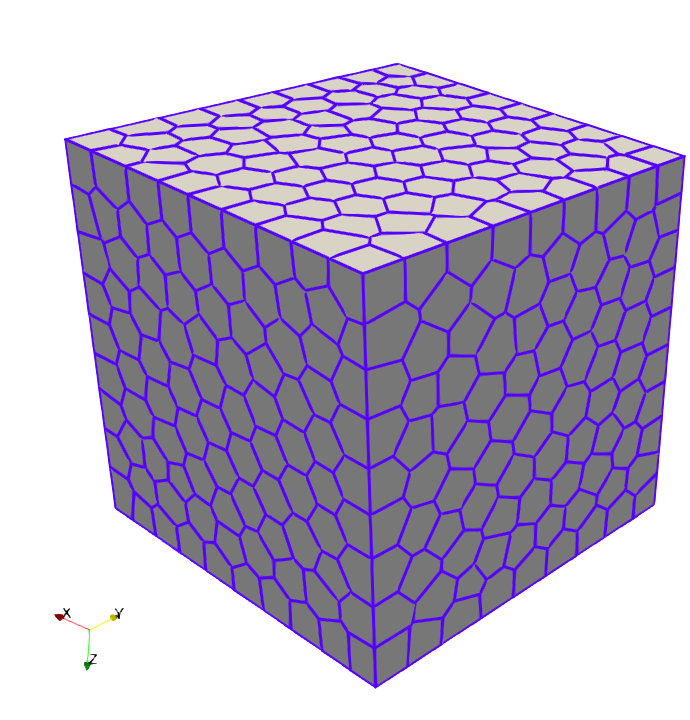}
     \put(80,15){$\Large \TauFHrom^2$}
    \end{overpic}
\end{subfigure} \hfill
\begin{subfigure}{0.325\textwidth} 
    \begin{overpic}[width=\textwidth]{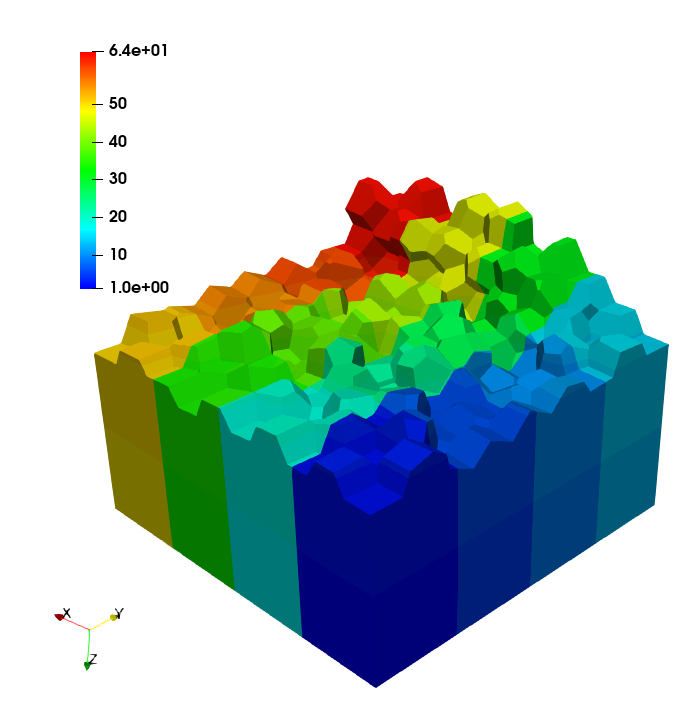}
    \put(80,15){$\Large \TauAHrom^2 \cap \TauFHrom^2$}
    \end{overpic}
    \end{subfigure}
    \caption{Intersection between $\TauAHrom^2$ (left) and $\TauFHrom^2$ (center). The fluid grid is made of polyhedral elements, while the acoustic grid is made of hexahedral elements. All the intersections between elements in $\TauFHrom^2$ and a single element in $\TauAHrom^2$ have the same color.}
    \label{fig:Test3D}
\end{figure}

To check the accuracy of the proposed intersection algorithm we consider the following mesh configurations.
In the first test, we set $\Omega = \Omega_{\textrm{f}} = \Omega_{\textrm{a}} = (-2,2)\times(-2,2)\times(-0.05,0.05)$ and define 
the acoustic grid $\TauAHrom^1$ (resp. fluid grid $\TauFHrom^1$) by extruding in the vertical direction distorted quadrilaterals (resp. polygons), cf. Figure~ \ref{fig:Test2D}. The acoustic grid has 64 elements and the original Cartesian mesh size was $h_{\textrm{a}}=0.5$, while the fluid grid has 109 elements and  $h_{\textrm{f}}=0.5$. 
In the second test, we consider $\Omega = \Omega_{\textrm{f}} = \Omega_{\textrm{a}} = (-0.5,0.5)^3$ and use a Cartesian grid $\TauAHrom^2$ with 64 elements and $h_{\textrm{a}} = 0.25$ in $\Omega_{\textrm{a}}$, while a Voronoi polyhedral grid $\TauFHrom^2$ with 1000 elements and $h_{\textrm{f}} = 0.1$ in $\Omega_{\textrm{f}}$, see Figure~\ref{fig:Test3D}. The computed intersections are shown in Figures~\ref{fig:Test2D} and \ref{fig:Test3D} (right).  
To have a quality check of the performed algorithm we color the resulting grid  $\TauAHrom^1 \cap  \TauFHrom^1$ in the following way. All the intersections between elements in $\TauFHrom^1$ and a single element in $\TauAHrom^1$ have the same color. 
It is possible to notice that even small elements are intersected properly by the proposed intersection algorithm, cf. Figures~\ref{fig:Test2D} and \ref{fig:Test3D} (right). 
To show the accuracy of the proposed intersection algorithm, we consider the following verification test that computes the integral of polynomial functions over the intersection grid $\TauAHrom\cap\TauFHrom$. In Tables~\ref{tab:checkIntersezioni_1} and  ~\ref{tab:checkIntersezioni_2} we report the relative errors \begin{equation} E_{rel}(f) = \dfrac{\left|\left(\int_{\Omega} f d\mathbf{x} - \int_{\mathcal{T}} f d\mathbf{x} \right)\right|}{ \left|\left(\int_{\Omega} f d\mathbf{x} \right)\right|}, \label{eq:ErrorIntegration} \end{equation}  computed  by employing the quadrature free method in Section~\ref{sec:quadrature_free}, by varying the mesh $\mathcal{T}$ of the domain $\Omega$. Here, $\phi$ is a generic monomial function.
Since the quadrature-free algorithm is exact for homogeneous functions, from the results it is possible to conclude that the intersection computation does not introduce any additional error.



\begin{table}[]
	\centering
	\bgroup
	\def\arraystretch{1.2}
	\begin{tabular}{ccccc}
		mesh  & $ E_{rel}(1)$ & $E_{rel}(x^2y^2) $ & $ E_{rel}(x^4y^4) $ \\
		\hline
		$\TauAHrom^1$ & $ 5.551\times10^{-16}$ & $  6.661\times10^{-16}$ & $8.882\times10^{-16}$ \\
		$\TauFHrom^1$  & $5.551\times10^{-16}$ & $4.441\times10^{-16}$ & $4.441\times10^{-16}$ \\
		$\TauAHrom^1\cap\TauFHrom^1$ & $5.551\times10^{-16} $ & $1.11\times10^{-16}$ & $2.22\times10^{-16}$ & \\
	\end{tabular}
	\egroup
	
	\caption{Computed \Ar{$E_{rel}(f)$, see eq.~\eqref{eq:ErrorIntegration}} for different meshes: acoustic grid $\TauAHrom^1$,  fluid grid $\TauFHrom^1$ and their intersection $\TauAHrom^1\cap\TauFHrom^1$. Here,  $\Omega = (-2,2)\times(-2,2)\times(-0.05,0.05)$.}
	\label{tab:checkIntersezioni_1}
\end{table}

\begin{table}[]
\centering
\bgroup
\def\arraystretch{1.2}
    \begin{tabular}{ccccc}
    mesh  & $E_{rel}(1)$ & $E_{rel}(x^2y^2z^2) $ & $E_{rel}(x^4y^4z^4) $ \\
    \hline
      $ \TauAHrom^2$ & $0$ & $0$ & $4.441\times10^{-16}$ \\
       $\TauFHrom^2$  & $4.441\times10^{-16}$ & $2.22\times10^{-16}$ & $1.11\times10^{-15}$ \\
      $\TauAHrom^2\cap\TauFHrom^2$ & $0$ & $0$ & $0$ & \\
    \end{tabular}
    \egroup
    \caption{Computed error \Ar{$E_{rel}(f)$, see eq.~\eqref{eq:ErrorIntegration}}, for different meshes: acoustic grid $\TauAHrom^2$,  fluid grid $\TauFHrom^2$ and their intersection $\TauAHrom^2\cap\TauFHrom^2$. Here,  $\Omega = (-0.5,0.5)^3$.} 
    \label{tab:checkIntersezioni_2}
\end{table}



To assess the scalability of the algorithm we consider $\Omega = \Omega_{\textrm{f}} = \Omega_{\textrm{a}} = (-0.5,0.5)^3$. We tessellate the domain $\Omega_{\textrm{a}}$ (resp. $\Omega_{\textrm{f}}$) with a grid made by $32^3$ (resp. $65^3$) elements.
The total number of computed intersections is 884736, with $\sum_{ K_{\textrm{a}}\in\TauAHrom} \card(\mathcal{C}_\textrm{a}) = 39304$ and $\sum_{ K_{\textrm{a}}\in\TauAHrom} \card(\mathcal{I}_\textrm{a}) = 845432$, see Section~\ref{sec:IntersectionAlgorithm}, that is where most of the computational time is spent by the algorithm concerns the evaluation of actual intersections, i.e., lines 17-19 of Algorithm \ref{alg:cutCell}. 

We perform a strong scalability test on the G100 cluster located at Cineca, by keeping the same computational grids while varying the number of available cores. 
\begin{figure}
    \centering
    \includegraphics[width=0.5\textwidth]{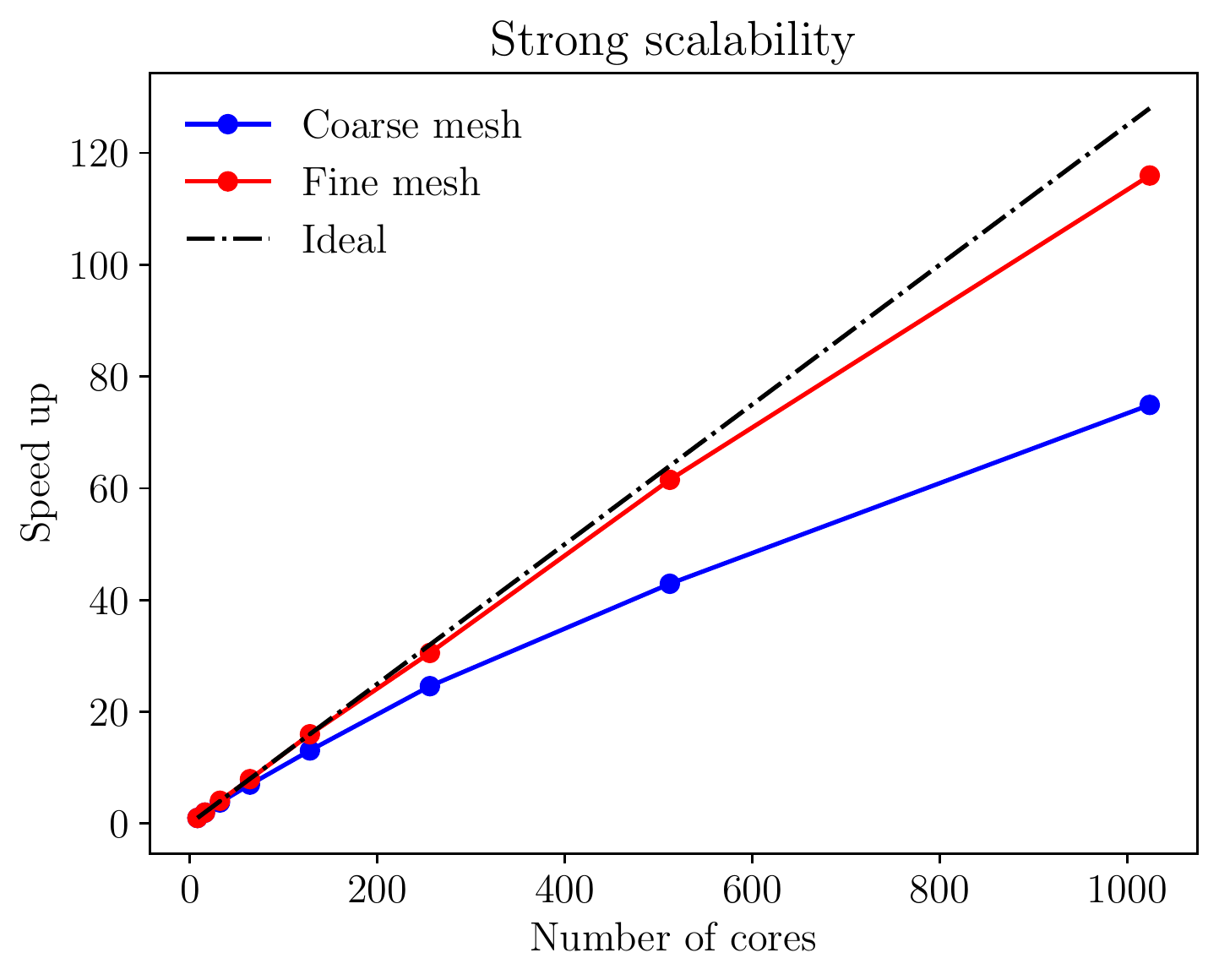}
    \caption{Scalability test. The speed-up is computed with respect to the test performed on 8 cores. The coarse mesh has 884736 intersections, while the fine mesh has 7077888 intersections.}
    \label{fig:scalabilita}
\end{figure}
From Figure~\ref{fig:scalabilita} it is possible to notice that the algorithm scales well up to 128 cores. Then, since the partitioning of the acoustic mesh 
is independent of the underlying fluid grid, the number of intersecting elements for larger decompositions might vary largely between the processors, leading to unbalance in the intersection computations. To verify this, we design a larger test where an acoustic Cartesian grid with $64^3$ elements and a fluid grid with $65^3$ elements are considered. The total number of computed intersections is $\sum_{ K_{\textrm{a}}\in\TauAHrom} \card(\mathcal{I}_\textrm{a}) = 7077888$. On this latter test, where more elements are employed, the balance of the intersection is good and the scalability is almost ideal.

\section{Convergence results for the $L^2$-projection method} \label{sec:ConvergenceResults}
In this section, we inquire about the convergence properties of the developed projection method. 
In particular, we verify the theoretical estimate in Theorem~\ref{theorem:approximationTheorem}  for the approximation error $\norm{f-f_\textrm{a}}{L^2(\Omega)}$ and compare our approach with the one presented  in \cite{Schoder2021}.

\subsection{Verification and validation test cases} 
\begin{figure}[h!]
	\begin{subfigure}{0.32\textwidth} 
		\includegraphics[width=\textwidth]{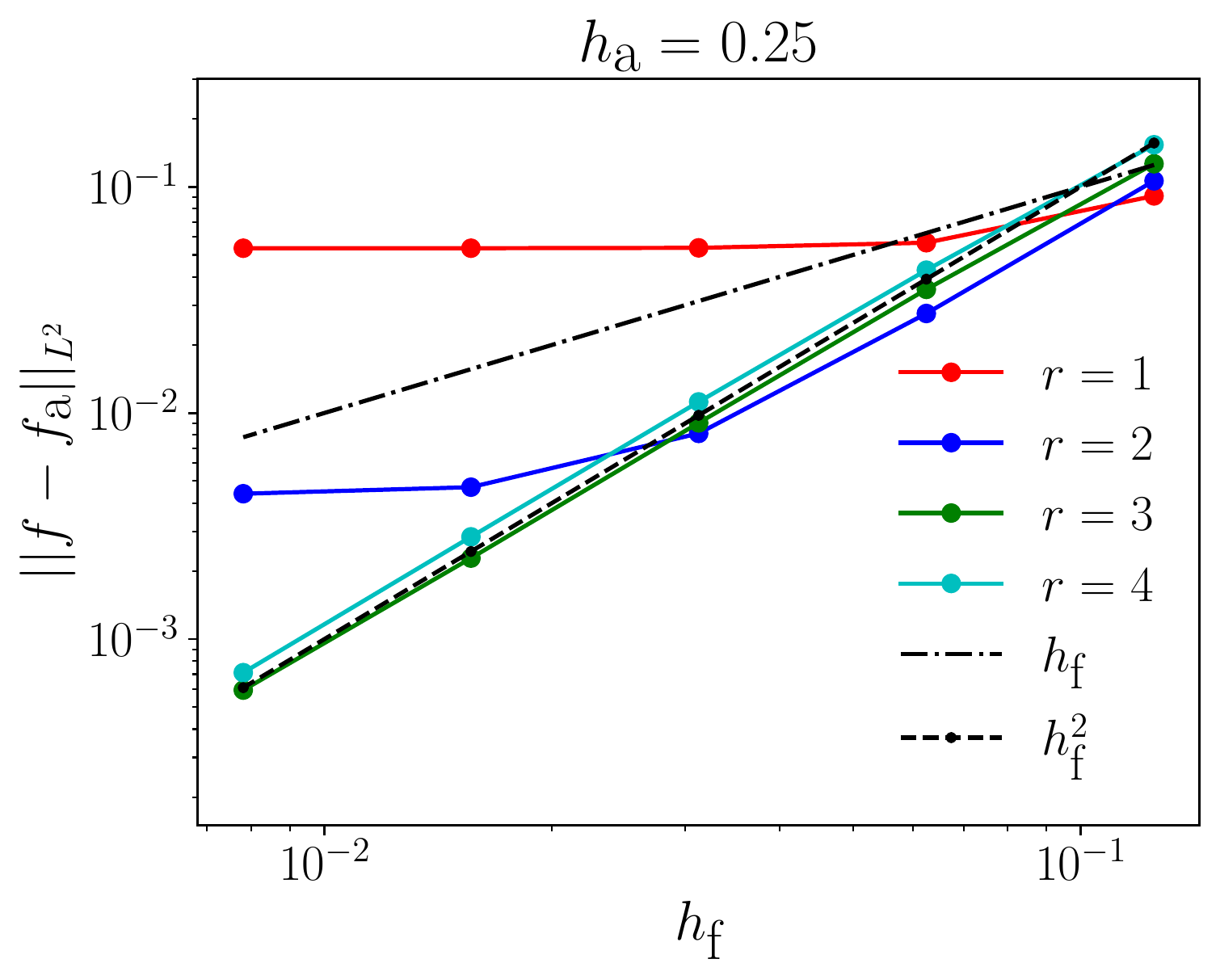}
	\end{subfigure}
	\begin{subfigure}{0.32\textwidth} 
		\includegraphics[width=\textwidth]{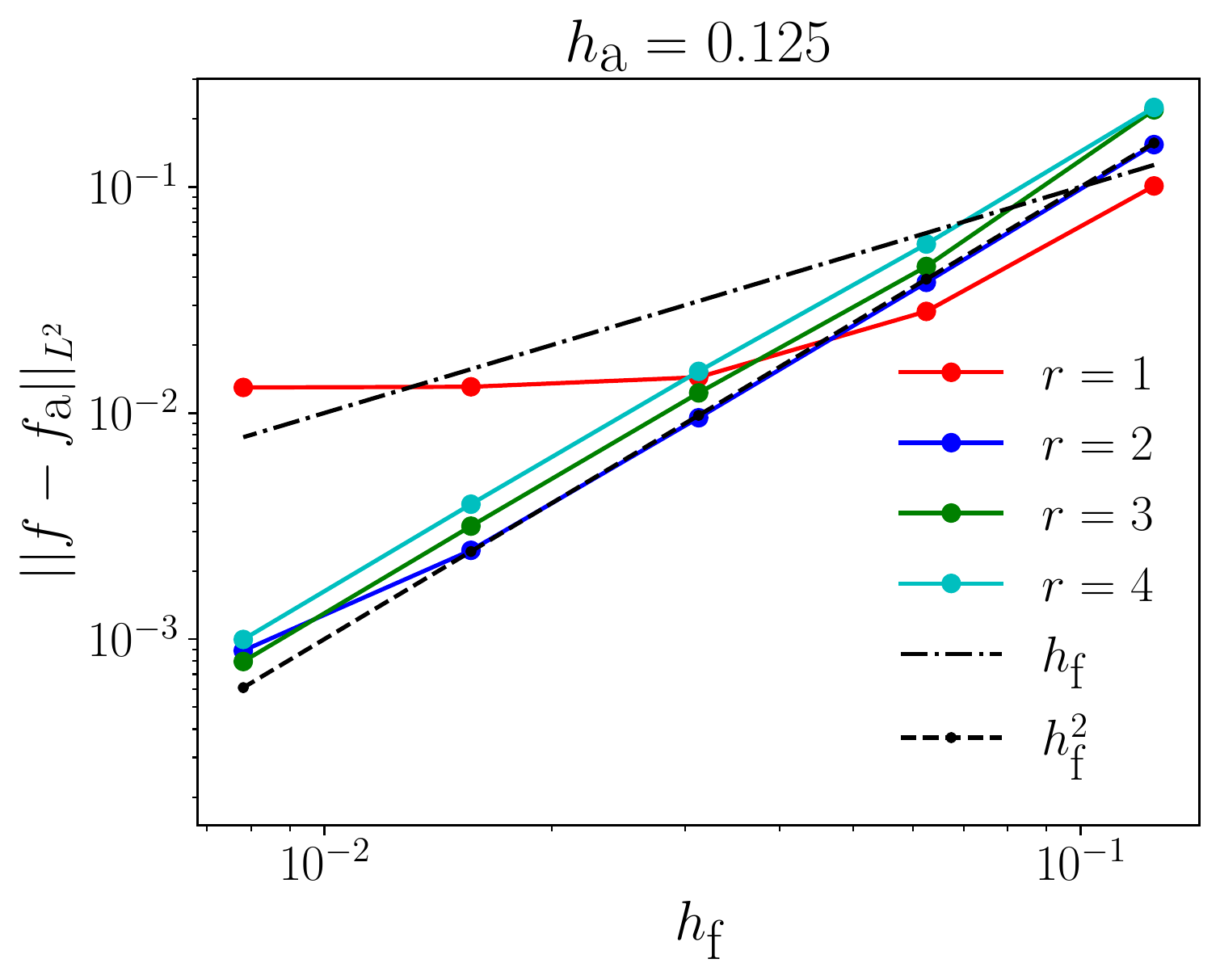}
	\end{subfigure}
	\begin{subfigure}{0.32\textwidth} 
		\includegraphics[width=\textwidth]{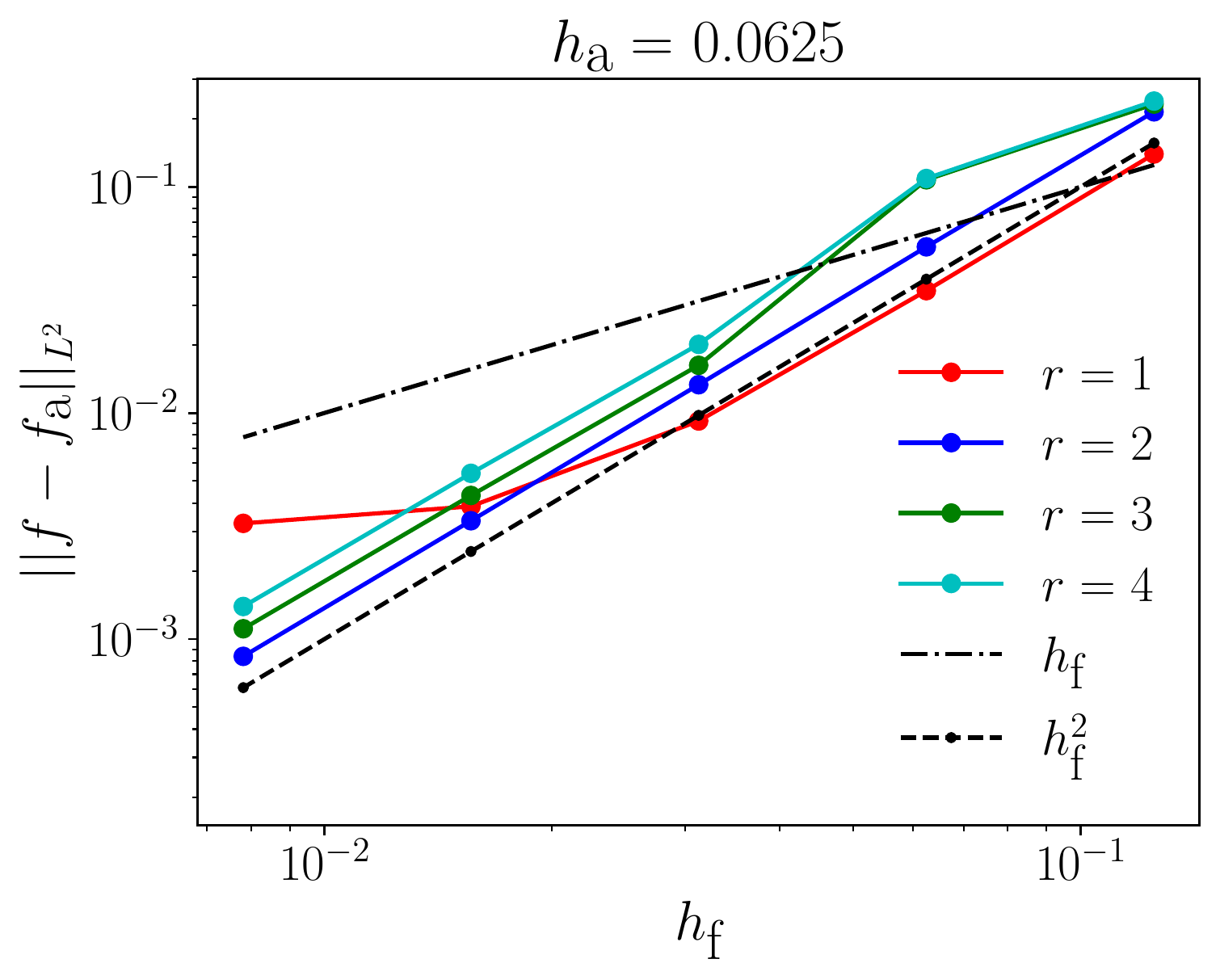}
	\end{subfigure}
	\caption{Computed errors $\norm{f - f_\textrm{a}}{L^2(\Omega)}$ versus $h_{\textrm{f}}$, for different polynomial degrees $r=1,2,3,4$ and different choices of $h_{\textrm{a}}=0.25,0.125,0.0625$ . \label{fig:ApproxF-FA-1}}
\end{figure}

\begin{figure}[h!]
	\begin{subfigure}{0.32\textwidth} 
		\includegraphics[width=\textwidth]{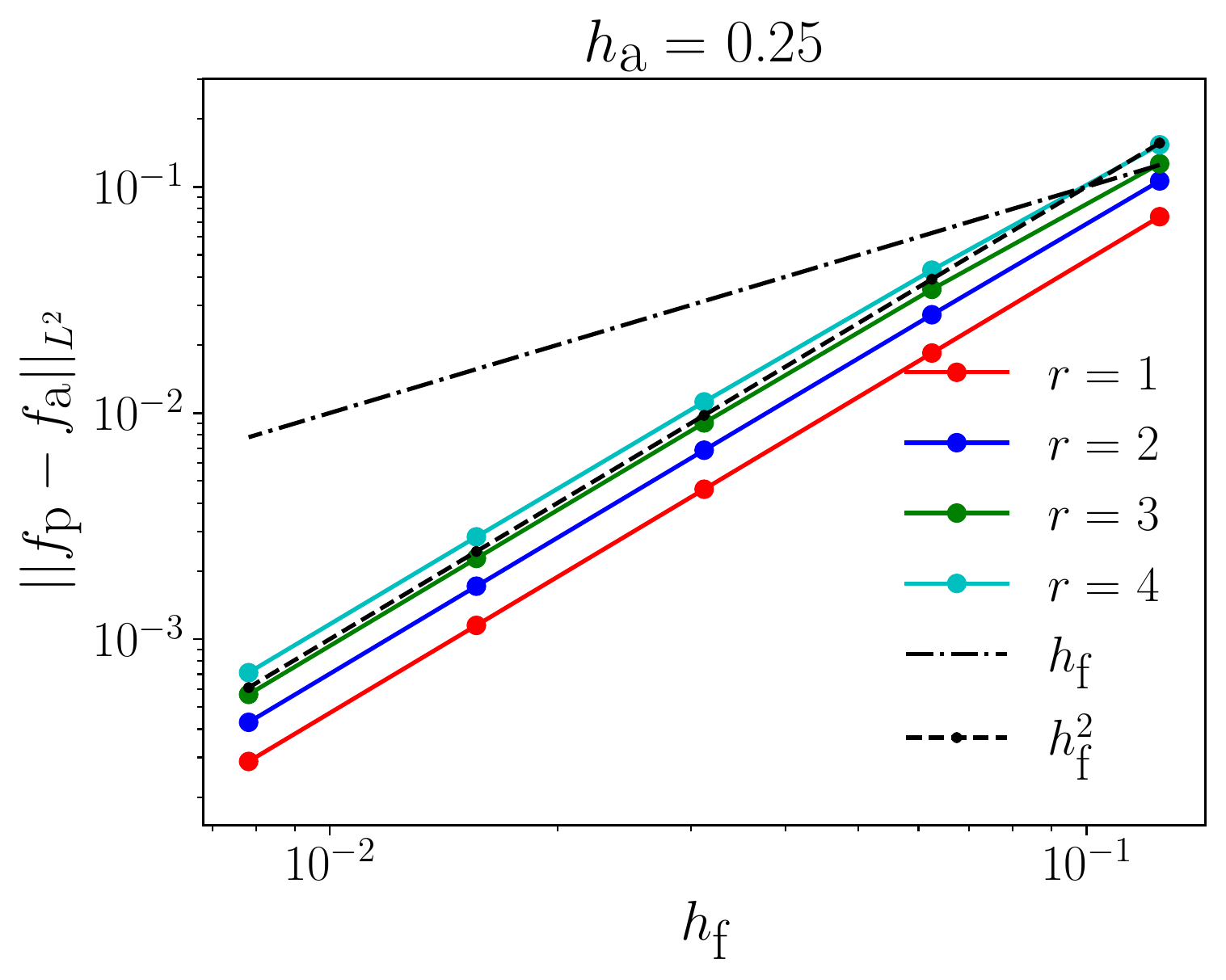}
	\end{subfigure}
	\begin{subfigure}{0.32\textwidth} 
		\includegraphics[width=\textwidth]{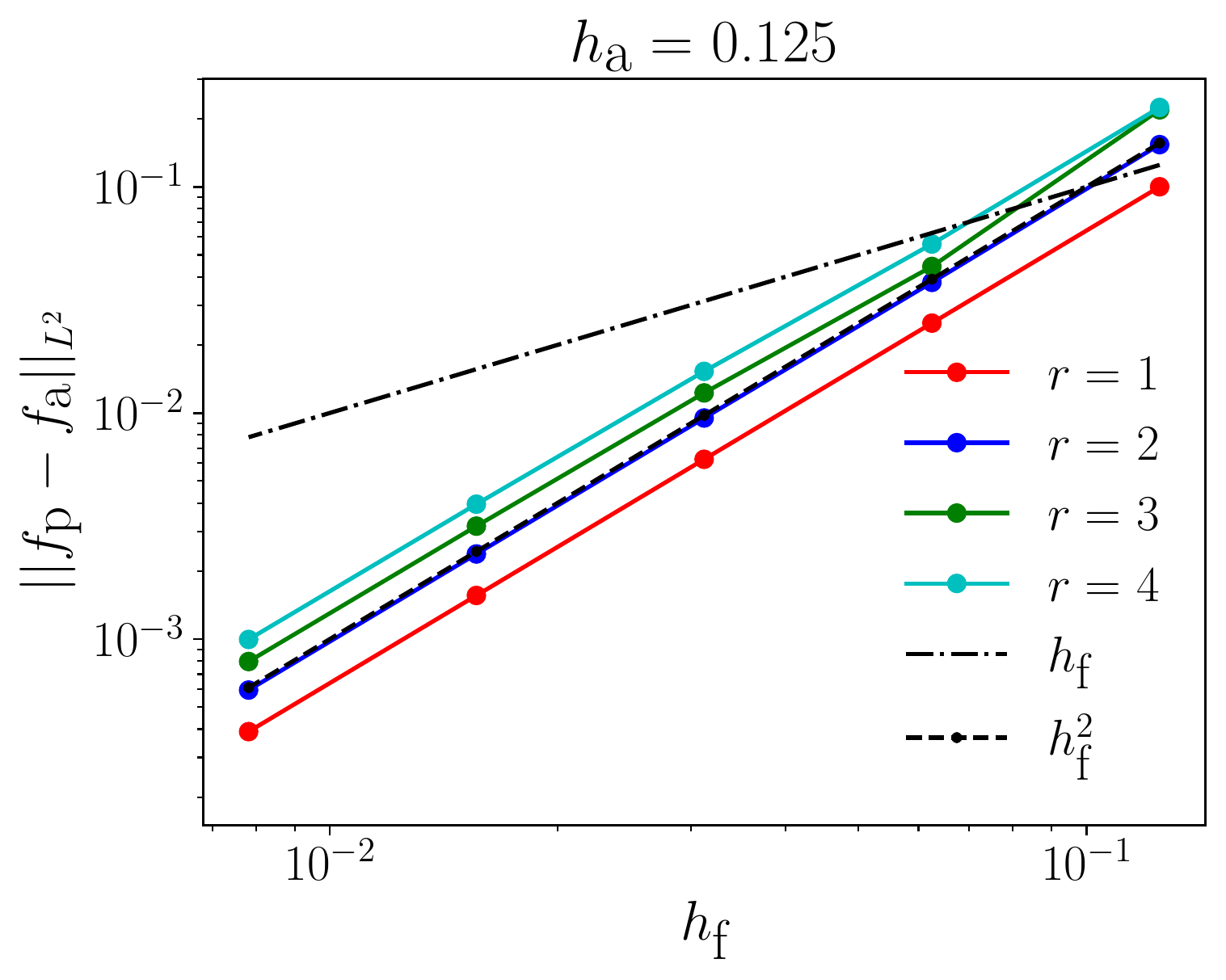}
	\end{subfigure}
	\begin{subfigure}{0.32\textwidth} 
		\includegraphics[width=\textwidth]{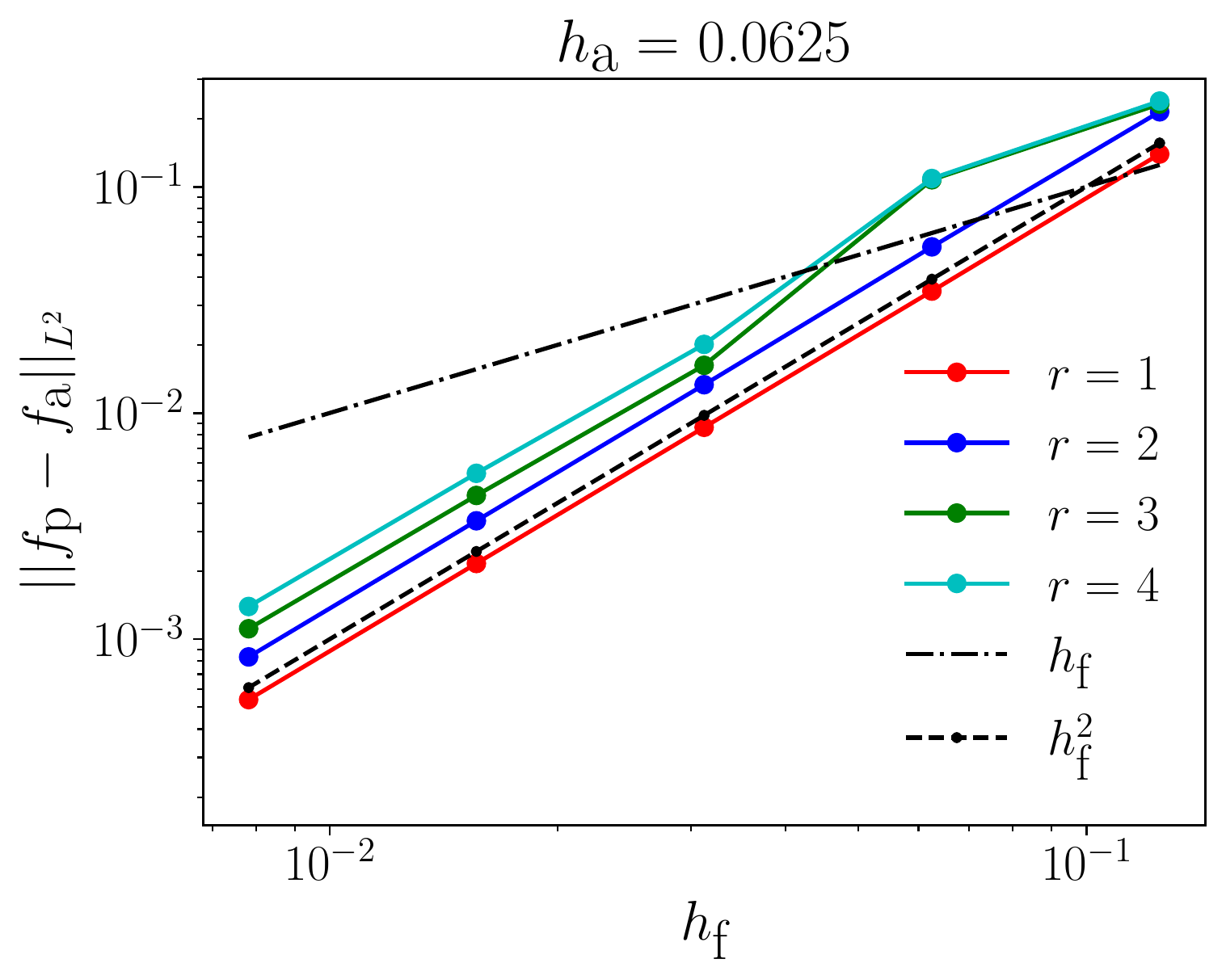}
	\end{subfigure}
	\caption{Computed errors $\norm{f_\textrm{p} - f_\textrm{a}}{L^2(\Omega)}$ versus $h_{\textrm{f}}$, for different polynomial degrees $r=1,2,3,4$ and different choices of $h_{\textrm{a}}=0.25,0.125,0.0625$ . \label{fig:ApproxFP-FA-1}}
\end{figure}

\begin{figure}[h!]
	\begin{center}
		\begin{subfigure}{0.45\textwidth} 
			\includegraphics[width=\textwidth]{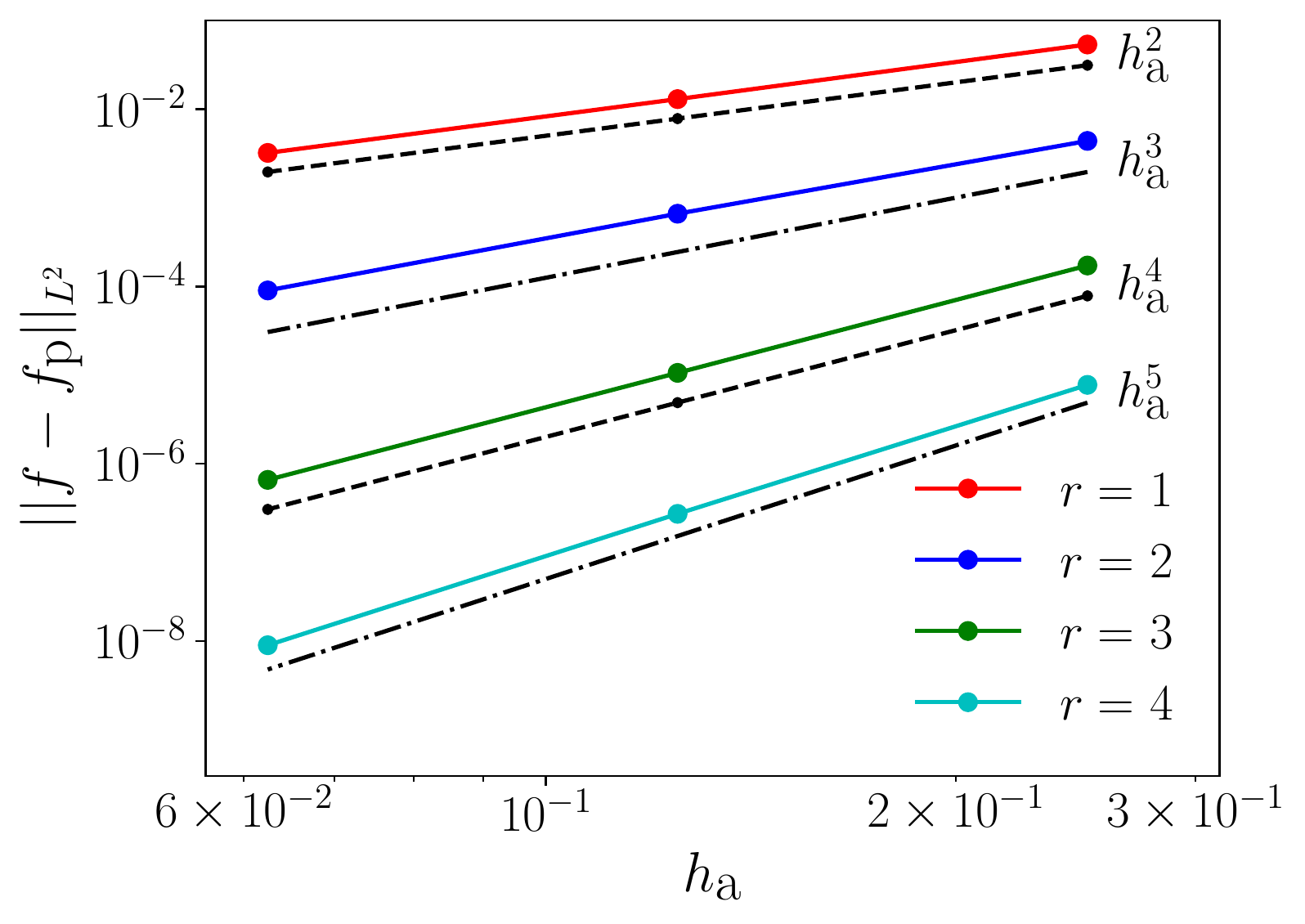}
		\end{subfigure}
		\begin{subfigure}{0.45\textwidth} 
			\includegraphics[width=\textwidth]{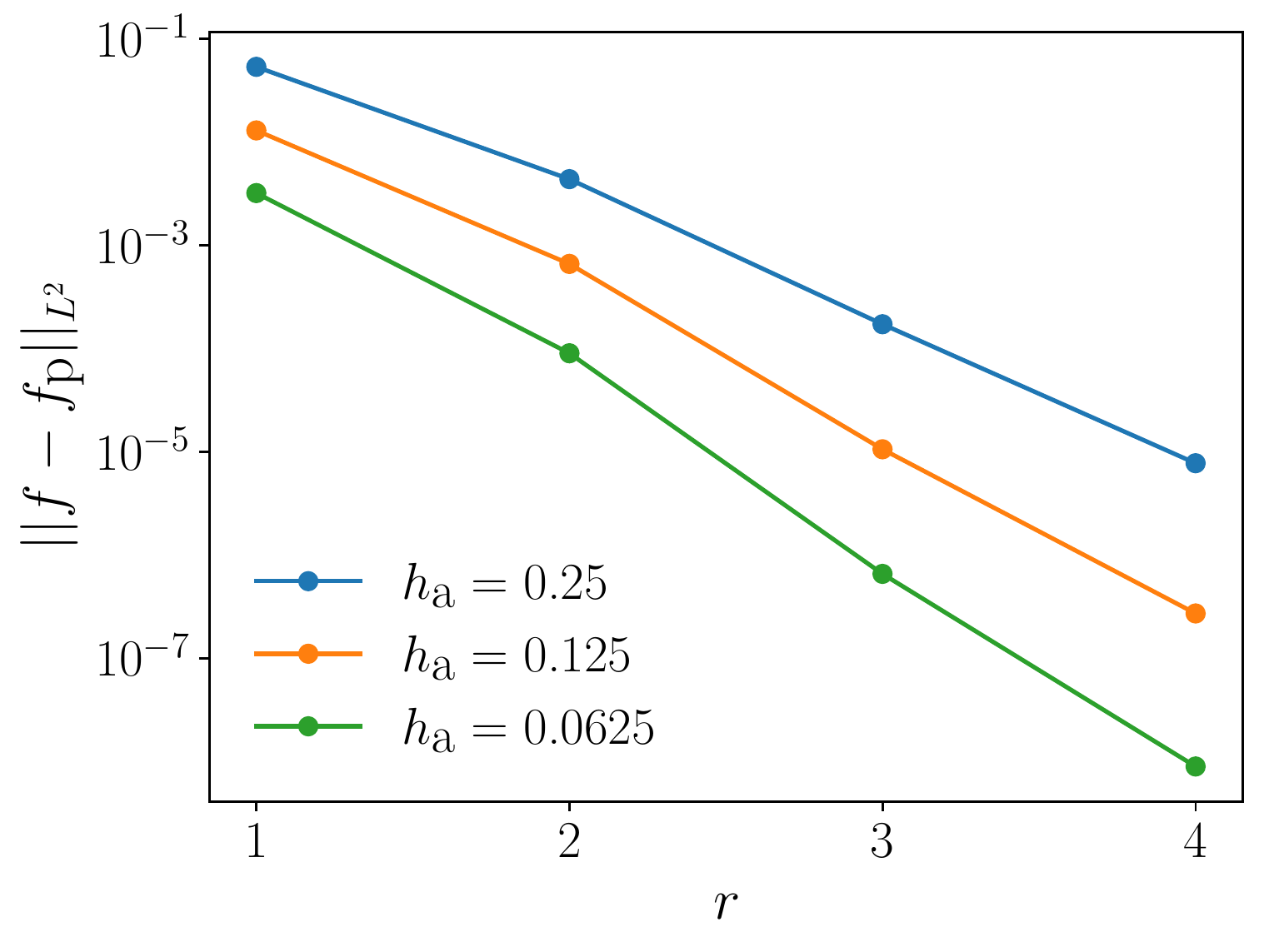}
		\end{subfigure}
		\caption{Computed errors $\norm{f - f_\textrm{p}}{L^2(\Omega)}$ versus $h_{\textrm{a}}$ (left) and $r$ (right), for different choices of  $r=1,2,3,4$ and $h_{\textrm{a}}=0.25,0.125,0.0625$ .   \label{fig:ApproxF-FP-1}}
	\end{center}
\end{figure}

We consider a cubic domain $\Omega = \Omega_{\textrm{a}} = \Omega_{\textrm{f}} = (-0.5,0.5)^3$ 
and two Cartesian nested tessellation $\TauAHrom$ and $\TauFHrom$, being the acoustic mesh size $h_{\textrm{a}}$ a multiple of fluid one $h_{\textrm{f}}$. 
Next, we consider $f = \cos(2\pi x)\cos(2\pi y)\cos(\pi z)$ and compute $E_\textrm{a} = \norm{f-f_\textrm{a}}{L^2(\Omega)}$, where $f_\textrm{a}$ is the projection defined as in Figure~\ref{fig:notationApprox} and computed as described in Eq.~\ref{eq:ProjectionProblemDiscrete} employing the quadrature method discussed in Sec.~\ref{sec:quadrature_free}. 
In Figure~\ref{fig:ApproxF-FA-1} we report the projection error $E_\textrm{a}$, by varying $h_{\textrm{f}}$ for fixed values of $h_{\textrm{a}}$ and the polynomial degree $r$. 
It is clear that the error $E_\textrm{a}$ saturates as we refine $h_{\textrm{f}}$. Indeed, by triangle inequality, we observe that
\begin{equation}
    E_\textrm{a}  \le \norm{f-f_\textrm{p}}{L^2(\Omega)} + \norm{f_\textrm{p}-f_\textrm{a}}{L^2(\Omega)},
\end{equation}
and that $\norm{f-f_\textrm{p}}{L^2(\Omega)}$ is the leading term of the error independent of $h_{\textrm{f}}$.
This is confirmed by the plots reported in Figure~\ref{fig:ApproxFP-FA-1} where we show the trend of the error $\norm{f_\textrm{p}-f_\textrm{a}}{L^2(\Omega)}$ as a function of 
$h_{\textrm{f}}$. The latter is proportional to $ h_{\textrm{f}}^2$
as predicted by \eqref{eq:fp-fa}.
On the other hand,  the error $\norm{f-f_\textrm{p}}{L^2(\Omega)}$ remains constant, cf. \eqref{eq:ApplicationInterpolantGLL}.
Moreover, we notice that increasing the polynomial degree $r$, keeping fixed $h_{\textrm{a}}$, reduces the saturation value reached by the error $\norm{f-f_\textrm{a}}{L^2(\Omega)}$.
Finally, in Figure~\ref{fig:ApproxF-FP-1} (left), we plot the error $\norm{f-f_\textrm{p}}{L^2(\Omega)}$ versus the mesh size $h_{\textrm{a}}$ while in  Figure~\ref{fig:ApproxF-FP-1} (right) the same quantity is shown as a function of $r$. The expected convergence rate given by estimate \eqref{eq:ApplicationInterpolantGLL} is confirmed by the numerical results.


We provide the following rule of thumb to decide how to relate the acoustic and fluid grid in terms of mesh sizes $h_{\textrm{a}},h_{\textrm{f}}$, and polynomial degree $r$.
Lower projection errors would be generally obtained if both grids have a similar number of degrees of freedom.
As seen from estimate \eqref{eq:estimateApprox}, the approximation error $E_\textrm{a}$ is lower employing for the acoustic problem a low order polynomial degree and a spatial resolution comparable to the fluid grid, namely $h_{\textrm{f}} \approx h_{\textrm{a}}$. However, this choice deteriorates the convergence error estimates provided by the Strang Lemma for the SEM-NI method, see e.g., \cite[Lemma 10.1]{Quarteroni2009}. The numerical tests presented above show that the dependency on the polynomial degree is not so severe as stated in \cref{eq:estimateApprox}, encouraging the use of high-order basis functions.




%
\begin{figure}[h!]
	\begin{subfigure}{0.32\textwidth} 
		\includegraphics[width=\textwidth]{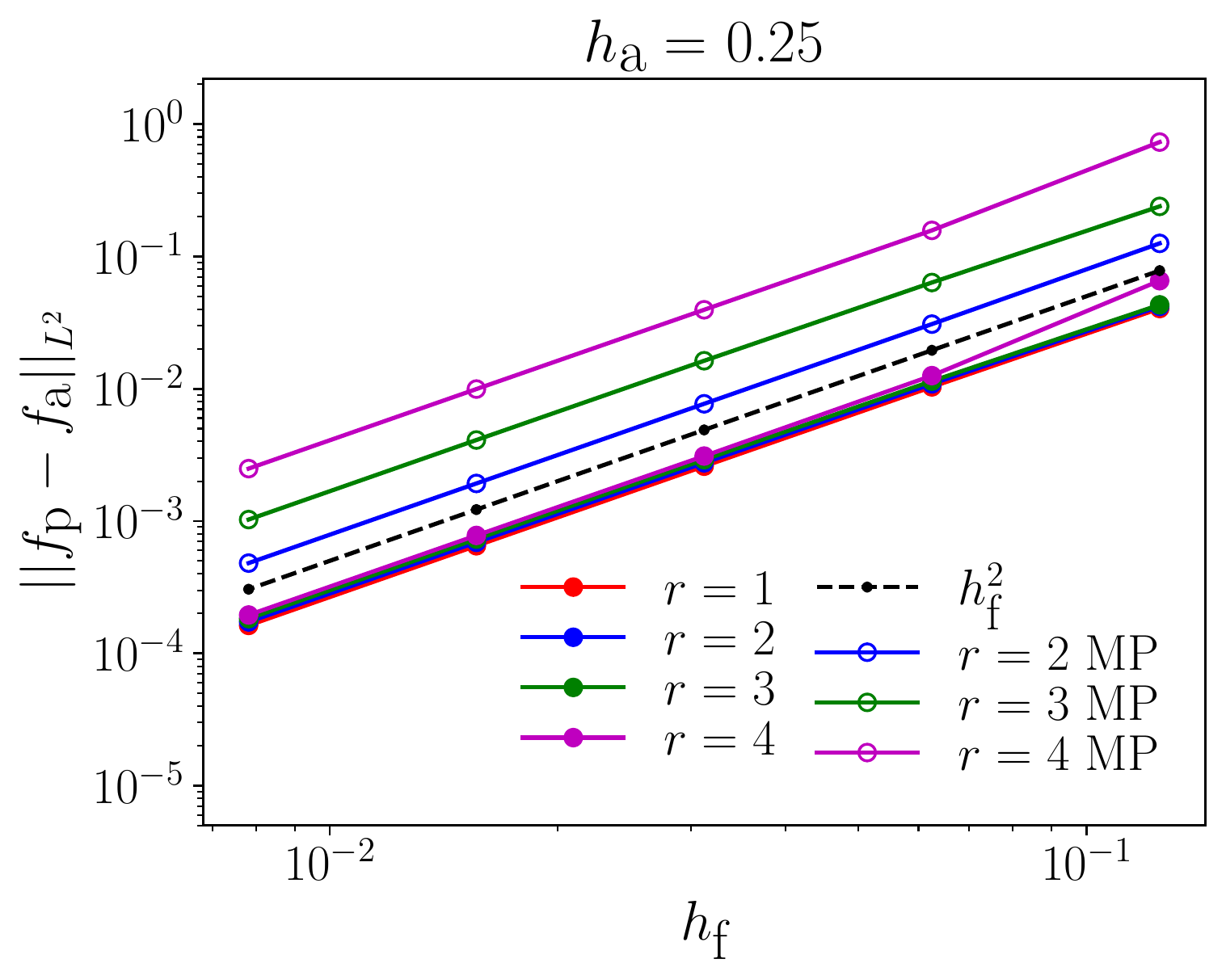}
	\end{subfigure}
	\begin{subfigure}{0.32\textwidth} 
		\includegraphics[width=\textwidth]{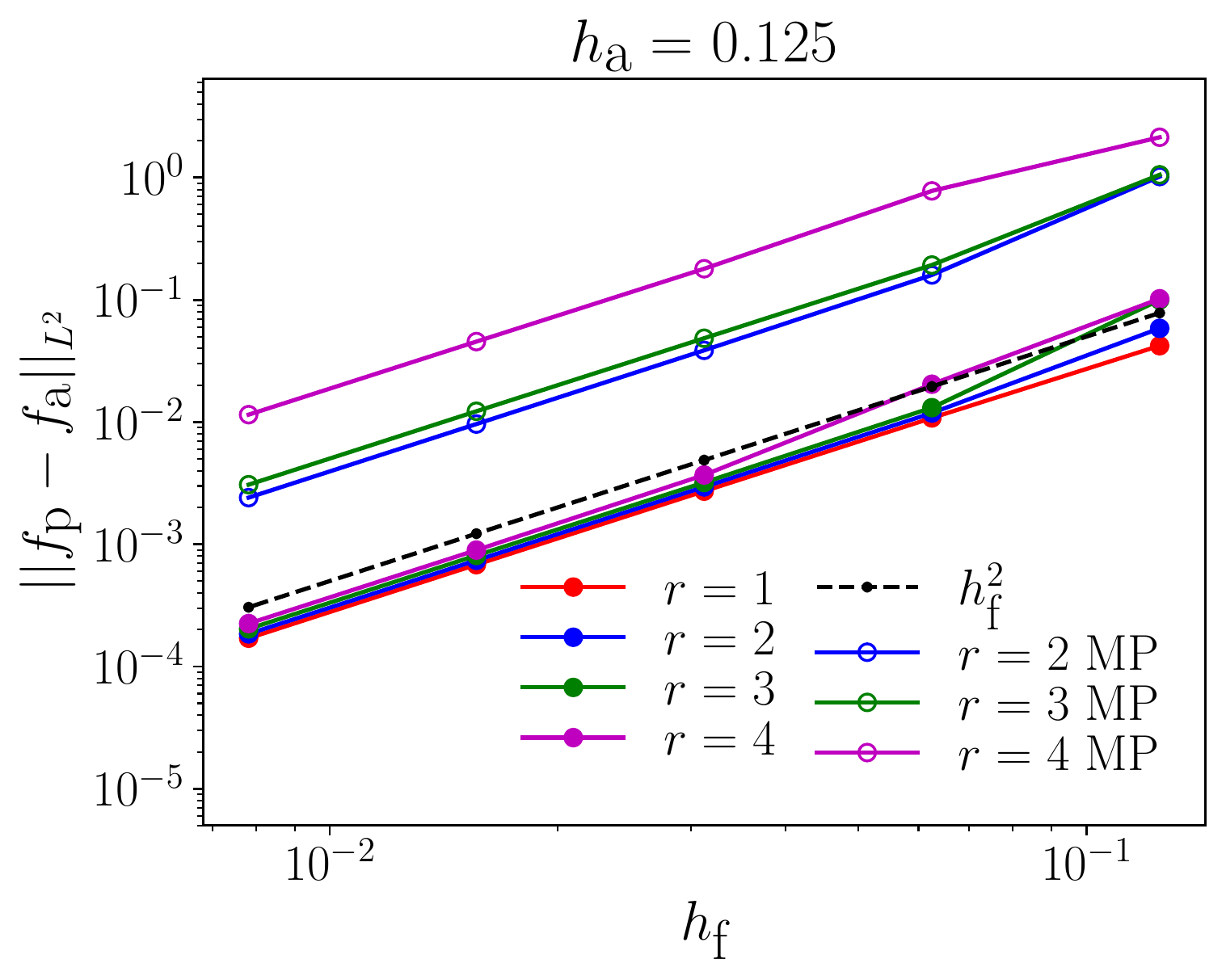}
	\end{subfigure}
	\begin{subfigure}{0.32\textwidth} 
		\includegraphics[width=\textwidth]{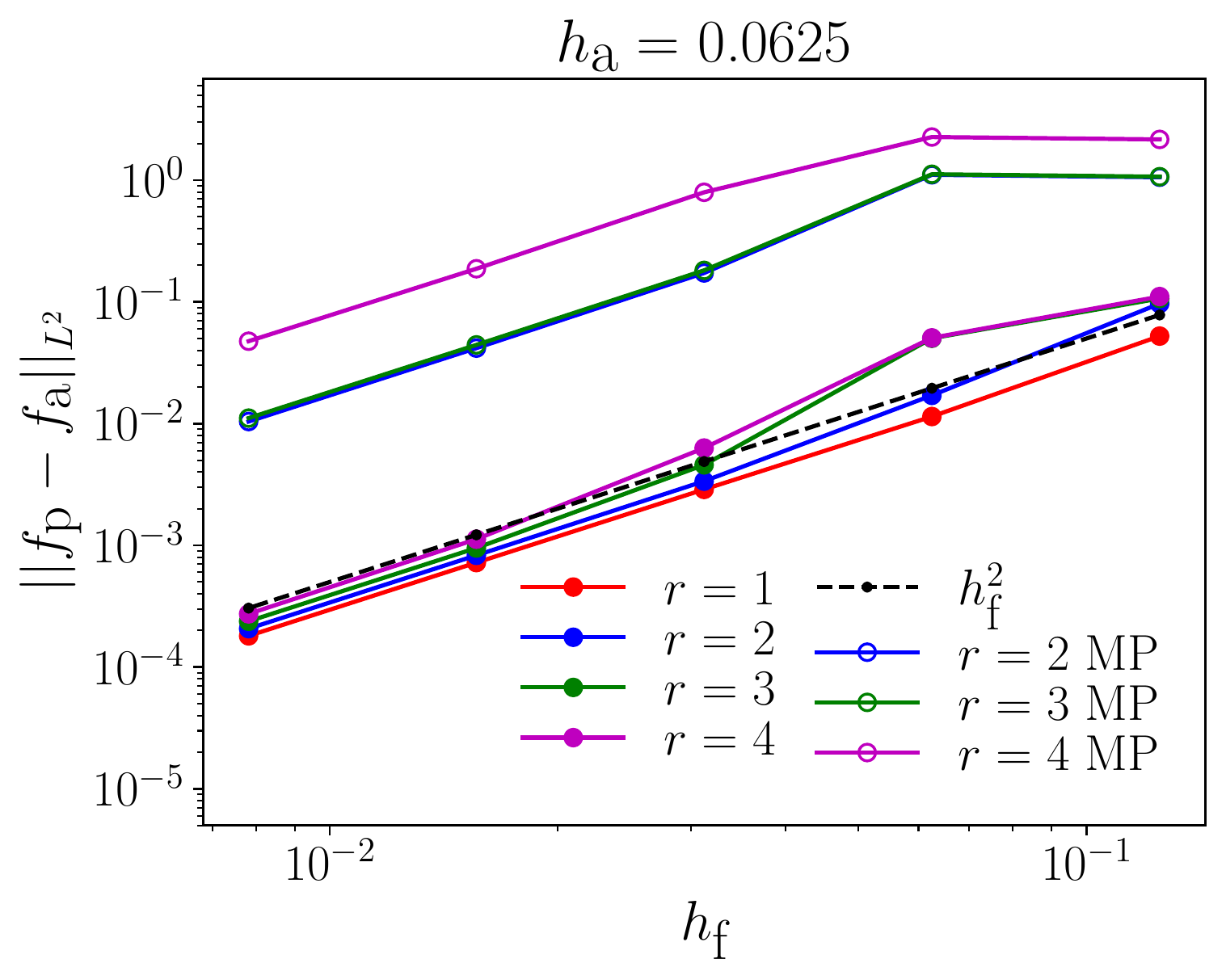}
	\end{subfigure}
	\caption{Computed errors $\norm{f_\textrm{p} - f_\textrm{a}}{L^2(\Omega)}$ versus $h_{\textrm{f}}$, for different polynomial degrees $r=1,2,3,4$ and different choices of $h_{\textrm{a}}=0.25,0.125,0.0625$  comparing our proposed projection method and mid point projection method (MP). \label{fig:ApproxFP-FA-PM}}
\end{figure}
We now consider the mid-point projection defined in Eq.~\eqref{eq:projKaltenbacher}. As already discussed, the difference between the projection $f_\textrm{a}$ defined in Eq.~\eqref{eq:ProjectionProblemDiscrete} and Eq.~\eqref{eq:projKaltenbacher} is the employed quadrature method.
We first observe that if the underlying map $\bm{\theta}_{ K_{\textrm{a}}}$ is linear for all $ K_{\textrm{a}} \in \TauAHrom$, for $r=1$ the two methods coincide. However, as we increase the polynomial degree or if we employ a trilinear map, the quadrature error increases showing that the projection computed with an exact quadrature formula is clearly more accurate, see Fig.~\ref{fig:ApproxFP-FA-PM}. 
In fact, if we want to employ high order polynomials on the acoustic grid, as already discussed before, it is crucial to minimize the projection error to exploit the accuracy provided by the spectral solver. 
From the numerical tests, it can be seen that the use of the midpoint projection method still provides accurate results since the convergence rate is $h^2_\textrm{f}$. However, the number of fluid elements required to have an accurate projection that does not interfere with the discretization error of the acoustic solver increases considerably. The latter makes the computational cost for the computation of the fluid solution, which is the real bottleneck of the workflow, very high.

\section{Aeroacoustic Applications} \label{sec:AeroResults}

In this section we apply the developed aeroacoustic hybrid strategy to relevant aeroacoustic benchmark problems. First, we test our strategy on a benchmark having an analytical solution, namely, the corotating vortex pair. This problem has been largely employed as a benchmark for aeroacoustic problems, see for instance \cite{mitchell1995}, \cite{lee1995} or \cite{Escobar2008}. Next, we consider the noise induced by the two-dimensional laminar flow around a squared cylinder.

\subsection{Corotating vortex pair}

\begin{figure}
    \begin{center}
\includegraphics[width=0.3\textwidth]{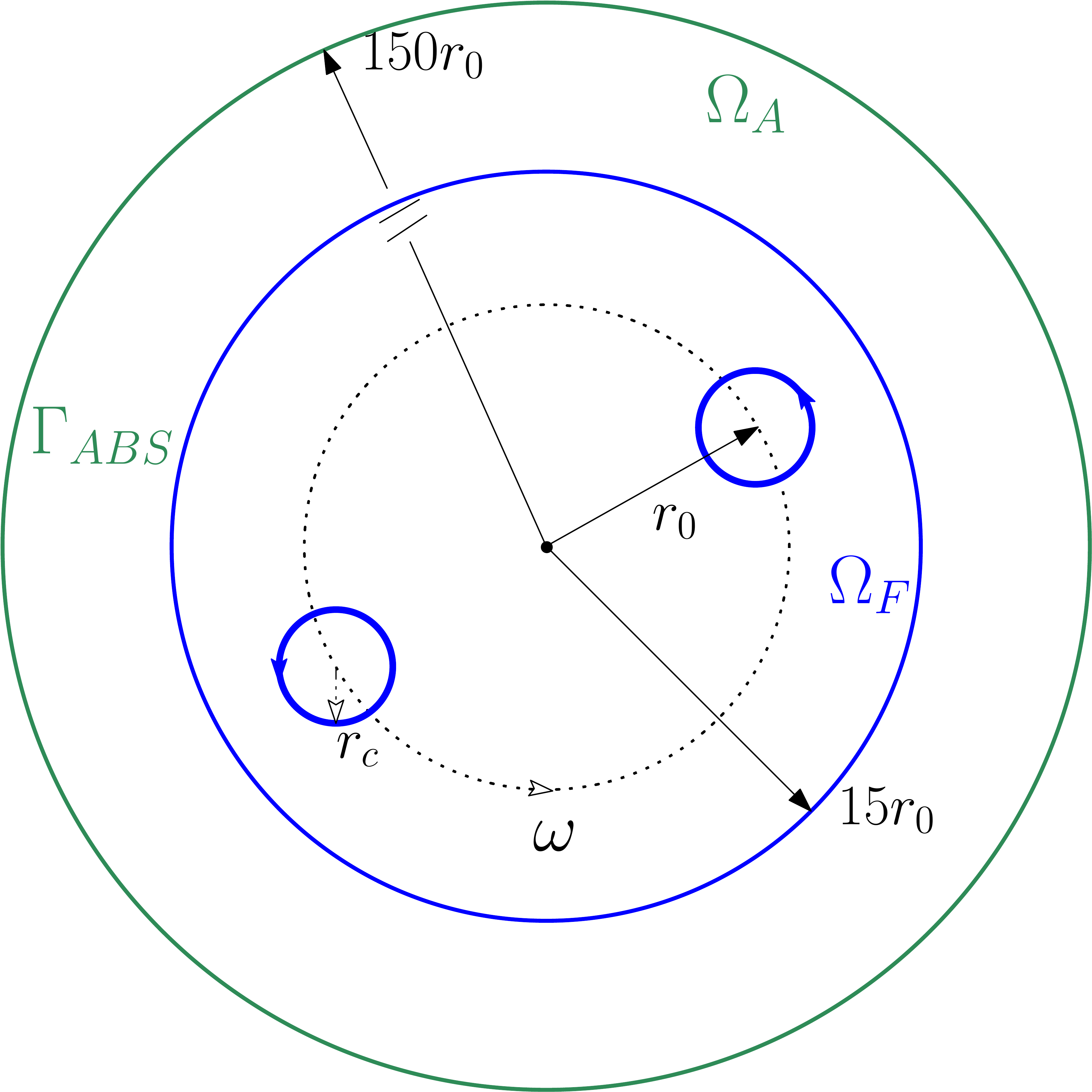}
    \end{center}
    \caption{Sketch of the domain for the corotating vortex pair problem. \label{fig:fluidVP}}
\end{figure}

We apply our hybrid aeroacoustic computational strategy to the corotating vortex pair problem. 
For this test case the fluid solution can be computed analytically based on potential flow theory. Furthermore, an analytical expression for the pressure fluctuations is obtained at the far field, for a detailed derivation of the analytical solution see \cite{Muller1967} or \cite{mitchell1995}. 
We assume that the flow field induced by the corotating vortex pair is inviscid and incompressible. This assumption allows us to employ a complex potential function $\Phi (z,t): \mathbb{C}\times(0,T] \rightarrow \mathbb{C}$ to describe the flow field, namely:
\begin{equation}
	\Phi (z,t) = \frac{\Gamma}{2\pi i} \ln(z - b(t)) + \frac{\Gamma}{2\pi i} \ln(z + b(t)), \label{eq:complexPotential}
\end{equation}
where $\Gamma $ is the circulation, $i$ is the imaginary unit and $b = r_0 \exp(i\omega t) \in \mathbb{C}$ are the rotating centres of the vortexes, where $\omega$ is the rotational speed defined as $\omega = \Gamma/(4\pi r_0^2)$ and $r_0$ is the distance with respect to the origin axes, see Figure~\ref{fig:fluidVP}. We introduce the rotating Mach number $\displaystyle M_r = \Gamma/(4\pi r_0 c_0)$, where $c_0$ is the speed of the wave. The period of the rotating monopoles is $\displaystyle T_{\textrm{f}} = 8\pi^2r_0^2/\Gamma$, while the emitted period of the acoustic wave is $T_{\textrm{a}} = T_{\textrm{f}} / 2$.
From the complex potential in \eqref{eq:complexPotential} we compute the two-dimensional fluid flow velocity $\mathbf{u} = [u;v]$ as
\begin{equation}
	u - iv = \frac{\partial}{\partial z}\Phi(z,t), \label{eq:complexVelocity}
\end{equation}
and then we compute the Lighthill's stress tensor. We report the far field solution for the pressure fluctuations $p' = p - \overline{p}$, see for instance  \cite{SchlottkeThesis}:
\begin{equation}
	p'(z,t) = -\frac{\rho_0 c_0^2}{64\pi^3}\left(\frac{\Gamma}{r_0c_0}\right)^4 [J_2(k \text{r})\sin(2(\theta-\omega t) + Y_2(k \text{r}) \cos(2\theta-\omega t)], \label{eq:analyticalSolution}
\end{equation}
where $J_2$ and $Y_2$ are respectively the first and second type Bessel functions, $k= 2\omega/c_0$ and $z = \text{r}\exp{(i\theta)}$.
As already showed in \cite{lee1995}, a desingularization model is required in order to avoid numerical issues in representing the source vortexes. Here, we employ the Scully model \cite{Scully1975} getting
\begin{equation}\label{eq:desingularization}
	u_{\theta}(r_v) = \frac{\Gamma r_v}{2\pi (r_c^2 + r_v^2)},
\end{equation} 
where $r_c$ is the desingularized core radius, $u_{\theta}(r_v)$ is the tangential velocity and $r_v$ is the distance with respect to the vortex core center.

\subsubsection*{Fluid Setup}
 We consider the corotating vortex pair problem with the parameters summarized in Table~\ref{tab:paramVP}. The fluid domain is a circle $\Omega_{\textrm{f}}$ with radius $15 r_0$. The flow solution is computed by employing the complex velocity in \eqref{eq:complexVelocity} and then the Lighthill's tensor  $\nabla \cdot \textbf{T} =  \rho_0\nabla\cdot(\mathbf{u}\otimes\mathbf{u})$ as a post-process of the flow velocity, see Section~\ref{sec:disc_NS}. The solutions are saved at each time instant with time step $\Delta t_{\textrm{f}} = 0.02$. Since the solution is periodic, we store the solutions up to $T_{\textrm{a}}$.
 
 
 \renewcommand{\arraystretch}{1.25}
 \begin{table} \centering
 \begin{tabular}{cccc cccc}
 $\Gamma~[\SI{}{\meter\square\per\second}]$ & $M_r$ & $T_{\textrm{a}}~[\SI{}{\second}]$ & $r_0~[\SI{}{\meter}]$ \\
 \hline
 0.98696 & 0.0785397 & 40  & 1  \\
 \hline\hline
 $\rho_0~[\SI{}{\kilo\per\meter\cubed}]$ & $c_0~[\SI{}{\meter\per\second}]$ & $\omega~[\SI{}{\per\second}]$ & $r_c~[\SI{}{\meter}]$\\
 \hline
 1  & 1  & 0.0785397  & 0.2  
\end{tabular}
\caption{Parameters employed for the rotating vortex pair test case.\label{tab:paramVP}}
\end{table}
\renewcommand{\arraystretch}{1} 

\subsubsection*{Acoustic Setup}
The acoustic domain $\Omega_{\textrm{a}}$ is a circle of radius $150r_0$. A circular domain has been chosen since the employed absorbing boundary conditions work better when the incident plane is parallel to the boundary, see \cite{Engquist1977}. On the external boundary $\Gamma_{\textrm{abs}}$ absorbing conditions are imposed. We employ a polynomial  degree $r=2$ for the SE discretization. We consider a structured grid meshing strategy with a total of 57500 elements. Each vortex has around at least 15 elements. For the time discretization, an implicit Newmark method is used, see for instance \cite{kaltenbacher2007book} or \cite{Mazzieri2020}, with $\beta = 0.5, \alpha = 0.25$  and  $\Delta t_{\textrm{a}} = 0.02$. 
In order to avoid spurious oscillation due to the non-consistent initial conditions, see \cite{EscobarThesis} or  \cite{Liow2006}, the following time ramp is multiplied by the source term  $ \displaystyle f(t) = \frac{1}{2}\left(1 - \cos\left(\pi\frac{t}{T_{end}}\right)\right)$,
 where $T_{end} = T_{\textrm{a}}$, up to $T_{end}$. 

\subsubsection*{Numerical results}
\begin{figure} \centering
 \includegraphics[width =0.45\textwidth]{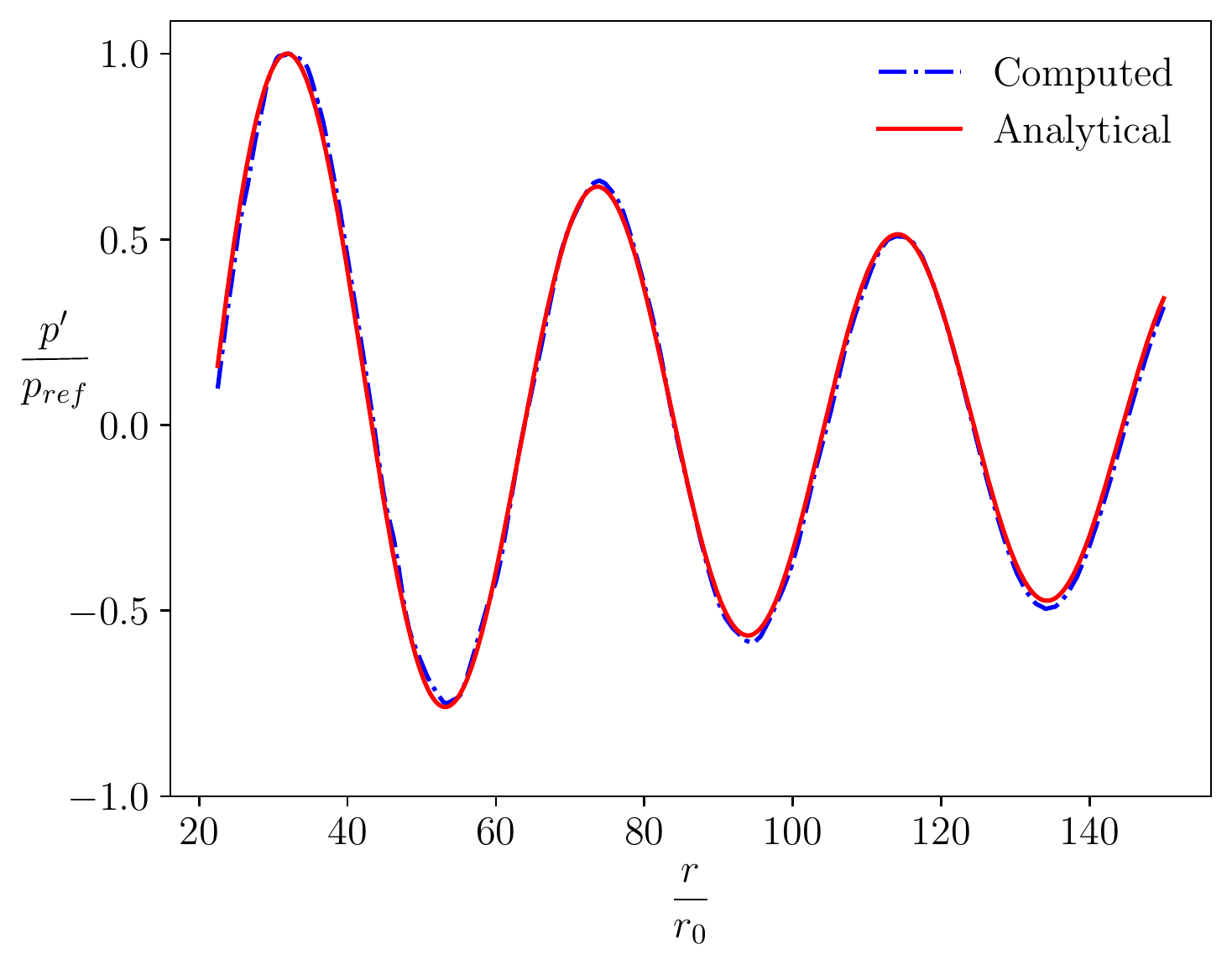} 
 \caption{ Comparison between the analytical far field solution and the computed numerical solution obtained with the hybrid approach and by employing the vortex core model. The results have been normalized with respect to $p_{ref} = \max(p-\overline{p})$ to take into account the energy disparity introduced by the vortex model. \label{fig:ProbeSolutionVP}}
\end{figure}

\begin{figure}
\begin{subfigure}{0.47\textwidth}  \centering
 \includegraphics[width =0.9\textwidth]{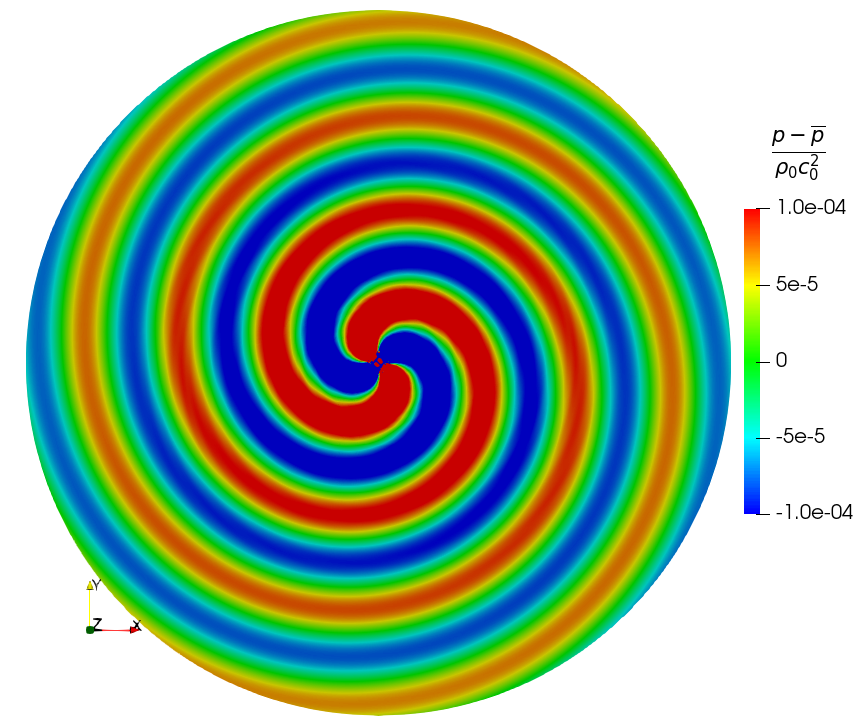}
 	\caption{$t = \SI{365}{\second}$.}
 \end{subfigure} \hfill
\begin{subfigure}{0.47\textwidth} \centering
 \includegraphics[width =0.9\textwidth]{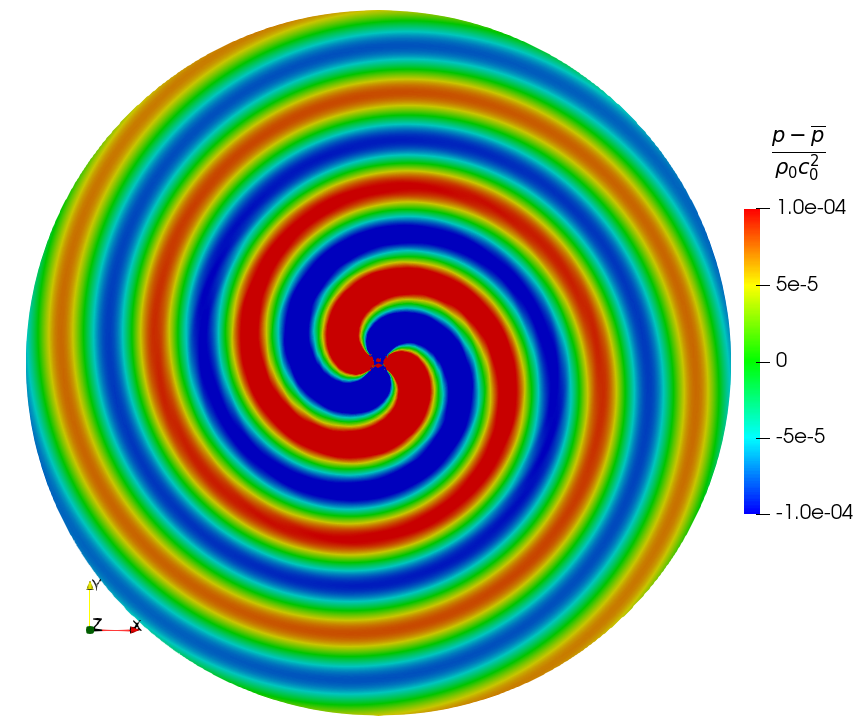} \\
 	\caption{$t = \SI{370}{\second}$.}
 \end{subfigure}
\begin{subfigure}{0.47\textwidth}  \centering
 \includegraphics[width =0.9\textwidth]{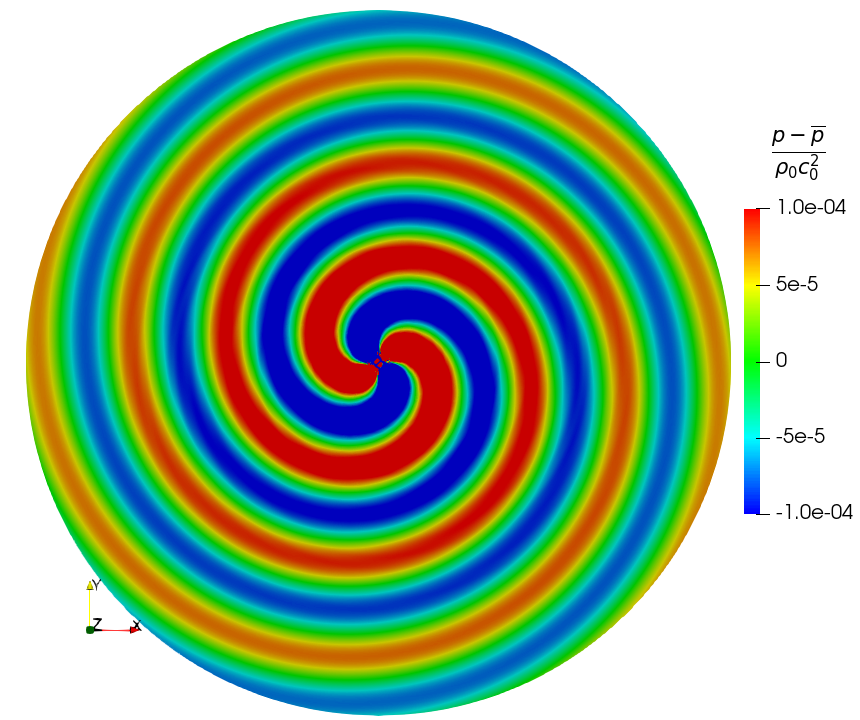} 
 	\caption{$t = \SI{375}{\second}$.}
 \end{subfigure} \hfill
\begin{subfigure}{0.47\textwidth}  \centering
 \includegraphics[width =0.9\textwidth]{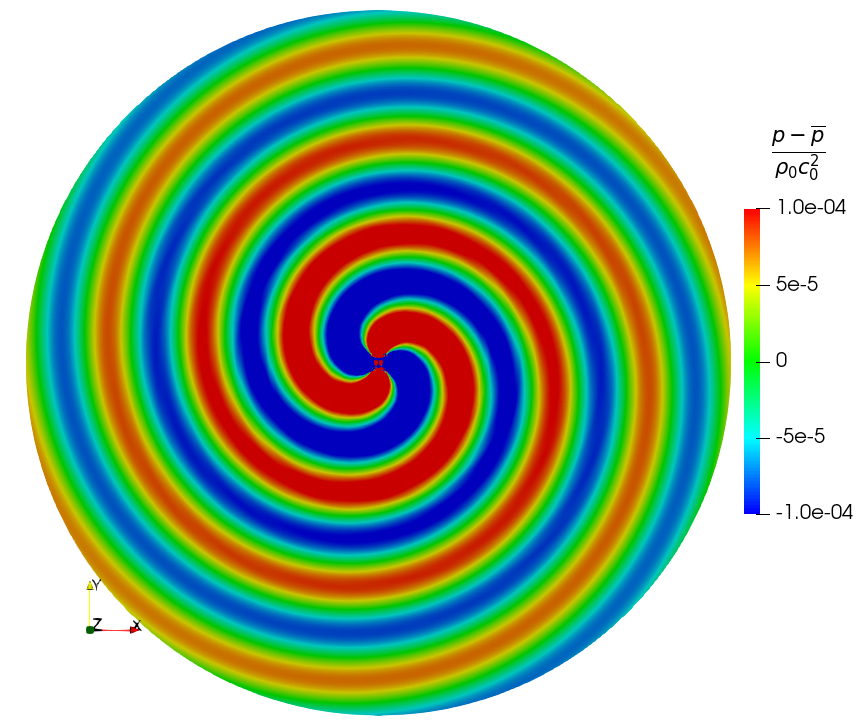} \\
 	\caption{$t = \SI{380}{\second}$.}
 \end{subfigure}
 \caption{Snapshots of the computed numerical solution for the corotating vortex pair for $t = 365,370,375,380 \ \SI{}{\second}$. \label{fig:SnapshotSolutionVP} } 
\end{figure}

The acoustic field generated by a pair of corotating vortices
is a rotating acoustic quadrupole as can be seen from Figure~\ref{fig:SnapshotSolutionVP}. The numerical solution obtained through the proposed algorithm matches the analytical solution as it is shown in Figure~\ref{fig:ProbeSolutionVP}, where we sampled the pressure fluctuations $p'$ along the line $y=0$ with $\displaystyle x>0$ at $T = \SI{380}{s}$. The results obtained with the analytical solution have been normalized  to a reference pressure $p_{ref} = \max(p')$ in order to take into account the desingularization effect in \eqref{eq:desingularization}, see for instance \cite{Schlottke2016}.

\subsection{Flow around a squared cylinder at low Reynolds number}

\label{sec:LaminarFlowCylinder}
Finally, we consider the case of a laminar flow around a square cylinder, see for instance the Direct Numerical Simulation (DNS) performed by \cite{Inoue2006} or the solution obtained with a Curle analogy in \cite{Ali2010}.
When a rigid squared cylinder is placed in a uniform flow, it exhibits strong vortex shedding, resulting in fluctuating forces due to the alternating pressure highs and drops at the wake. These forces and the turbulence in the wake generate noise. For laminar flows, the main frequency radiated from the body is associated with the Strouhal number and the intensity of the observed noise is proportional to the fluctuation of the forces. The flow solution has been computed by employing the Pressure Implicit Splitting Operator (PISO) method implemented in OpenFOAM \cite{OpenFOAM}.

\subsubsection*{Fluid Setup}

\begin{figure}
    \centering
    \includegraphics[width=0.5\textwidth]{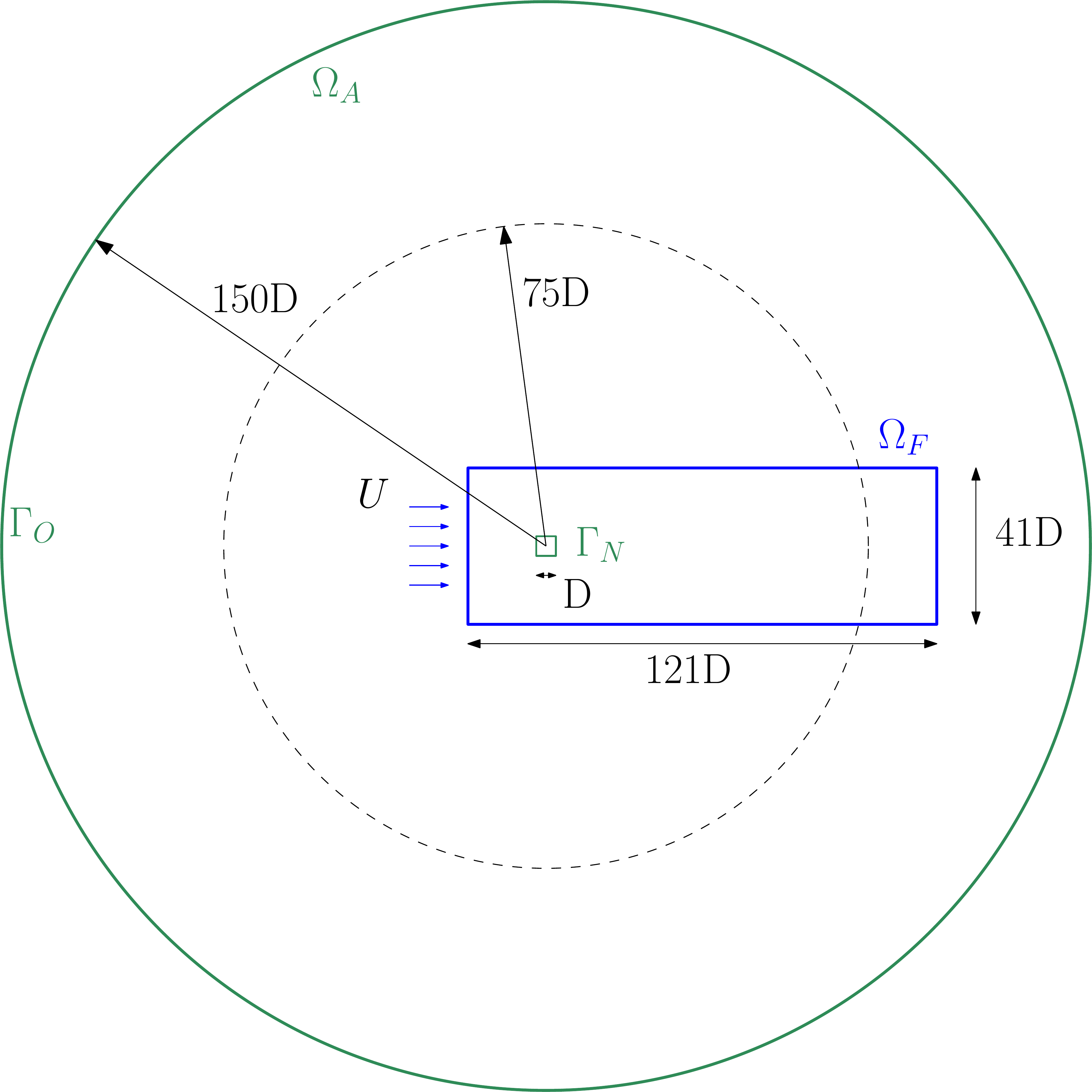}
    \caption{Computational domain of the fluid problem and the acoustic problem. The square cylinder has a diameter $D=3.28\times10^{-5}\SI{}{\meter}$ and the fluid domain is a rectangle of size $121D\times41D$. The acoustic domain is a circle of radius $150D$, centered at the centre of the square. The dotted line represent the sampled probes in the acoustic domain employed to compute the directivity, see Fig.~\ref{fig:DirectivityCyl}.}
    \label{fig:aeroCilindroDominio}
\end{figure}

\begin{figure}
    \centering
    \includegraphics[width=0.5\textwidth]{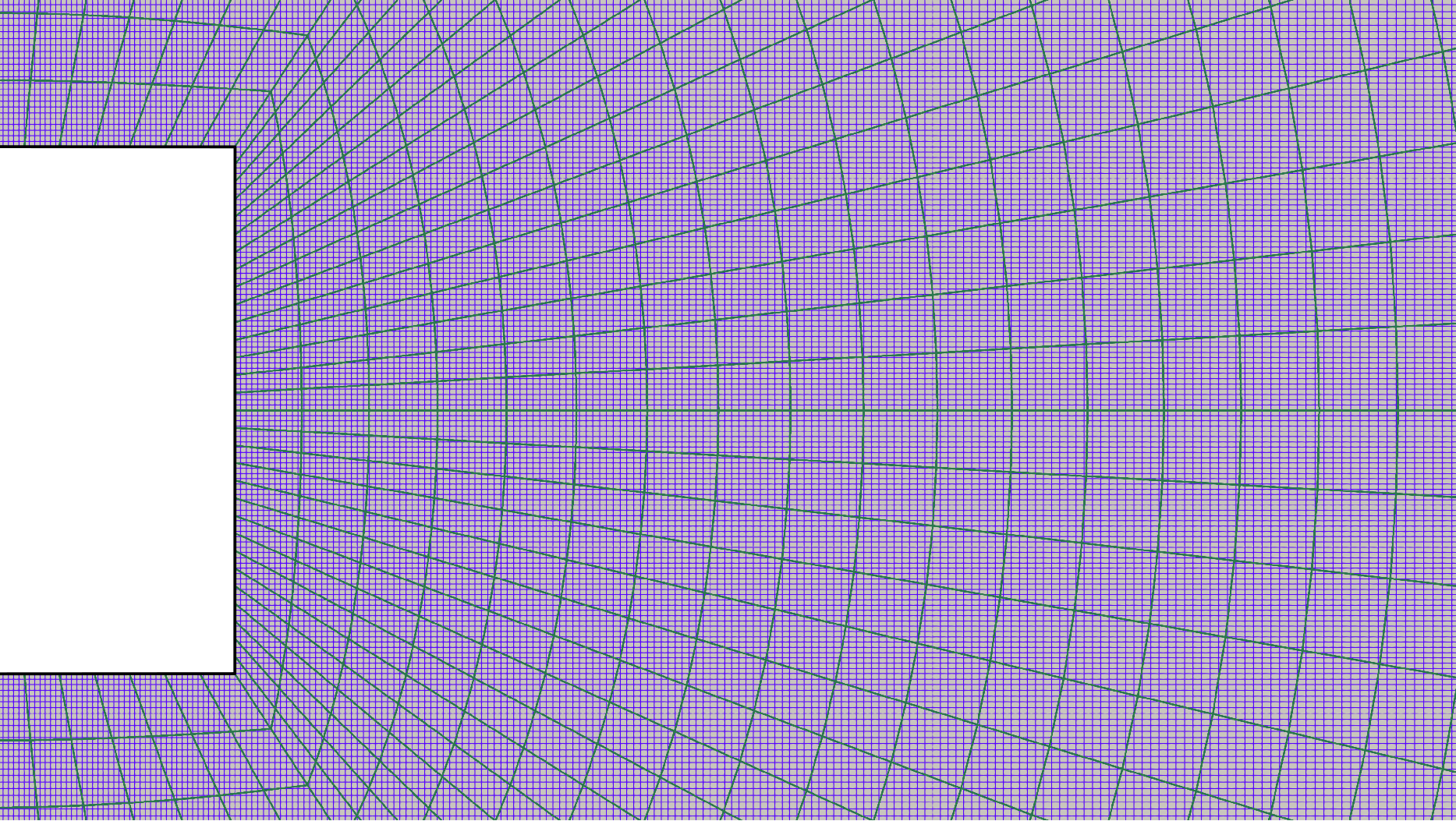}
    \caption{Detail of the acoustic (green) and fluid (blue) computational grids around the square cylinder for the aeroacoustic test case.}
    \label{fig:Zoom}
\end{figure}

A laminar two-dimensional incompressible simulation of a fluid flow around a square cylinder is performed. Let $D=3.28\times10^{-5}~\SI{}{m}$ be the length of the square cylinder, $U = \SI{68.7}{\meter\per\second}$ be the inlet velocity  and  $\nu = 1.5\times10^{-5}\SI{}{\meter\squared\per\second}$ the kinematic viscosity. The Reynolds number is $\textrm{Re}=150$ and the Mach number is $\textrm{Ma}=0.2$.
The fluid computational domain $\Omega_{\textrm{f}}$ is $(-20.5D,100.5D)\times(-20.5D,20.5D)$, see Figure~\ref{fig:aeroCilindroDominio}. A fixed velocity $U$ is prescribed at the inlet. On the upper and lower wall symmetry conditions are employed. No slip conditions are applied on the cylinder walls. Zero gradient pressure conditions are applied at the wall of the cylinder and at the inlet. On the outlet, the pressure is set to zero, while a zero gradient condition is imposed for the velocity. A block structured h-grid around the square cylinder is used, employing $970000$ elements. In Fig. \ref{fig:Zoom} a zoom of the fluid grid is shown. The computational time step is $\Delta t = 10^{-9}\SI{}{s}.$

\subsubsection*{Acoustic Setup}

To minimize spurious reflections, the acoustic domain $\Omega_{\textrm{a}}$ is a circle of radius $150D$, with an internal square hole of side $D$, see Figure~\ref{fig:aeroCilindroDominio}. On the external boundary $\Gamma_{\textrm{abs}}$, absorbing conditions are employed. On the solid wall $\Gamma_{\textrm{b}}$ Neumann boundary conditions are imposed. The coupling region is given by $\Omega_{\textrm{a}}\cap\Omega_{\textrm{f}} = \Omega_{\textrm{f}}$. A smoothing function is employed in order to let the sound source term decay to avoid the well known spurious noise generation due to the abrupt domain cut on the wake, see for instance \cite{Oberai2000} or \cite{MartinezLera2008}. We used the following spatial smoothing function:
\begin{equation*}
	g(x) = \begin{cases}
		1 \quad &x < r_i, \\
		\displaystyle \frac{1}{2} \left( 1 + \cos(\pi \frac{x - r_i}{r_o - r_i})\right)  \quad &x \ge r_i,
	\end{cases}
\end{equation*}
where $r_i$ is the initial filtering position, while $r_o$ is the end of the fluid domain. In this case, it is sufficient to apply the smooth function only downstream and along the $x$ direction, so that $r_i = 65D$ and $r_o = 105.5D$.
The fluid solution is sampled every 10 fluid time steps, meaning that $\Delta t_{\textrm{a}} = 10^{-8} \SI{}{s} = 10 \Delta t_{\textrm{f}}$. The expected main frequency is the Strouhal frequency. The acoustic discretization close to the square is h-type, with $h_{\textrm{a}} = D/10$. Then the grid is unstructured and a o-type grid is employed. The  polynomial degree chosen is $r=4$. The whole acoustic grid has around $1261500$ degrees of freedom. The main wavelength associated to the lift force is $\lambda\approx32.5D$ and around 40 nodes per wavelength where placed in the far field. The acoustic simulation was run for $0.00005\SI{}{s}$, starting from a fluid time of $t_{\textrm{f}} = \SI{0.0001}{s}$, hence with a fully developed flow field. A zoom of the acoustic grid is shown in Fig~\ref{fig:Zoom}. Note that the acoustic element size is larger than the size of the fluid elements.

\subsubsection*{Numerical results}
\begin{figure} \centering
    \includegraphics[width=0.5\textwidth]{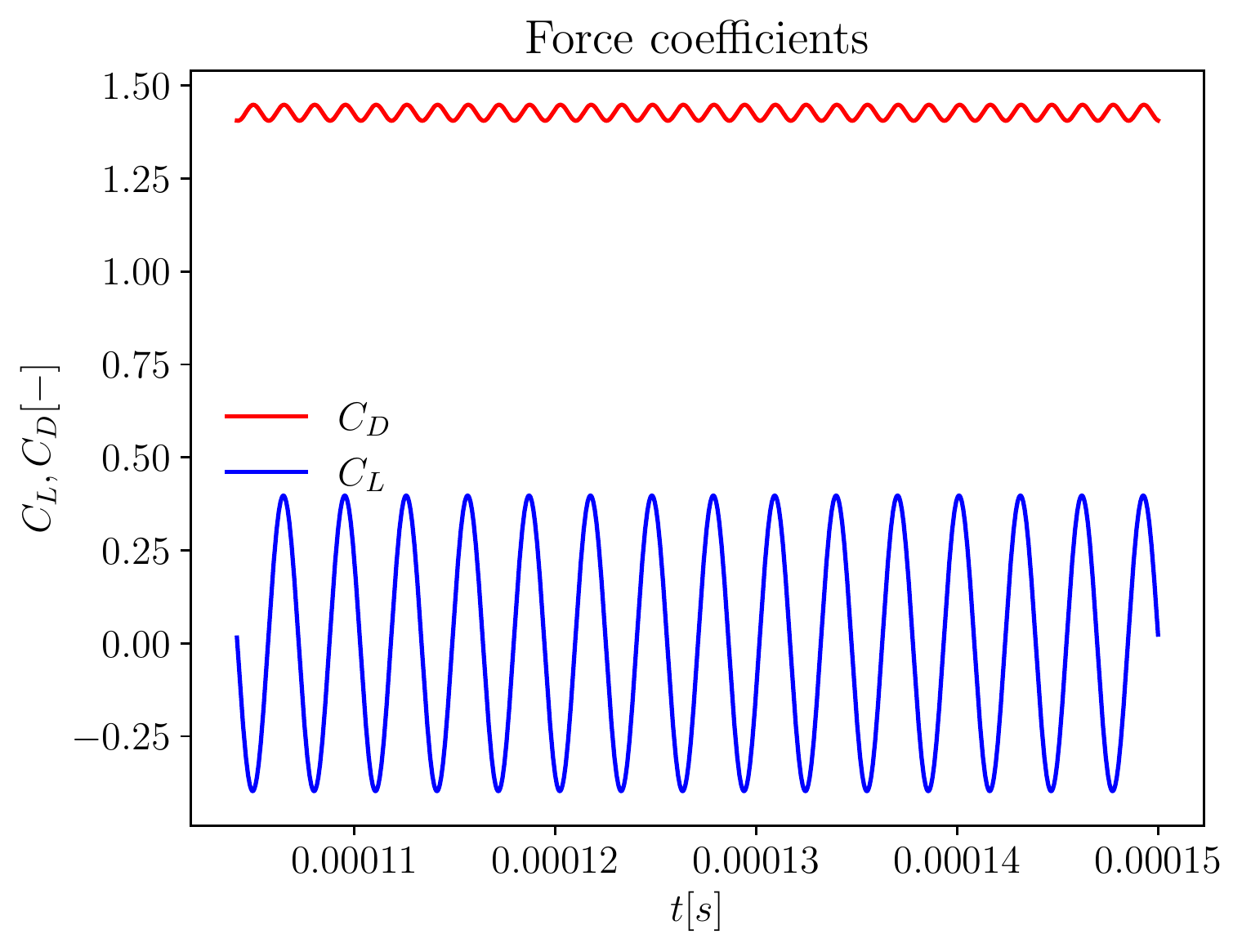}
    \captionof{figure}{Computed $C_L$ and $C_D$ coefficients. 
    \label{fig:forces}}
\end{figure}

\begin{figure}
    \begin{subfigure}{0.47\textwidth}
    \centering
		\includegraphics[width=0.9\textwidth]{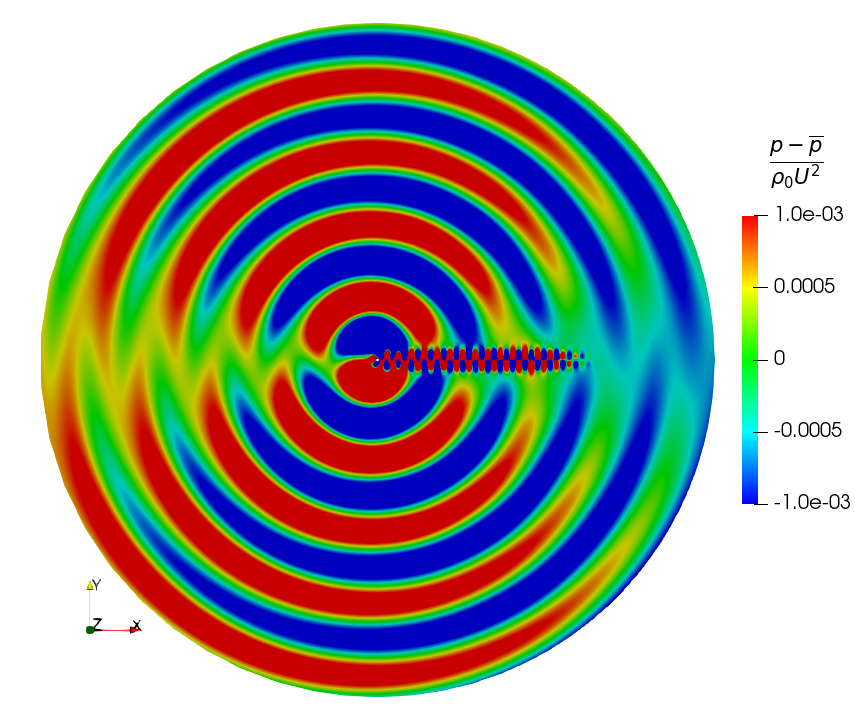}
  \caption{$t=\SI{0.000147}{\second}$.}
	\end{subfigure}
\hfill
    \begin{subfigure}{0.47\textwidth} \centering
	\includegraphics[width=0.9\textwidth]{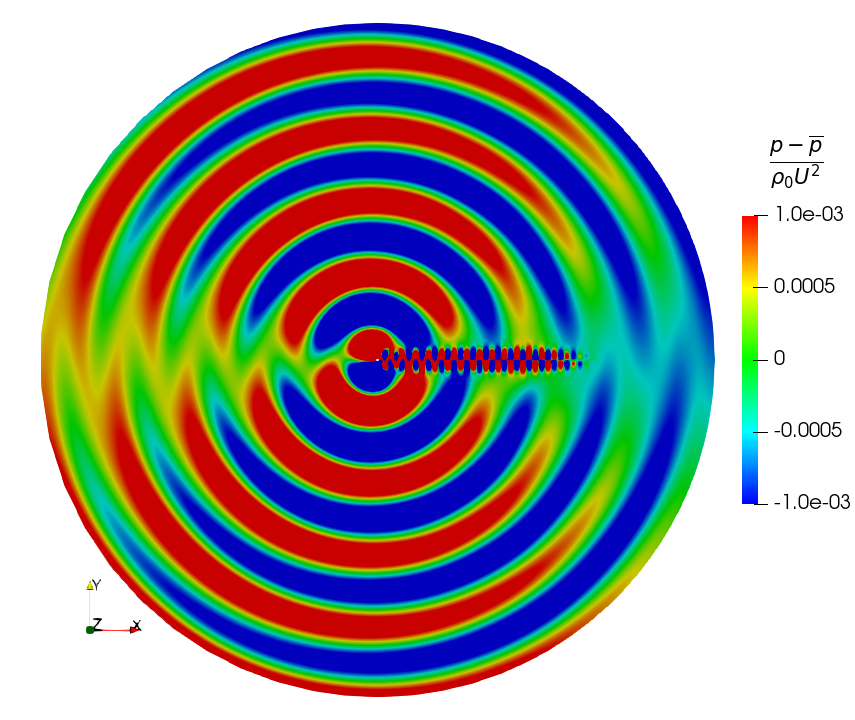}
   \caption{$t=\SI{0.000148}{\second}$.}

	\end{subfigure}
    \begin{subfigure}{0.47\textwidth}
    \centering
		\includegraphics[width=0.9\textwidth]{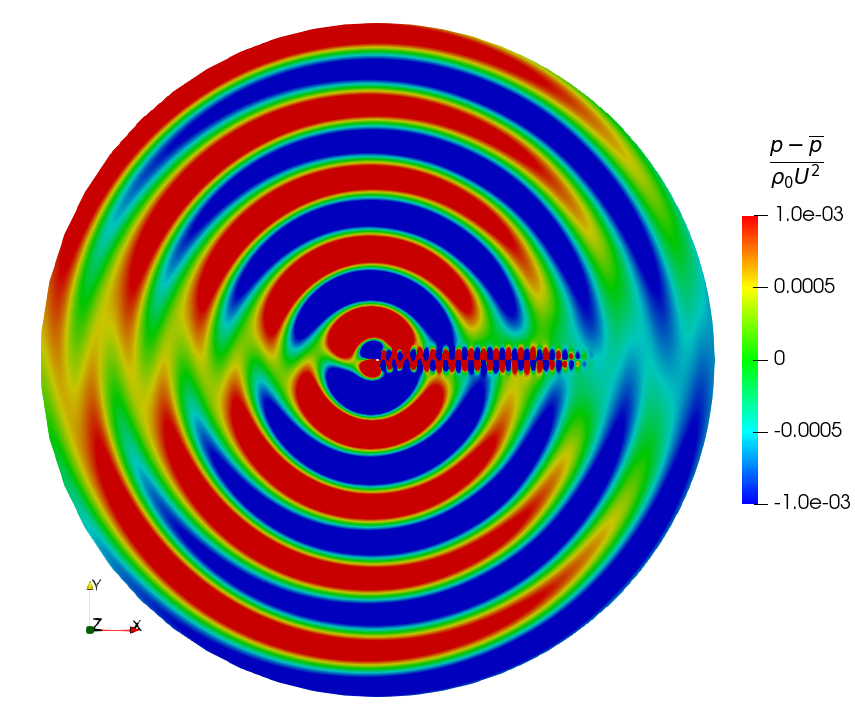}
    \caption{$t=\SI{0.000149}{\second}$.}

	\end{subfigure}
\hfill
    \begin{subfigure}{0.47\textwidth} \centering
	\includegraphics[width=0.9\textwidth]{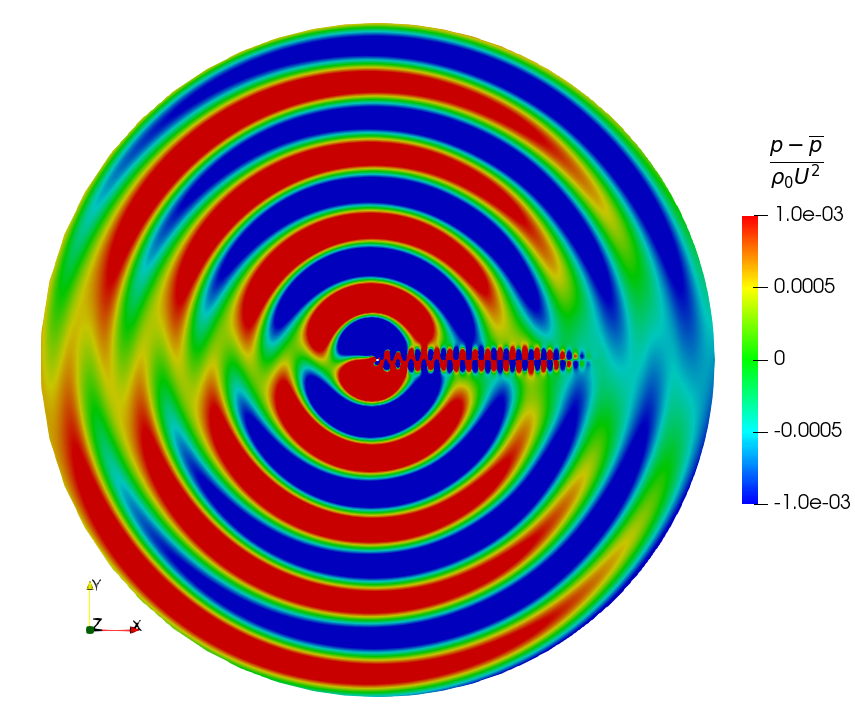}
   \caption{$t=\SI{0.00015}{\second}$.}
	\end{subfigure}

 		\caption{Snapshot of the computed acoustic pressure field at t = $0.000147,\ 0.000148,\ 0.000149,\  0.00015\ \SI{}{\second}$. \label{fig:SnapshotSolutionCyl}}
\end{figure}

 $\displaystyle C_D = \frac{F_D}{\frac{1}{2}\rho_0 U^2 A}$ and $\displaystyle C_L = \frac{F_L}{\frac{1}{2}\rho_0 U^2 A}$, where $F_D$ and $F_L$ are the drag and lift forces respectively, with $A = D\times H$, being $H$ the width of the domain and having chosen $H=D$, we plot $C_L$ and $C_D$ in Figure~\ref{fig:forces}. We introduce the Strouhal number $\displaystyle St = f\frac{D}{U}$, with $f$ being the frequency of the $C_L$. In Table~\ref{table:CFD} we compare our results with those available in the literature. The obtained Strouhal number matches the results obtained by \cite{Doolan2009} and they are aligned with the experiments \cite{Okajima1982,Sohankar1999} and the compressible DNS performed by \cite{Inoue2006}. 
The intensity of the noise emitted by the square cylinder depends mainly on the fluctuations of the forces. Hence, during the flow computation it is critical to match the root means squared ($rms$) values.
By defining $\overline{C}_D$ as the average of $C_D$, respectively $C_L$, we compute the $rms$ values as $\displaystyle C_{L,rms} = \sqrt{\overline{(C_L - \overline{C}_L)^2}}$ and we also identify $C_{L,peak} = \max(|C_L|)$. Again, from Table~\ref{table:CFD}  we see that out results are in agreement with the ones available in literature.
Finally, we compute the acoustic field, namely the noise induced by the flow around the square cylinder.
From Figure~\ref{fig:SnapshotSolutionCyl}, we see the characteristic dipole pattern, that is mainly due to the lift force acting on the cylinder. The obtained results are comparable with the compressible simulations, see for instance \cite{DAlessandro2019} and  \cite{Inoue2006}. To further validate the obtained acoustic results, we computed the directivity from $p'_{rms} = \sqrt{\overline{(p - \overline{p})^2}}$. The obtained directivity pattern is in good agreement with the references, see \cite{Inoue2006} and \cite{Ali2010}.

  \begin{table} \centering 
		\begin{tabular}{lcccc}	
			& $St$ & $\overline{C}_D$ & $C_{L,rms}$ & $C_{L,peak}$\\
			\hline
			Experiments \cite{Okajima1982,Sohankar1999} &  0.148-0.155 & 1.4 & - & - \\
			Doolan  \cite{Doolan2009}& 0.156 & 1.44 & 0.296 & - \\
			Ali \cite{Ali2010} & 0.16 & 1.47 & 0.285 & - \\
			Inoue \cite{Inoue2006}& 0.151 & 1.4 &   - & 0.4 \\
			Current study & 0.156 & 1.43& 0.281 & 0.3976
		\end{tabular}
		\captionof{table}{Comparison of the flow results with analogous results available in literature. \label{table:CFD}}
\end{table}

\begin{figure}
	\begin{center}
		\includegraphics[width=0.5\textwidth]{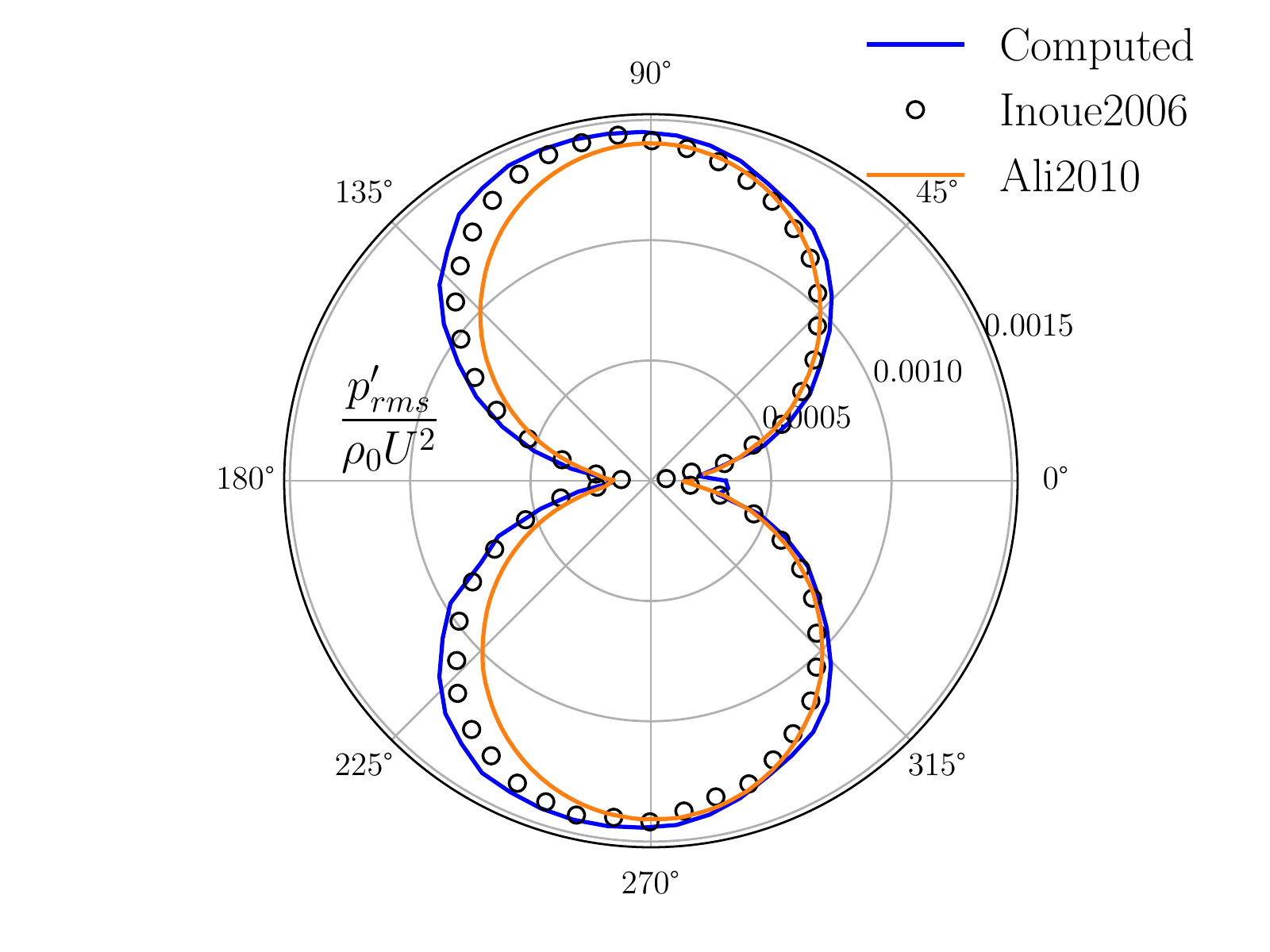}
		\caption{Directivity pattern. The adimensionalized $p'_{rms}$ has been sampled on a circumference of radius 75D, see Fig.~\ref{fig:domainAeroProblem}. Comparison with the DNS in \cite{Inoue2006} and the Curle computations of \cite{Ali2010}. \label{fig:DirectivityCyl}}
	\end{center}
\end{figure}



\section{Conclusion}
We proposed a hybrid computational strategy that couples a finite volume flow solver with a high order acoustic solver for aeroacoustic simulations. First, the fluid flow solution is computed employing the open-source finite volume library OpenFOAM. Then, a post-processing of the flow solution computes the sound source term on the fluid grid, by means of the Lighthill's acoustic analogy. Next, a projection method is used to map the flow source term from the fluid to the acoustic grid. Finally, an inhomogeneous wave equation is solved by employing a high-order spectral element method. The employed projection method exploited a robust intersection algorithm that is able to perform the intersection between the two computational grids. Furthermore, we employed a quadrature free method to integrate polynomial functions over the generic polyhedral elements stemming after the intersections computation. 
We explored the computational aspects of the proposed intersection algorithm both from a theoretical and numerical point of view.  
Finally, we applied the developed computational strategy to different aeroacoustic problems, showing the effectiveness of the proposed method.



\section{Acknowledgements}
The authors thank prof. M. Verani, prof. R. Corradi and Dr. P. Schito for the insightful discussions on the topic.  \Ar{We also thank the anonymous reviewers for carefully reading the manuscript and for the insightful comments and suggestions. }
The simulations have been partly run at Cineca thanks to the computational resources made available through the HO-AERO HP10C9XBN9 ISCRA-C project.  A.A., P.F.A., I.M. e N.P. are members of the INdAM Research Group GNCS. P.F.A. has been partially funded by the research projects PRIN n. 201744KLJL, funded by MIUR, and P.F.A. and N.P. have been partially supported by PRIN n. 20204LN5N5 research grant funded by MIUR.
P.F.A., I.M. and N.P. have been partially supported by ICSC—Centro Nazionale di Ricerca in High Performance Computing, Big Data, and Quantum Computing funded by European Union—NextGenerationEU.

\appendix
\section{Appendix: finite volume approximation} \label{appendix:discFV}
	\begin{figure}[h!]
		\begin{center}
	\includegraphics[width=0.35\textwidth]{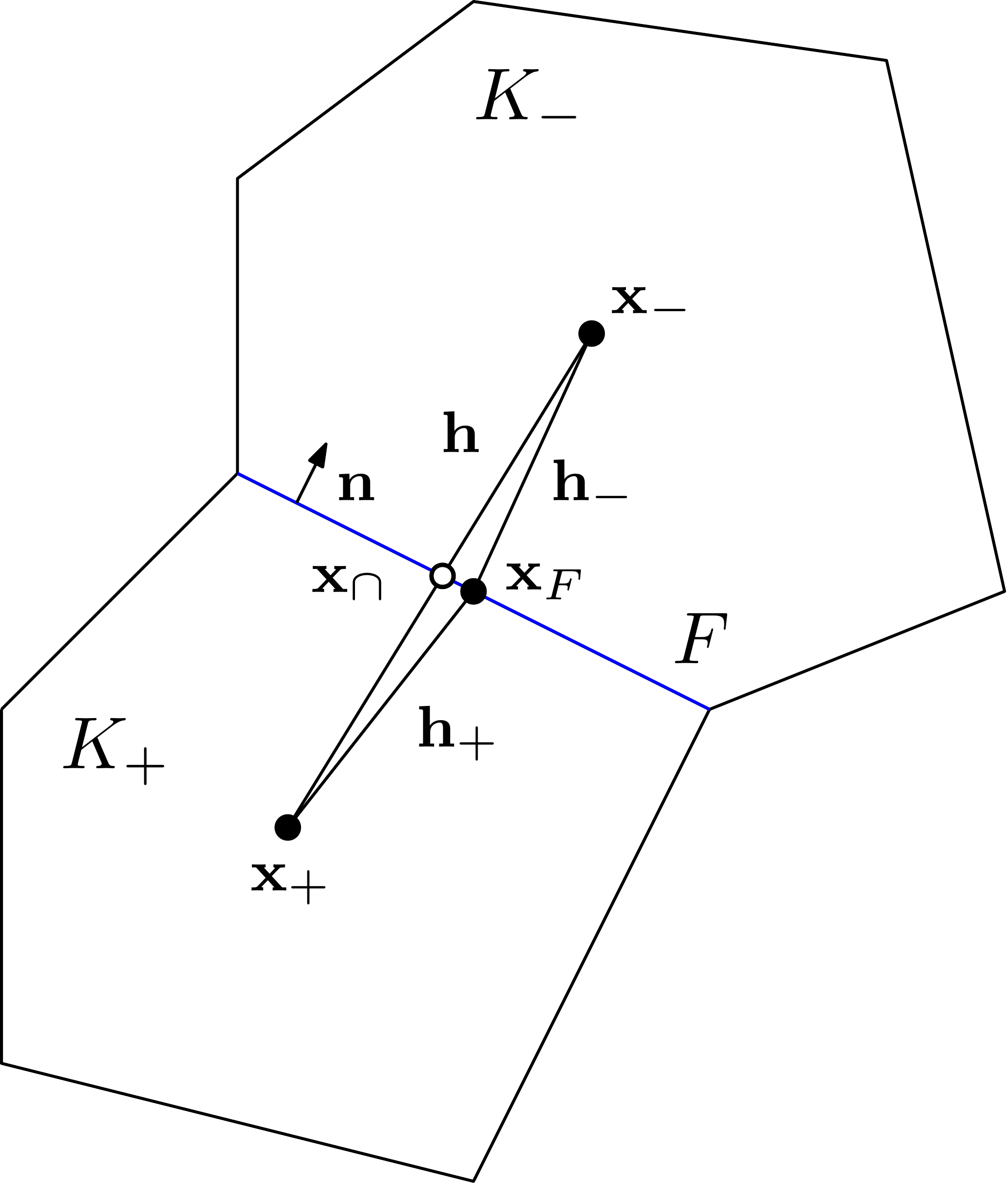}
	\end{center}
	\caption{\textcolor{blue}{Sketch of the geometrical notation for the finite volume discretization.  } \label{fig:FiniteVolume}}
\end{figure}
\Ar{
	Let us introduce some further geometrical notation, necessary to handle the discretization of the finite volume method on non-structured grids, see Fig.~\ref{fig:FiniteVolume}. \\ 
	Given two neighbouring cells $ K_+$ and $ K_-$, we denote with $ F$ their common face. Also, let $\mathbf{x}_+$ and $\mathbf{x}_-$ be the barycentres of the polyhedrons $ K_+, K_-$, with $\mathbf{x}_{ F}$ the barycentre of the face, and with $\mathbf{n}$ the unit normal to the face $ F$, outward with respect to the element $ K_+$. Now, let $\mathbf{h}$ be the vector connecting the cell barycentres $\mathbf{h} = \mathbf{x}_+ - \mathbf{x}_-$, and let $h = |\mathbf{h}|$.
	Now, let $\mathbf{x}_{\cap}$ be the intersection point between the face $F$ and the vector $\mathbf{h}$. 
 	We introduce the interpolation weights $w$:
	\begin{equation} \label{eq:interpolationWeights}
		w = \frac{h_+}{h_+ + h_-},
	\end{equation}
	where $h_+ = |\mathbf{x}_{\cap} - \mathbf{x}_+|$ and $h_- = |\mathbf{x}_{\cap} - \mathbf{x}_-|$.
	First, recall the following Gauss gradient reconstruction on the barycentre of the fluid element. Given an element $K$, its Gauss gradient at the cell centre is approximated by:
	\begin{equation} \label{eq:GaussGrad}
		\nabla \mathbf{u} \approx \frac{1}{|K|}\sum_{F \in \partial K} \mathbf{u}_F \mathbf{n}_F |F|,
	\end{equation}
	where $\mathbf{n}_F$ is the outward normal face to $F$, and the value of $\mathbf{u}_F$ is computed with a linear interpolation:
	\begin{equation}
		\mathbf{u}_F = w \mathbf{u}_+ + (1-w) \mathbf{u}_-,
	\end{equation}
	where $\mathbf{u}_{\pm}$ is the velocity field evaluated at the center of the cells $K_\pm$.
	The gradient $\nabla\mathbf{u}_F \mathbf{n}$ coming from the discretization of the diffusion term in eq.~\eqref{eq:DiffusionFV} is computed linear approximation including a non-orthogonal correction (see \cite{Jasak1996}):
	\begin{equation} \label{eq:DiffusionGrad}
		\nabla\mathbf{u}_F\mathbf{n} \approx \frac{\mathbf{u}_+ -  \mathbf{u}_-}{\mathbf{h}\cdot \mathbf{n}} + (w \nabla \mathbf{u}_+ + (1-w) \nabla \mathbf{u}_-) \left(\mathbf{n} - \frac{\mathbf{h}}{\mathbf{h}\cdot \mathbf{n}} \right),
	\end{equation}
	where $\nabla\mathbf{u}_\pm$ are the gradients computed with a Gauss formula at the cell centres of the elements $K_\pm$, see eq.~\eqref{eq:GaussGrad}. Note that in the case of structured orthogonal grids, the formula simply reduces to $\displaystyle \nabla\mathbf{u}_F\mathbf{n} \approx \frac{\mathbf{u}_+ - \mathbf{u}_-}{h}$. \\
	Concerning the convective term in  \eqref{eq:finite_vol_1}, the usual algorithms employed by OpenFOAM treat the convective term in an explicit way. So, we assume that the flux $\mathbf{u}\cdot\mathbf{n}$ is known, we denote it with $\mathbf{b}\cdot\mathbf{n}$, and we discretize the convective term with a linear upwind:
\begin{align} \label{eq:linearUpwind}
		\sum_{F \in \partial K_{\textrm{f}}} \mathbf{u}_F (\mathbf{b}\cdot\mathbf{n} )|F| & = 
		\sum_{F \in \partial K_{\textrm{f}}}
		a\mathbf{u}_+ + (1-a)\mathbf{u}_- + {\nabla\mathbf{u}}_{\textrm{up}} (\mathbf{x}_F - \mathbf{x}_{\textrm{up}})
\end{align}
where the weight $a$ is defined as
\begin{equation}
	a = 
	\begin{cases}
		1 \text{ if } \mathbf{b}\cdot\mathbf{n} \ge 0, \\
		0 \text{ if } \mathbf{b}\cdot\mathbf{n} < 0,
	\end{cases}
\end{equation}
and the index $\textrm{up}$ is 
\begin{equation}
	\textrm{up} = 
	\begin{cases}
		+ \text{ if } \mathbf{b}\cdot\mathbf{n} \ge 0, \\
		- \text{ if } \mathbf{b}\cdot\mathbf{n} < 0,
	\end{cases}
\end{equation}
where $\nabla \mathbf{u}_\pm$ is computed with a Gauss formula.
}

\newpage

\printbibliography

\end{document}